\documentclass{article}
\usepackage{amsmath,amssymb}
\usepackage[all]{xy}
\newtheorem{thm}{Theorem}[section]
\newtheorem{prop}{Proposition}[section]
\newtheorem{cor}{Corollary}[section]
\newtheorem{lem}{Lemma}[section]
\newtheorem{de}{Definition}[section]
\newtheorem{ex}{Example}[section]
\newtheorem{rem}{Remark}[section]
\newenvironment{proof}{
                        \noindent{Proof.}}
                                       {\hfill {$\mathbf \Box$}\medskip}

\newcommand{\got}[1]{\mathfrak{#1}}

\newcommand{\K}{\mathbb{K}}
\newcommand{\N}{\mathbb{N}}

\newcommand{\C}{\mathbb{C}}

\newcommand{\Q}{\mathbb{Q}}

\newcommand{\HH}{\mathrm{H}}
\newcommand{\Hom}{\mathrm{Hom}}

\newcommand{\Def}{\mathrm{Def}}
\newcommand{\J}{\mathrm{J}}
\newcommand{\I}{\mathrm{I}}

\newcommand{\Der}{\mathrm{Der}}

\newcommand{\D}{\mathrm{D}}

\newcommand{\m}{\mathfrak{m}}
\newcommand{\n}{\mathfrak{n}}
\newcommand{\Gl}{\mathrm{Gl}}
\newcommand{\G}{\mathrm{G}}
\newcommand{\M}{\mathrm{M}}
\newcommand{\Aut}{\mathrm{Aut}}

\newcommand{\gl}{\mathfrak{gl}}
\newcommand{\g}{\mathfrak{g}}
\newcommand{\V}{\mathrm{V}}
\newcommand{\A}{\mathrm{A}}

\renewcommand{\r}{\mathfrak{r}}
\renewcommand{\L}{\got{L}}
\renewcommand{\O}{\mathcal{O}}
\makeatletter
\@addtoreset{equation}{section}

\makeatother
\title{Deformations of Lie algebras and Induction of Schemes}
\author{Roger Carles\\UMR 6086 du CNRS, Laboratoire de Math\'{e}matiques, Universit\'e de Poitiers\\F-8692 Futuroscope Chasseneuil France\\{carles@math.univ-poitiers.fr}
\\ Toukaiddine Petit \footnote{\tt Supported by the EC project Liegrits MCRTN 505078.}
\\Department Wiskunde en Informatica, Universiteit Antwerpen\\
B-2020 (Belgium)\\{toukaiddine.petit@ua.ac.be}.}    
\date{}
\begin{document}
\maketitle
\begin{abstract}
Let $\L_m$ be the scheme of the laws defined by the identities of Jacobi on $\K^m$. The
local studies of an algebraic Lie algebra $\g=\mathrm{R}\ltimes\n$ in $\L_m$ and its
nilpotent part $\n$ in the scheme $\L_n^{\mathrm{R}}$ of $\mathrm{R}$-invariant Lie algebras on
$\K^n$ are
linked. This comparison is made by means of slices, which are transversal subschemes
 to the orbits of $\g$ and $\n$ under the classical groups acting
on $\L_m$
and $\L_n^{\mathrm{R}}$ respectively. We prove a reduction theorem saying that, under
certain conditions on $\g$, the local rings of the slices at $\g$ and $\n$ are
isomorphic. In particular, $\g$ is rigid if and only if is $\n$.
In the formalism developed at beginning of this paper, a deformation of $\g$
with base a local ring $\A$ is a local morphism from the local ring of $\L_m$ at $\g$ to $\A$. So the
study of deformations for a large class of Lie algebras $\g$ in $\L_m$ is equivalent to that of $\n$ in
$\L_n^{\mathrm{R}}$ "modulo" the actions of groups, which is a more simple problem.
The laws of $\L_n^{\mathrm{R}}$ are nilpotent with the choice of $\mathrm{R}$ and then we can construct
these laws by central extensions. This corresponds to an induction on the
schemes themselves $\L_n^{\mathrm{R}}\rightarrow\L_{n+1}^{\mathrm{R}}$. We restrict this study to a torus
$\mathrm{R}=\mathrm{T}$ for certain slices. This leads to a concept of continuous
families with the possibility to have nilpotent parameters $t$ (the schemes are generally
not reduced). This gives an alternative formalism for the problem of
obstructions classes in the theory of formal deformations of M.Gerstenhaber.
Examples are given with $t^2=0$ ($t\neq 0$) and $t^5=0$ ($t^4\neq 0$). 
\end{abstract}
\newpage
%\tableofcontents
\section*{Introduction}
The formal deformations (i.e. with base $\K[[t]]$) of associative algebras and Lie algebras over a closed field $\K$ of characteristic zero have been first studied by M. Gerstenhaber \cite{G}, A. Nijenhuis and R.W. Richardson \cite{NR}. M. Schlessinger generalized the theory of deformations of algebras from the base $\K[[t]]$ to a commutative local $\K$-algebra, \cite{S}. 
A Lie algebra of dimension $m$ is viewed as a point $\phi_0$ of the scheme $\L_m$ defined by the identities of Jacobi. In the prolongation of these ideas, we consider in \cite{C3}, a deformation of a Lie algebra $\phi_0$, parametrized by a local ring $\A$, as a morphism of local rings $\O\rightarrow\A$, where $\O$ is the local ring of $\L_m$ at $\phi_0$. 
The canonical deformation or deformation identity, $\mathrm{id}:\O\rightarrow\O$, which is defined by germs of coordinate functions of the scheme $\L_{m}$ at $\phi_0$, is an initial universal object in the category of the deformations at the point $\phi_0$. Each deformation of $\phi_0$ with base $\A$ may be deduced from the canonical deformation by a base changing such that we can restrict the study of deformation of $\phi_0$ on this one. There is a subgroup $\G_m(\A)$ of the linear group $\Gl_m(\A)$, canonically acting on the set of deformations $\Def(\phi_0,\A)$ and we denote by $\overline{\Def}(\phi_0,\A)$ the set of orbits. In this paper we study the deformations of Lie algebras with base $\A$ a local ring and we develop a theory of sequences of schemes. The deformations of Poisson polynomial algebras and the enveloping algebras of Lie algebras will be separately treated in a forthcoming work, \cite{CP}. We might find the study of deformations of the enveloping algebras of Lie algebras in (\cite{BTM},\cite{P})
. This paper is organized as follows.\\
In Section $1$, we give the classical material of deformations in this new definition where the properties depend on the ring $\A$. If $\A$ is integral we obtain a theory of obstruction which generalizes that of Gerstenhaber, \cite{G}, in the case where $\A$ is equal to $\K[[t]]$. However, if $\A$ is equal to $\O$ the obstructions will disappear with the profit of nilpotent elements; what we see on the examples. We introduce the schemes $\L_n^\mathrm{D}$ consisting of laws which are invariant under $\D\subset\gl_n(\K)$ and $\Def(\varphi_0,-)^{\mathrm{D}}$ the  associated functor of deformations of $\varphi_0\in\L_n^\mathrm{D}(\K)$ which is representable by the local ring $\O_{\varphi_0}^{\mathrm{D}}$ of $\L_n^{\D}$ at $\varphi_0$. We assume now $\mathrm{D}$ is a completely reducible subalgebra in $\gl_n(\K)$. There is a subgroup $\G_n(\A)^{\D}$ consisting of elements of $\G_n(\A)$ which commute with $\mathrm{D}$ and acting on $\Def(\varphi_0,-)^{\mathrm{D}}$. Denote by  $\overline{\Def}(\varphi_0,\A)^{\mathrm{D}}$ the set of orbits under $\G_n(\A)^{\D}$.\\
Finally, we introduce a subscheme $\L_{m,\phi_0}^{\mathcal{A}}$ of $\L_m$ consisting of laws $\phi$ such that $\phi^k_{ij}=(\phi_0)^k_{ij}$ for all $(^k_{ij})\in\mathcal{A}$, where $\mathcal{A}$ is a subset of $\mathcal{I}$. The associated functor of deformations of $\phi_0\in\L_m^\mathcal{A}(\K)$ is representable by the local ring $\O_{\phi_0}^{\mathcal{A}}$ of $\L_{m,\phi_0}^{\mathcal{A}}$ at $\phi_0$. We show that if $n<m$ and $\D$ is completely reducible then there is a set $\mathcal{A}$ such that $\L_{m,\phi_0}^{\mathcal{A}}$ is isomorphic to $\L_n^{\D}$ as schemes. \\
In Section $2$, we suppose that $\A$ is a Noetherian complete local ring. We extend some results shown in the analytical case in \cite{NR} to the local case which some of them were announced in \cite{C3}. In particular we develop a process which permits to fix a maximal family of parameters which are related to the orbit of $\phi_0$, leaving thus free only the parameters said to be essential. This is obtained by a criterion of linear nature which allows to select a family of multi-indices $\mathcal{A}$, called admissible at $\phi_0$, describing these parameters that one must fix. Let $\mathcal{A}$ be such an admissible set at $\phi_0$. Then for each $\phi\in\Def(\phi_0,\A)$ a deformation of $\phi_0$, there exists an element $\Phi$ of the orbit $[\phi]$ under $\G_m(\A)$ such that $\Phi^k_{ij}=(\phi_0)^k_{ij}$ for all $(^k_{ij})\in\mathcal{A}$.\\ 
All these processes can be applied directly to any scheme endowed with an algebraic action of an algebraic group, in particular the schemes defined by associative laws. This method can be programmed on a computer. This leads to the classical concept of versal
deformation, (\cite{K},\cite{S}). One can parametrize such a deformation on $\m\otimes\mathrm{H}^2(\g,\g)$ where $\mathrm{H}^2(\g,\g)$ is identified with the quotient of the tangent to the scheme $\L_m$ at the point $\phi_0$ by that of the orbit at the same point.
For instance, for any admissible set $\mathcal{A}$ at $\phi_0$, then all deformation $h:\O_{\phi_0}\longrightarrow\A$, is equivalent to a deformation $h_0$ which is defined by
$h_0=\overline{h}_0\circ\pi:\O_{\phi_0}\stackrel{\pi}{\longrightarrow}\O^{\mathcal{A}}_{\phi_0}\stackrel{\overline{h}_0}{\longrightarrow}\A,$
where $\pi$ will call the versal deformation of $\phi_0$ relative to $\mathcal{A}$. \\
In Section $3$, we show the Reduction Theorem (Theorem \ref{N31}):
let $(\g,\phi_0)=\mathrm{R}\ltimes\n$ be an algebraic Lie algebra with $(\n,\varphi_0)$ the nilpotent radical, $\mathrm{R}$ a maximal reductive Lie subalgebra and $n$ the dimension of $\n$. We suppose that $\g$ satisfies the following properties\\
1. $\overline{\mathrm{i}}_1:\mathrm{H}^1\left(\n,\n\right)^{\mathrm{R}}\rightarrow\mathrm{H}^1\left(\g,\g\right)$
is an epimorphism,\\
2. $\overline{\mathrm{i}}_2:\mathrm{H}^2\left(\n,\n\right)^{\mathrm{R}}\rightarrow\mathrm{H}^2\left(\g,\g\right)$
is an isomorphism, and \\
3. $\overline{\mathrm{i}}_3:\mathrm{H}^3\left(\n,\n\right)^{\mathrm{R}}\rightarrow\mathrm{H}^3\left(\g,\g\right)$
is a monomorphism. \\
For all Noetherian complete local ring $\A$ and an admissible set $\mathcal{A}$ (resp. $\mathcal{A'}$) at $\phi_0$ (resp. $\varphi_0$) then there is a bijection (resp. a local isomorphism) 
$$\overline{\Def}\left(\varphi_0,\A\right)^\mathrm{R}\simeq\overline{\Def}\left(\phi_0,\A\right)\quad(\mathrm{resp.}\quad\O^{\mathcal{A'},\mathrm{R}}_{n,\varphi_0}\simeq\O^{\mathcal{A}}_{m,\phi_0}).$$
Then the local study of $\phi_0$ in the slice $\L^{\mathcal{A}}_{m,\phi_0}$ is equivalent to that of $\varphi_0$ in the slice $\L^{\mathrm{R},\mathcal{A}'}_{n,\varphi_0}$. 
In Section $4$, we give two new criteria of formal rigidity of Lie algebras. The first criterion (Theorem \ref{t4.1}) characterizes the formal rigidity of $\phi_0$ by the vanishing of the dimension of Krull of $\widehat{\O^{\mathcal{A}}_{\phi_0}}$ for all admissible set $\mathcal{A}$ at $\phi_0$, where $\widehat{\O^{\mathcal{A}}_{\phi_0}}$ is the completion ring of $\O^{\mathcal{A}}_{\phi_0}$ with respect to its maximal ideal. 
The second criterion (Theorem \ref{t6.2}) characterizes the formal rigidity for a class of solvable Lie algebra in the Zariski closed subset $\mathfrak{R}_m$ of $\L_m$ consisting of solvable laws: \\
we suppose that there is a valuation on $\K$. If $\g=\mathrm{T}\ltimes\n$ is a solvable Lie algebra satisfying $\mathrm{H}^1(\g,\g)=\mathrm{H}^0(\g,\g)$. Then the two following properties are equivalent:\\
	i) $\n$ is formal rigid in $\L_n^{\mathrm{T}}$;\\
	ii) we have $[\g,\n]=\n$ and $\g$ is formal rigid in $\mathfrak{R}_m$.\\
This result is the most elaborate form which we could expect for the formal rigidity. 
Under the hypotheses of the Reduction Theorem, the local study of $\phi_0$ in $\L_m$ is equivalent to that of $\varphi_0$ in $\L^\mathrm{R}_n$. If the weights of the center of $\mathrm{R}$ were different from zero, the scheme $\L_n^{\mathrm{R}}$ would be consisted of nilpotent laws.
Consequently, in this case, we would limit the study to $\L^\mathrm{R}_n$. Passing from $\L_m$ to $\L^\mathrm{R}_n$ would have the advantage to use methods which would be specific for nilpotent Lie algebras.  \\
The study of rigid complete solvable Lie algebras in \cite{C8} had led quite naturally to that a sequence of algebraic solvable Lie algebras $\g^{(n)}:=\mathrm{T}\ltimes\varphi_n$ where the passage from $\varphi_n$ to $\varphi_{n+1}$ is made by central extension, the maximal torus $\mathrm{T}$ on $\varphi_n$ being extended on $\varphi_{n+1}$ with an additional weight $\alpha_{n+1}$, and $n$ is the dimension of $\varphi_n$. We denote again by $\mathrm{T}$ this new torus. \\
In Section $5$, one develops a process under the hypotheses of a simple path of weights $(\alpha_1,...,\alpha_n,...)$ related to $\mathrm{T}$, which defines a sequence of schemes $\L_n^{\mathrm{T}}$ for $n\geq n_0$. Under the hypotheses chosen on the weights, the local study of the Lie algebras $\g^{(n)}$ will profit from a double method:\\
1) the theorem of reduction which is applied here and reduces the local study of $\g^{(n)}$ to that of $\varphi_n$ in $\L_n^{\mathrm{T}}$;\\
2) the local study of the laws $\varphi_n$ is done by induction on $n$ starting from an initialization $n_0$ where the torus $\mathrm{T}$ appears.\\
Let $\sum_n\left(\mathrm{T}\right)$ be the subset of $\L_n^{\mathrm{T}}\left(\K\right)$ consisting of laws such that $\mathrm{T}$ is a maximal torus of derivations. If a central extension $\varphi_{n+1}$ of $\varphi_n\in\sum_{n}(\mathrm{T})$ belongs to $\sum_{n+1}(\mathrm{T})$ we will say that $\varphi_{n+1}$ is obtained by a direct filiation of $\varphi_n$. This means that there are two admissible sets $\mathcal{A}_n$ and $\mathcal{A}_{n+1}$ at $\varphi_{n}$ and $\varphi_{n+1}$ respectively such that 
$\mathcal{A}_{n}\subset\mathcal{A}_{n+1}$. Then this inclusion induces a process of construction by extension denoted by $\L_n^{\mathrm{T},\mathcal{A}_n}\longrightarrow\L_{n+1}^{\mathrm{T},\mathcal{A}_{n+1}}$, and
is said to be a direct filiation. This leads to a concept of continuous
families (Theorem \ref{n51}) with the possibility to have nilpotent parameters (the schemes are generally
not reduced). One gives some general theorems for this construction (Theorem \ref{p6.1} and Theorem \ref{t6.1}). More general cases of filiations $\L_n^{\mathrm{R},\mathcal{A}_n}\longrightarrow\L_{n+1}^{\mathrm{R},\mathcal{A}_{n+1}}$ could be studied later on.\\ 
In Section $6$, we give some examples and one sees how, in the filiations, the parameters of the deformation appear, and how they can become nilpotent or disappear when $n$ increases. The calculations of the first example, 6.1, are very simple, which permits to study another problem: that of the research of the nilpotent elements in the global ring of the scheme $\L_n^{\mathrm{T}}$.  The technique employed here, entirely different from the local treatment, is based on the elimination of the variables in the sequence of the polynomial equations given by the identities of Jacobi. The second example, 6.3  is still obtained by simple calculations for its local aspect: a local nilpotent element with order 5 appears in dimension 12. 
It remains to understand this phenomenon of appearance of the nilpotent elements.
\section{Schemes and Deformations}
\subsection{Schemes}
In this work we consider a commutative algebraically closed field $\K$ of characteristic $0$. Let $\left\{e_i,i=1,\cdots,m\right\}$ be the canonical basis of $\V:=\K^m$. Let $\A$ be a commutative associative $\K$-algebra with unity $1=1_\A$, and let $\L_m\left(\A\right)$ denote the set of all Lie $\A$-algebra multiplications on $\A\otimes_\K\V$. An element $\phi$ of $\L_m\left(\A\right)$ is defined by the structure constants $\phi^k_{ij}\in\A$: $\phi\left(e_i,e_j\right)=\sum_{k=1}^{m}\phi^k_{ij}e_k$. The statement that $\phi$ forms a Lie multiplication is equivalent to the following relations for the structure constants:
\begin{equation}\label{n1}
\phi^k_{ij}+\phi^k_{ji}=0,\quad J_{ijk}^l=\sum_{s=1}^{p}\left(\phi^s_{ij}\cdot\phi^l_{sk}+\phi^s_{jk}\cdot\phi^l_{si}+\phi^s_{ki}\cdot \phi^l_{sj}\right)=0.
\end{equation}
Let $\mathrm{B}$ be another commutative associative $\K$-algebra with unity. If $f:\A\rightarrow\mathrm{B}$, is a $\K$-algebra morphism, then it induces a map $\L_m\left(f\right):\L_m\left(\A\right)\rightarrow\L_m\left(\mathrm{B}\right)$, defined by $\left(a^k_{ij}\right)\mapsto\left(f\left(a^k_{ij} \right)\right)$. Thus we obtain a functor of the category of commutative associative $\K$-algebras to the category of sets. We denote by $\J_m$ the ideal of the ring of polynomials $\mathrm{P}_m:=\K\left[X^k_{ij}:1\leq i<j\leq m,1\leq k\leq m\right]$ generated by the relations (\ref{n1}), and $\I_m$ the quotient algebra $\mathrm{P}_m/\J_m$. We can see that, giving a point $\phi$ of $\L_m\left(\A\right)$ (a Lie multiplication on $\A\otimes_\K\V$) is equivalent to giving a $\K$-algebra morphism, $f=f_{\phi}:\I_m\rightarrow\A$. The $\K$-functor $\L_m$ is representable by $\I_m$ and we have $\L_m\simeq\mathrm{Spec}\left(\I_m\right)$, cf  \cite{R}. The $\K$-functor $\L_m$ is an affine algebraic $\K$-scheme with algebra $\I_m$. 
The set of rational points of $\L_m\left(\A\right)$ is identified with $\L_m\left(\K\right)$ the set of all Lie $\K$-algebra multiplications on $\V$. An element $\phi_0$ of $\L_m(\K)$ is defined by a ring $\K$-morphism $f_0:\I_m\rightarrow\K$, by its maximal ideal $\m_0:=\ker(f_0)$ or the injection $\mathrm{Spec}(f_0):\mathrm{Spec}(\K)\rightarrow\L_m$.
\begin{prop}The Zariski tangent space $\mathrm{T}(\got{L}_{m})_{\phi_0}(\K)$ to $\got{L}_{m}$ at $\phi_0$ is equal to the space $\mathrm{Z}^2(\g,\g)$ of two cocycles of $\g$ with values in the adjoint $\g$-module $\g$. 
\end{prop}
\subsection{Deformations} 
Let $\A$ be a local ring with maximal ideal $\m=\m\left(\A\right)$, residue field $\K=\A/\m$, and augmentation map $\mathrm{pr}:\A\rightarrow\K$. Denote by $\widehat{\A}$ its completion for the $\m$-adic topology.
\begin{de} 
A deformation of a point $\phi_0$ of $\L_m\left(\K\right)$ with base $\A$ is a point $\phi$ of $\L_m\left(\A\right)$ such that $\phi_0=\L_m\left(\mathrm{pr}\right)\left(\phi\right)$. We extend $\Def\left(\phi_0,-\right)$ to $\widehat{\A}$ by taking inverse limits.
\end{de}
For all fixed point $\phi_0$, we obtain a local functor, $\Def\left(\phi_0,-\right)$, of the category of commutative associative local rings $\A$ to the category of the sets. If $\O$ is the local ring to the scheme $\L_m$ at $\phi_0$, then it is equal to the localization ring of $\I_m$ by the maximal ideal $\m_0$.\\
Giving a deformation $\phi$ of $\phi_0$ with base $\A$ is equivalent to giving a ring $\K$-morphism, $f:\I_m\rightarrow\A$ such that $f_0=\mathrm{pr}\circ f$.
This is equivalent to $ f$ sends $\m_0$ to $\m$. Then $f$ induces a local morphism from $\O$ to $A$. We deduce that giving a deformation $\phi$ of $\phi_0$, is equivalent to give a local $\K$-ring morphism, 
$$f:\O\rightarrow\A.$$
\begin{thm}The functor $\mathrm{Def}(\phi_0,-)$ is representable by $\O$, i.e. for all local ring $\A$ there exists a bijection,
$$\mathrm{Def}(\phi_0,\A)\rightarrow\Hom_\K\left(\O,\A\right).$$
\end{thm}

\begin{de}The deformation $\mathrm{id}:\O\rightarrow\O$ obtained by the point $\mathrm{X}=\left(X_{ij}^k\right)$ consisting of germs of coordinate functions of the scheme $\L_{m}$ at $\phi_0$, is called deformation identity or canonical deformation. It is an initial universal object in the category of the deformations at the point $\phi_0$.
\end{de}
If $h:\A\rightarrow\mathrm{B}$, is a local ring $\K$-morphism, then it induces a map,
 $\Hom_\K\left(\O,\A\right)\rightarrow\Hom_\K\left(\O,\mathrm{B}\right)$, defined by $f\rightarrow h\circ f$, we will call it, a base changing of deformation. Each deformation of $f_{0}$ with base $\A$ may be deduced from the canonical deformation by a base changing.\\
\begin{rem}A deformation $\phi\in\Def\left(\phi_0,\A\right)$ is defined with the structure constants, $\phi_{ij}^k=f\left(X_{ij}^k\right)$. The smallest subring of $\A$ containing them, is equal to $f\left(\O\right)$. It is a Noetherian local ring since it is isomorphic to 
$\O/\ker\left(f\right)$. Next we can assume that $\A$ is Noetherian.
\end{rem}
\textbf{Notation.} Let $\got{R}$ be the category of Noetherian local $\K$-rings, $\widehat{\got{R}}$ that of Noetherian complete local $\K$-rings.
\subsection{Lie brackets and Chevalley-Eilenberg complex} 
The Chevalley-Eilenberg complex $\mathrm{C}\left(\g,\g\right)$ of a Lie algebra $\g$ defined on an underlying vector space $\V$ endowed with a law $\phi_0$ has the components $\mathrm{C}^m\left(\g,\g\right):=\Hom_\K\left(\wedge^m\g,\g\right)$ in degrees $m\geq 0$ and vanishing components in negative degrees. Its differential is defined by
$$\left(d f\right)\left(x_1,...,x_{m+1}\right)=\sum_{1\leq i< j\leq m+1}(-1)^{i+j}f\left(\phi_0\left(x_i,x_j\right),x_1,...,\widehat{x_i},...,\widehat{x_j},...,x_{m+1}\right)$$
	$$+\sum_{1\leq k\leq m+1}(-1)^{k+1}\phi_0\left(x_k,f\left(x_1,...,\widehat{x_k},...,x_{k+1}\right)\right),$$
where $f\in\mathrm{C}^m\left(\g,\g\right),x_1,..,x_{m+1}\in\g$ and the symbol $\widehat{x}$ indicates that $x$ must be omitted. 
The Richardson-Nijenhuis product of a $m$-cochain $f$ by a $q$-cochain $g$ is the $(m+q-1)$-cochain defined by
$$(f\bullet g)\left(x_1,...,x_{m+q-1}\right)=\sum\mathrm{sign}(\sigma)f\left(g\left(x_{\sigma(1)},...,x_{\sigma(q)}\right),x_{\sigma(q+1),...,x_{\sigma(q+m-1)}}\right),$$
where $\sigma$ runs through the permutations which are increasing on $\left\{1,..,q\right\}$ and $\left\{q+1,...,m+q-1\right\}$. The Nijenhuis-Richardson bracket of a $m$-cochain $f$ by a $q$-cochain $g$ is the commutator
$$\left[f,g\right]=f\bullet g-(-1)^{(m-1)(q-1)}g\bullet f.$$
It may be checked that $d f=(-1)^{m+1}\left[\phi_0,f\right]$.\\ 
We denote by $\mathrm{Z}^m\left(\g,\g\right)$, $\mathrm{B}^m\left(\g,\g\right)$ and $\mathrm{H}^m\left(\g,\g\right)$, the space of $m$-cocycles, $m$-coboundaries, and $m$-cohomologies respectively.
\subsection{Obstructions to extending deformations}
Let $f\in\Hom_\K\left(\O,\A\right)$ be a deformation with $\A$ assumed complete, then, for each $p\geq 1$, its reduction modulo $\m^p$ is a $\A/\m^{p}$-deformation. We obtain a map
\begin{equation}\label{e14}
\Hom_\K\left(\O,\A\right)\quad\longrightarrow\quad\underleftarrow{\mathrm{lim}_p}\Hom_\K\left(\O,\A/\m^p\right)
\end{equation}

\begin{de}
We call a truncated deformation with order $p$ associated with a deformation $f:\O\rightarrow\A$, the deformation $f_p=\pi_p\circ f$ with base $\A/\m^{p+1}$ and $\pi_p:\A\rightarrow\A/\m^{p+1}$, the canonical homomorphism. In particular $f_1$ is called the infinitesimal or tangent deformation. 
\end{de}
It follows from the bijection (\ref{e14}) that for each element $f$ of $\Hom_\K\left(\O,\A\right)$, there is a family $\left(f_p\right)_{p\in\N}$ with $f_p\in\Hom_\K(\O,\A/\m^{p+1})$ such that $f=\underleftarrow{\mathrm{lim}}\left(f_p\right)$ as $p\rightarrow\infty$. The statement that $f$ forms a deformation is equivalent to lift each local morphism $f_p$ to the local morphism $f_{p+1}$. This is equivalent to the following diagram is commutative:
\begin{equation}
	\xymatrix{
	\O
	\ar@{->}[r]^{f_{p+1}} \ar@{->}[dr]^{f_{p}}
	& \A/\m^{p+2}
		\ar@{->}[d]^{\pi}\\
	&\A/\m^{p+1}}
	\end{equation}
where $\pi$ is the projection, see \cite{R} for $\A=\K[[t]]$. We said that there is an obstruction with order $p+1$ if it is impossible to find such a solution.\\
Let $\A\in\widehat{\got{R}}$ be a Noetherian complete local ring. If $r$ is the dimension of $\m/\m^2$ over $\K$, then $\A=\K\left[\left[T_1,...,T_r\right]\right]/\mathfrak{a}$ with $\mathfrak{a}$ an ideal, cf \cite{B}. 
Let $f\in\Hom_\K\left(\O,\A\right)$ be a deformation defined by the structure constants $f\left(X_{ij}^k\right)=\phi_{ij}^k$. It follows that the structure constants $(\phi_{ij}^k)$ may be written as formal series of generators $t_1,...,t_r$ of $\m\left(\A\right)$, say $$\phi_{ij}^k=\phi_{ij}^k\left(t_1,...,t_r\right),$$ 
such that 
$$\L_{m}(\mathrm{pr})\left(\phi_{ij}^k\right)=(\phi_0)^k_{ij}=\phi^k_{ij}\left(0,...,0\right).$$ 
This deformation may be written as 
\begin{equation}\label{e16}
	(\phi_{ij}^k)=\phi\left(t\right)=\sum_{\mu\in\N^r}t^\mu\phi_\mu,
\end{equation}
where $\phi_\mu\in\mathrm{C}^2\left(\g,\g\right)$, $\mu=\left(\mu_1,...,\mu_r\right)\in\N^r$ and $t^\mu=t_1^{\mu_1}...t_r^{\mu_r}$. The decomposition (\ref{e16}) of $(\phi_{ij}^k)$ is not unique since the presence of the ideal $\mathfrak{a}$. It will become when $\mathfrak{a}=0$. The law $(\phi_{ij}^k)$ is a Lie multiplication iff 
\begin{equation}\label{e17}
	[\phi(t),\phi(t)]=\sum_{\gamma\in\N^r} t^\gamma(\sum_{\alpha+\beta=\gamma}[\phi_\alpha,\phi_\beta])=0.
\end{equation}
We then obtain
\begin{equation}\label{E18}
	\sum_{\alpha+\beta=\gamma}\left[\phi_\alpha,\phi_\beta\right]=-2d_{\phi_0}\phi_\gamma+2\omega_\gamma\quad\mathrm{with}\quad \omega_\gamma:=\frac{1}{2}\sum^{\alpha+\beta=\gamma}_{\alpha\neq 0,\beta\neq 0}[\phi_\alpha,\phi_\beta].
\end{equation}
%Since [$[\phi(t),\phi(t)],\phi(t)]=0$, then $[\phi_0,\sum_{\gamma\in(\N^*)^r}t^\gamma\omega_\gamma]=0$, i.e. 
%\begin{equation}
%	\sum_{\gamma\in(\N^*)^r}t^\gamma\omega_\gamma\in\A\otimes\mathrm{Z}^3(\g,\g).
%\end{equation}
It follows from (\ref{E18}) that the (\ref{e17}) is equivalent to the following equation in $\A\otimes\mathrm{C}^3(\g,\g)/\mathrm{B}^3(\g,\g)$
\begin{equation}\label{E19}
\sum_{\gamma\in(\N^*)^r}t^\gamma\overline{\omega}_\gamma=0,
\end{equation}
where $\overline{\omega}_\gamma$ is the class of $\omega_\gamma$ modulo $\mathrm{B}^3(\g,\g)$.
If the sequence of the $\overline{\omega}_{\alpha}$ is nonzero, let $\overline{\omega}_1,...,\overline{\omega}_i,..$ denote the maximal linear independent family on $\K$ which can be extracted from them. We index the finite family $(\omega_i)_i$ on $I=\left\{1,...,q\right\}$ such that the natural order relation on $I$ is compatible with the order relation given by the length $|\gamma_i|$ of multi-indices $\gamma_i$ satisfying $\omega_i=\omega_{\gamma_i}$, i.e, 
$$i\leq j\Leftrightarrow|\gamma_i |\leq |\gamma_j|.$$ 
Each $\overline{\omega}_\gamma$ may be uniquely written as $\overline{\omega}_\gamma=\sum_{i\in I}c^i_\gamma\overline{\omega}_i,$ 
where $c^i_\gamma\in\K$. \\
The Eq. (\ref{E19}) is equivalent to the following system of equations called obstructions:
\begin{equation}
	\sum_{\gamma}c^i_\gamma t^\gamma=0,\quad\mathrm{with}\quad i\in I
\end{equation}
The law $\phi$ forms a Lie multiplication iff the $|I|$-formal series $\sum_{\alpha} c^i_\alpha T^\alpha$ for $i\in I$ successively have to belong to the ideal $\got{a}$, where $|I|$ is the cardinal of $I$. The particular case where the ideal $\mathfrak{a}$ is contained in $T^{q+1}\K\left[\left[T\right]\right]$ then the elements $\overline{\omega}_\gamma$ are zero for all $\gamma$ such that $|\gamma|\leq q$. In the case $\mathfrak{a}=0$ the existence of the deformation $\phi$, i.e., $\phi\in\L_{m}\left(\A\right)$, is then equivalent to the vanishing of the sequence of $\overline{\omega}_\gamma$.
\begin{rem} If $\overline{\omega}_\gamma=0$ for all $|\gamma|<n$, then we have $\overline{\omega}_\gamma\in\mathrm{H}^3(\g,\g)$, for all $\gamma$ such that $|\gamma|=n$, \cite{G}.
\end{rem}
\subsection{Equivalence of Deformations}
Let $\A$ be a local ring and $\Gl_{m}(\A)$ be the group of invertible square matrices with order $m$ with coefficient in $\A$.
The augmentation map $\mathrm{pr}:\A\rightarrow\K$ induces a group $\K$-morphism, ${\mathrm{pr}}:\Gl_{m}\left(A\right)\rightarrow\Gl_{m}\left(\K\right)$.
Let $\G_{m}\left(\A\right)$ denote the kernel of ${\mathrm{pr}}$, then it is the subgroup of $\Gl_{m}\left(\A\right)$ which is equal to $\mathrm{id}+\M_{m}\left(\m\right)$ where $\M_{m}\left(\m\right)$ the set of matrices with coefficients in $\m$. The group $\Gl_{m}\left(\A\right)$ canonically acts on $\L_{m}\left(\A\right)$ by
\begin{equation}
	\ast:\Gl_{m}\left(\A\right)\times\L_{m}\left(\A\right)\rightarrow\L_{m}\left(A\right)
\end{equation}

\begin{equation}\label{A}
	s\ast\phi=s\circ\phi\circ(s^{-1}\times s^{-1})
\end{equation}
where $s\in\Gl_{m}\left(\A\right)$ and $\phi\in\L_{m}\left(\A\right)$. 
The canonical action of $\Gl_{m}\left(\K\right)$ on $\L_{m}\left(\K\right)$ is obtained by projection of the action (\ref{A}). We denote by $\left[\phi\right]$ the $\Gl_{m}\left(\A\right)$-orbit of $\phi$.
\begin{de}Two deformations of $\phi_0$ with same base $\A$ are said to be equivalent if they lie on the same orbit under $\G_{m}\left(\A\right)$. Denote by $\overline{\Def}\left(\phi_0,\A\right)$ the set of equivalence classes of deformations of $\phi_0$ with base $\A$.
\end{de}
\begin{de}A deformation is said to be $\A$-trivial if it is equivalent to the constant deformation $c:\O\rightarrow\A$ defined by $c\left(u\right):=f_{0}\left(u\right)\cdot 1_{\A}$, where $u\in\O$.\\
A point $\phi_0$ of $\L_{m}\left(\K\right)$ is said to be $\A$-rigid if each deformation of $\phi_0$ is $\A$-trivial, i.e., $\Def\left(\phi_0,\A\right)$ consists of only one orbit under $\G_{m}\left(\A\right)$.
\end{de}
The notion of $\A$-rigidity depends on the base $\A$. The geometric rigidity which means that the orbit of $\phi_0$ under $\G_{m}\left(\K\right)$ is a Zariski open in $\L_{m}\left(\K\right)$, corresponds to the rigidity with base $\A:=\K[[t]]$ which is equal to the ring of formal power series, since $\K$ is algebraically closed, \cite{NR}.

\begin{prop}The Zariski tangent space $\mathrm{T}(\left[\phi_0\right])_{\phi_0}(\K)$ to  the $\G_m(\K)$-orbit $[\phi_0]$ of $\phi_0$ at $\phi_0$ is equal to the space $\mathrm{B}^2(\g,\g)$ of two coboundaries of $\g$ with values in the adjoint $\g$-module $\g$.
\end{prop}
We deduce the classical result.
\begin{thm}
The statement $\mathrm{H}^2(\phi_0,\phi_0)=0$ is equivalent to the $\G_m(\K)$-orbit of $\phi_0$ is a Zariski open set in $\L_{m}\left(\K\right)$ and the scheme $\L_{m}$ is reduced at point $\phi_0$.
\end{thm}
We suppose that $\A\in\widehat{\got{R}}$, then each element $s$ of $\G_{m}\left(\A\right)$ may be written as $\sum_{\mu\in\N^r}t^\mu s_\mu$, with $s_0=\mathrm{id}$ and $s_\mu\in\M_{m}\left(\K\right)$. Let $k$ be the smallest nonzero length of $\mu$ such that $\sum_{|\mu|=k}t^\mu s_\mu$, is nonzero. The action of $s$ on a deformation $\phi=\sum_{\mu\in\N^r}t^\mu\phi_\mu$ of $\phi_0$ with base $\A$, is given by
\begin{equation}\label{e4}
	s\ast\phi\equiv\sum_{|\mu|<k}t^\mu\phi_\mu+\sum_{|\mu|=k}t^\mu\left(\phi_\mu-d_{\phi_0} s_\mu\right)\mod\m^{k+1}
\end{equation}
\subsection{Automorphisms and Derivations of Deformations}
Let $\A\in\widehat{\got{R}}$ be a Noetherian complete local ring.
\begin{de}Let $\phi\in\L_{m}(\A)$ be a point, we define $\Aut\left(\phi,\A\right)$ the group of automorphisms of $\phi$ as the set of matrices $s$ of $\Gl_{m}\left(\A\right)$ such that $s\ast\phi=\phi$. 
\end{de}
If $\phi\in\Def\left(\phi_0,\A\right)$ is a deformation and $s$ is an element of $\Aut\left(\phi,\A\right)$ then we have $s_0:=\mathrm{pr}\left(s\right)\in\Aut\left(\phi_0,\K\right)$ since $s_0\ast\phi_0=\phi_0$. 
\begin{de}Let $\phi\in\L_{m}(\A)$ be a point, we define $\Der\left(\phi,\A\right)$ the Lie $\A$-algebra of derivations of $\phi$ as the set of matrices $\delta$ of $\M_{m}\left(\A\right)$ such that $\delta\cdot\phi=0$.
\end{de}
Let $\phi\in\Def\left(\phi_0,\A\right)$ be a deformation. To lift  a derivation $\delta_0\in\Der\left(\phi_0,\K\right)$ of $\phi_0$ to a derivation $\delta\in\Der\left(\phi,\A\right)$  of $\phi$ is a difficult problem which is linked to obstructions. However we have the following properties.
\begin{lem}\label{n6}Let $\phi\in\Def\left(\phi_0,\A\right)$ be a deformation. If $\delta\in\Der(\phi,\A)$ then $\delta$ is a deformation of $\delta_0:=\mathrm{pr}(\delta)\in\Der(\phi_0,\K)$. If $\delta_0$ is an inner derivation of $\phi_0$ then it may be lifted to a derivation of $\phi$. 
\end{lem}
\begin{proof} We can check that $\mathrm{pr}\left(\delta\right)\cdot\phi_0=0$ and $\delta_0=\mathrm{pr}\left(\delta\right)\in\Der\left(\phi_0,\K\right)$. If $x$ belongs $\g$ such that $\delta_0=\phi_0\left(x,.\right)$ then the map $\delta:=\phi\left(x,.\right)$ is a derivation of $\phi$.
\end{proof}
\begin{lem}\label{n7}Let $\phi\in\Def\left(\phi_0,\A\right)$ be a deformation. If $\delta$ is an element of $\m\otimes\Der(\phi,\A)$, then $\exp\left(\delta\right)=\sum^{\infty}_{n=0}\frac{1}{n!}\delta^n$ is an automorphism of $\phi$.
\end{lem}
\begin{proof}
Since $\A$ is complete, hence this formula converges in the sens of Krull in $\M_{m}\left(\A\right)$. We can reason by induction on the integer $n$ and we will have:
\begin{equation*}
	\delta^n\circ\phi=\sum^{n}_{i=0}(^n_i)\phi(\delta^i(-),\delta^{n-i}(-)).
\end{equation*}
\end{proof}
\begin{prop}\label{p1.13}Let $\phi\in\Def\left(\phi_0,\A\right)$ be a deformation. Let $\left\{L_i\right\}$ be a basis of $\M_{m}\left(\K\right)$ such that there are derivations of $\phi_0$, $(L_i)_{1\leq i\leq q}$, which may be lifted to  derivations of $\phi$. For all $\phi'\in\left[\phi\right]$  there exists an element $s$ belonging to the subgroup generated by $\mathrm{id}+\sum_{i>q} \m\otimes L_i$ of $\G_{m}\left(\A\right)$ such that $\phi'=s\ast\phi$.
\end{prop}
\begin{proof} Since each $L_i$ with $1\leq i\leq q$ may be lifted to a derivation of $\phi$, denoted by $L_i\left(t\right)$, then they may be written as, 
$$L_i\left(t\right)=\sum a_{ji}L_j\in\A\otimes_\K\M_{m}\left(\K\right),$$ 
where $a_{ji}$ is equal to $\delta_{ji}+t_{ji}$, $\delta_{ji}$ the Kronecker symbol and $t_{ji}\in\m$. Consider the matrix $c$ defined by $c_{ji}=a_{ji}$ if $i\leq q$ and $c_{ji}=\delta_{ji}$ if $i>q$, hence it is a invertible square matrix. Its inverse matrix $\left(b_{ji}\right)$ satisfies $b_{ji}\in\delta_{ji}+\m$. Hence,
\begin{equation}\label{e7}
	L_i=\sum_{j\leq q} b_{ji}L_j\left(t\right)+\sum_{j>q} b_{ji}L_j,\quad\mathrm{for\,all }\quad i
\end{equation}
From (\ref{e7}) then for any element $\phi'$ of the orbit $\left[\phi\right]$ under $\G_m\left(\A\right)$ there is $s\in\G_m\left(\A\right)$ such that
$$\phi'=s\ast\phi\quad\mathrm{and}\quad s=\mathrm{id}+\delta+\sum_{j>q}\tau^j L_j$$ 
where $\tau_j\in\m$ and $\delta\in\m\otimes\Der\left(\phi,\A\right)$. We will prove by induction on the integer $r\in\N$ the following property:
there exists $s_{\left(r\right)}\in\G_m\left(A\right)$ such that
\begin{equation}
	\phi'=s_{\left(r\right)}\ast\phi\quad\mathrm{and}\quad s_{\left(r\right)}=\mathrm{id}+\delta_{\left(r\right)}+\sum_{j>q}\tau^j_{\left(r\right)}L_j
\end{equation}
where $\delta_{\left(r\right)}\in\m^r\otimes\Der\left(\phi,\A\right)$ and $\tau^j_{\left(r\right)}\in\m$. It is obvious if $r=0$. From Lemma \ref {n7}, the series $\exp(\delta_{\left(r\right)})$ converges to an automorphism of $\phi$. Using the induction hypothesis we set
$$s_{\left(r+1\right)}:=s_{\left(r\right)}\circ\exp\left(-\delta_{\left(r\right)}\right).$$
We obtain that $s_{\left(r+1\right)}=\mathrm{id}+\sum_{j>q}\tau^j_{\left(r\right)}L_j+\mathrm{terms\,of\,degree}>r$. It follows from the equation (\ref{e7}) that 
$$s_{\left(r+1\right)}=\mathrm{id}+\delta_{\left(r+1\right)}+\sum_{j>q}\tau^j_{\left(r+1\right)}L_j,$$ 
$\delta_{\left(r+1\right)}\in\m^{r+1}\otimes\Der\left(\phi,\A\right)$ and  $\tau^j_{\left(r+1\right)}\equiv\tau^j_{\left(r\right)}\mod\m^{r+1}$. The sequence $\left(\tau^j_{\left(r\right)}\right)_{j>q}$ being convergent in the sense of Krull by construction, thus $$\phi'=\lim_{r\rightarrow\infty}s_{\left(r\right)}\ast\phi=\left(\lim_{r\rightarrow\infty}s_{\left(r\right)}\right)\ast\phi=s_\infty\ast\phi$$ with $s_\infty=\mathrm{id}+\sum_{j>q}\tau^j_{\left(\infty\right)}L_j\in\mathrm{id}+\sum_{i>q} \m\otimes L_i$.
\end{proof}
\subsection{$\D$-schemes and $\D$-deformations}
Let $\D$ be a subset of $\mathrm{M}_m(\K)$, and let $\L_m^{\D}\left(\K\right)$ denote the set of laws $\phi$ of $\L_m\left(\K\right)$ satisfying, $\delta\cdot\phi=0$, for all $\delta=(\delta^i_j)\in\D$. They are all Lie multiplications on $\V$ which admit $\D$ as derivations.
Let $\Delta_m$ be the ideal of $\K\left[X^k_{ij}\right]$ generated by the polynomials
\begin{equation}\label{E1.14}
\sum^p_{l=1}\left(X^l_{ij}\delta^k_l-\delta^l_iX^k_{lj}-\delta^l_jX^k_{il}\right),\left(i,j,k\in\left\{1,...,m\right\},i<j\right)
\end{equation}
for all $\delta\in\D$. Denote by $\J^{\D}_m$ the sum of $\J_m$ and $\Delta_m$, by $\I_m^{\D}$ the quotient algebra $\K\left[X^k_{ij}\right]/\J_m^{\D}$. 
Let $\L_m^{\D}$ denote the subscheme of $\L_m$ which is canonically isomorphic to $\mathrm{Spec}\left(\I_m^{\D}\right)$, then it is a functor from the category of commutative $\K$-algebras to the category of sets. Let $\phi_0$ be a point of $\L_m^{\D}(\K)$, $\A$ a local ring and denote by $\Def^{\D}(\phi_0,A)$ the set of deformations of $\phi_0$ with base $\A$ in $\L_m^{\D}(\A)$. Then we obtain again a local functor $\Def_{\D}(\phi_0,-)$, from the category of commutative associative local rings to the category of the sets and is representable by the local ring $\O_{\phi_0}^\D$ of $\L_m^{\D}$ at $\phi_0$ , i.e. for all local ring $\A$ there is a bijection from $\Def_{\D}(\phi_0,\A)$ onto $\Hom_\K(\O_{\phi_0}^{\D},\A)$.
The ring quotient $\K$-morphism from $\I_m$ to $\I_{m}^{\D}$, induces a local ring $\K$-morphism $\eta:\O_{\phi}\rightarrow\O_{\phi_0}^{\D}$ and a scheme embedding from $\L_{m}^{\D}$ onto $\L_m$, denoted by $\mathfrak{I}$.
It follows that it induces also a map denoted again by $\mathfrak{I}$ and defined by
$$\Hom_\K(\O_{\phi_0}^\D,\A)\longrightarrow\Hom_\K(\O_{\phi_0},\A),\quad f\longmapsto \mathfrak{I}(f)=f\circ\eta.$$
\begin{prop}
Let $\phi_0\in\L_m\left(\K\right)$ be a point and let $\D$ be a subspace of $\Der\left(\phi_0,\K\right)$. Then the tangent space at $\phi_0$ of $\L_{m}^{\D}\left(\K\right)$ is equal to $\mathrm{Z}^2\left(\g,\g\right)^\D$.
\end{prop}
Let $\left[\phi_0\right]$ be an orbit of a law $\phi_0$ of $\L_m\left(\K\right)$ under $\Gl_m\left(\K\right)$ and denote by $\left[\phi_0\right]^{\D}=\left[\phi_0\right]\cap\L_{m}^{\D}\left(\K\right)$ the orbit trace of $\left[\phi_0\right]$ in $\L_{m}^{\D}\left(\K\right)$.\\  
Assume that $\D$ is a completely reducible Lie subalgebra of $\mathrm{C}^1(\V,\V)$, i.e., its natural representation on $\V$ is semisimple. The normalizer subgroup $\HH$ of $\D$ in $\Gl\left(\V\right)$ stabilizes $\L_m^\D\left(\K\right)$ as well as the trace orbits. Since the adjoint representation of $\D$ is semisimple, hence the Lie algebra of $\mathrm{H}$ is given by $\mathrm{Lie}\left(\mathrm{H}\right)=\D+\mathrm{gl}(\V)^{\mathrm{D}}$. Denote by $\mathrm{H}_0$ the identity component of $\mathrm{H}$.
\begin{prop} \label{p} Let $\D$ be a completely reducible Lie subalgebra of $\mathrm{C}^1(\V,\V)$. Then
\begin{enumerate}
	\item $[\phi_0]^{\D}$ is a finite union of orbits under $\mathrm{H}_0$.
	\item The tangent space to $[\phi_0]^\D$ at $\phi_0$, is equal to the tangent space to its orbit under $\mathrm{H}_0$ at $\phi_0$. It is also equal to $\mathrm{B}^2\left(\g,\g\right)^{\D}$.
	\item The orbits under $\mathrm{H}_0$ are the same that under the identity component $\Gl\left(\V\right)^{\D}_0$ of the group  $\Gl\left(\V\right)^{\D}$ which is formed of the elements of $\Gl\left(\V\right)$ commuting with elements of $\D$.
\end{enumerate}
\end{prop}
\begin{proof} Classical proof, cf \cite{C7}.
\end{proof}\\
We deduce the analogical classical result:
\begin{thm}
For $\K=\C$, the statement $\mathrm{H}^2\left(\phi_0,\phi_0\right)^{\D}=0$ is equivalent to the fact that $\mathrm{H}_0$-orbit of $\phi_0$ is Zariski open set in $\L_{m}^{\D}\left(\K\right)$ and the scheme $\L_{m}^{\D}$ is reduced at $\phi_0$.
\end{thm}
Each deformation of $\phi_0$ in $\L^{\D}_m$ may be written as $\sum_{\alpha}t^{\alpha}\phi_{\alpha}$, where $\phi_{\alpha}\in\mathrm{C}^2(\V,\V)^{\D}$ and the analogical elements $\overline{\omega}_\gamma$ defined as in Eq. (\ref{E19}), belong to\\ $\mathrm{C}^3(\g,\g)^{\D}/\mathrm{B}^3(\g,\g)^\D=(\mathrm{C}^3(\g,\g)/\mathrm{B}^3(\g,\g))^\D$. If $\overline{\omega}_\gamma$ is zero for all $|\gamma|<n$, then $\overline{\omega}_\gamma$ belongs to $\mathrm{H}^3(\g,\g)^\D$, for all $|\gamma|=n$. In particular the condition $\mathrm{H}^3(\g,\g)^\D=0$, implies that all obstructions are zero.
\subsection{Fixing of coordinates in $\L_m$}
We will give a method of construction of some subschemes of $\L_m$ by the following way. It consists to fix some coordinates which permits to give nice properties. \\
Let $\mathcal{I}=\left\{(^k_{ij}):1\leq i<j\leq,1\leq k\leq m\right\}$ be the set of all multi-indices. Let $\mathcal{A}$ be a subset of $\mathcal{I}$ and let $\phi_0$ be a point of $\L_m\left(\K\right)$. Set 
\begin{equation}\label{E1.15}
	\J^\mathcal{A}_{m,\phi_0}:=\J_m+\left\langle  X^{\alpha}-(\phi_0)^{\alpha}:\alpha\in\mathcal{A}\right\rangle,\ \I^\mathcal{A}_{m,\phi_0}:=\K \left[X^{\alpha}:\alpha\in\mathcal{I}\right]/\J^\mathcal{A}_{m,\phi_0},
\end{equation}
we define a subscheme $\L^\mathcal{A}_{m,\phi_0}$ of $\L_m$ as 
\begin{equation}
\L^\mathcal{A}_{m,\phi_0}(\A):=\left\{\phi\in\L_m(\A)|\phi^\alpha=(\phi_0)^\alpha,\forall\alpha\in\mathcal{A}\right\}.
\end{equation}
Then
\begin{equation}
\L^\mathcal{A}_{m,\phi_0}\simeq\mathrm{Spec}\left(\I^{\mathcal{A}}_{m,\phi_0}\right).
\end{equation}
$\mathrm{\textbf{Semi-direct product}}$. We suppose now that $\mathcal{A}$ is the set of $\left(^k_{ij}\right)\in\mathcal{I}$ such that $i,j$ or $k>n$ with $n\leq m$. We assume that the law $\phi_0$ defines a semi-direct product $\g$ of a reductive part $\mathrm{R}=\left\langle e_{n+1},..,e_m\right\rangle$ by the nilpotent part $\mathfrak{n}=\left\langle e_1,...,e_n\right\rangle$. The adjoint action of $\mathrm{R}$ on $\n$ defines derivations $\delta_i$ of the restriction $\varphi_0$ of $\phi_0$ on $\n$ by
\begin{equation}\label{E1.18}
	\delta_i e_j=\phi_0\left(e_{n+i},e_j\right)=\sum_{k=1}^nC^k_{n+i,j} e_k,\quad\mathrm{with}\quad i\leq r\quad\mathrm{and}\quad j\leq n\quad\mathrm{if}\quad r=m-n.
\end{equation}
The elements $(\delta_i)_{i\leq r}$ generate the reductive Lie algebra 
$\mathrm{D}=\mathrm{ad}_{\n}\mathrm{R}$.
\begin{prop}\label{m12} Let $	\L^\mathcal{A}_{m,\phi_0}$ be the subscheme of $\L_m$ defined by $\mathcal{A}$ and $\phi_0$ as above. Then there is a scheme bijection, $\L^\mathcal{A}_{m,\phi_0}\simeq\L^\mathrm{R}_{n}.$
\end{prop}
\begin{proof} Set $\mathcal{B}:=\mathcal{I}-\mathcal{A}$. The algebra $\I^\mathcal{A}_{m,\phi_0}$ of definition of the scheme $\L^\mathcal{A}_{m,\phi_0}$ is isomorphic to $\K[X^\alpha:\alpha\in\mathcal{B}]/\J^\mathcal{B}_{m,\phi_0}$ where $\J^\mathcal{B}_{m,\phi_0}$ is the ideal defined by 
$$\J_m/\left\langle  X^{\alpha}-(\phi_0)^{\alpha}:\alpha\in\mathcal{A}\right\rangle\cap\J_m.$$
It suffices to set $X^\alpha=(\phi_0)^\alpha$ for all $\alpha\in\mathcal{A}$ in the Jacobi polynomials $J^d_{abc}$. For $a,b,c\leq n$ that gives the Jacobi polynomials in function of coordinates $(X^\alpha)_{\alpha\in\mathcal{B}}$ since $\n$ remains a Lie subalgebra. If $a=n+i$, $b$ and $c$ $\leq n$, then the terms $\delta_i\varphi(e_b,e_c)-\varphi(e_b,\delta_i e_c)-\varphi(\delta_i e_b,e_c)$ generate the ideal $\Delta_m$ defined by (\ref{E1.14}). The other Jacobi polynomials give null constants. We obtain $$\I^\mathcal{A}_{m,\phi_0}\cong\K[X^\alpha:\alpha\in\mathcal{B}]/(\J_m\cap\left\langle X^\alpha, \alpha\in\mathcal{B}\right\rangle+\Delta_m)=\I_{n,\varphi_0}^{\mathrm{R}}.$$
\end{proof}\\
\textbf{Notation}. This applies to a Lie algebra $(\g,\phi_0)=\mathrm{R}\ltimes\n$ semi-direct product of a reductive Lie algebra $\mathrm{R}$ by a nilpotent ideal $(\n,\varphi_0)$. We will denote simply by $\L_n^{\mathrm{R}}$ instead of $\L_n^\mathrm{D}$ for $\mathrm{D}=\mathrm{ad}_\n\mathrm{R}$.
The $\K$-epimorphism algebra from $\I_m$ to $ \I^\mathrm{R}_{n}$ induces a scheme embedding from $\L^\mathrm{R}_{n}$ to $\L_m$, denoted by $\mathfrak{I}$. Let $\A$ be a ring, thus the morphism 
$$\mathfrak{I}:\L^\mathrm{R}_{n}(\A)\longrightarrow\L_m(\A),$$
maps $(\varphi^k_{ij})$ to a law $(\phi^k_{ij})$, defined by
$$\phi^k_{ij}=\varphi^k_{ij}\quad\mathrm{if}\quad i,j,k\leq n \quad\mathrm{and}\quad\phi^k_{ij}=\left(\phi_0\right)_{ij}^k\quad\mathrm{if}\quad i,j\quad\mathrm{or}\quad k\quad\mathrm{is}>n.$$
If $\A=\K$, the morphism $\mathfrak{I}$ sends a Lie algebra $\n'$ to the semi-direct product $\mathrm{R}\ltimes\n'$. 
The Lie algebra $\n$ chosen at beginning does not play a privileged role and $\mathrm{R}$ remains fixed. The differential map $\mathrm{i}_2$ of $\mathfrak{I}$ at $\varphi_0$ is an injection of Zariski tangent spaces, i.e. $\mathrm{i}_2:\mathrm{Z}(\n,\n)^\mathrm{R}\rightarrow\mathrm{Z}^2(\g,\g)$ which maps $\mathrm{B}^2(\n,\n)^{\mathrm{R}}$ into $\mathrm{B}^2(\g,\g)$ and induces a map on the quotient, $\overline{\mathrm{i}}_2:\mathrm{H}(\n,\n)^\mathrm{R}\rightarrow\mathrm{H}^2(\g,\g)$, which play an important role in the section $4$. Denote by $\O^{\mathrm{R}}_{\varphi_0}$ the local ring of $	\L^\mathrm{R}_{n}$ at $\varphi_0$, then the morphism $\mathfrak{J}$ induces an epimorphism local ring $\eta$ at $\phi_0$,
\begin{equation}\label{e1.18}
	\eta:\O_{\phi_0}\rightarrow\O_{\varphi_0}^\mathrm{R}.
\end{equation}
We denote again by $\mathfrak{I}$ the map defined by
$$\Hom_\K(\O_{\phi_0}^\mathrm{R},\A)\longrightarrow\Hom_\K(\O_{\phi_0},\A),\quad f\longmapsto \mathfrak{I}(f)=f\circ\eta.$$
\section{Versal Deformation in $\L_m$}
In this section, we extend some results shown in the analytical case in \cite{NR} to the local case which some of them were announced in \cite{C3}. However we interpret them in term of parameters. We give also a method to built a versal deformation.
\subsection{Parameters of deformations}
We keep the above notations. Let $\A$ be an element of $\widehat{\got{R}}$. Let $(e_i)_{1\leq i\leq m}$ be a canonical basis of $\V$. The elements $e^{ij}_k$ of $\mathrm{C}^2\left(\V,\V\right)$ defined by  $e^{ij}_k\left(e_{i'},e_{j'}\right):=\delta^i_{i'}\delta^j_{j'}e_k$ form a basis of $\mathrm{C}^2\left(\V,\V\right)$ with $\delta^i_j$ the Kronecker symbol. A law $\phi\in\L_m\left(\A\right)$ is defined by the structure constants $\phi^k_{ij}\in\A$: $\phi=\left(\phi^k_{ij}\right)=\sum_{i<j,k}\phi^k_{ij}e^{ij}_k$ and satisfying Jacobi's identities. More generally we can take a basis $(e_{\alpha})_{\alpha\in\mathcal{I}}$ of $\mathrm{C}^2\left(\V,\V\right)$, an element $\phi=\sum_{\alpha\in\mathcal{I}}\phi^\alpha e_{\alpha}$ of $\mathrm{C}^2\left(\V,\V\right)$ is a Lie multiplication iff $\left[\phi,\phi\right]=2\phi\circ\phi=0$. We suppose now that $\phi$ is a deformation of a law $\phi_0$ of $\L_m\left(\K\right)$, hence the components of $\phi$ may be written as $\phi^{\alpha}=\phi_0^{\alpha}+\xi^{\alpha}$, with $\alpha\in\mathcal{I}$ and $\xi^\alpha\in\m$. Let $\overline{\m}$ denote the quotient of $\m$ by $\m^2$.
Each deformation $f:\O\rightarrow\A$ is defined in fact on the subring $f\left(\O\right)$ of $\A$ which is the quotient of $\O$ by $\ker\left(f\right)$.
\begin{de} One calls a family of parameters of a deformation $f$ of $\phi_0$, a set $\left\{t_i:1\leq i\leq r\right\}$ of $\m\left(f\left(\O\right)\right)$ such that the set $\left\{\overline{t_i}:1\leq i\leq r\right\}$ forms a basis of $\overline{\m}\left(f\left(\O\right)\right)$ on $\K$, where $\overline{t_i}$ is the class of $t_i$ in $\overline{\m}\left(f\left(\O\right)\right)$. The dimension of $\overline{\m}\left(f\left(\O\right)\right)$ over $\K$ is said to be the number of parameters of $f$.
\end{de}
The number of parameters of the canonical deformation is equal to the dimension of the Zariski tangent vector space of scheme $\L_m$ at $\phi_0$.
\begin{lem} The number of parameters of a deformation of $\phi_0$ is majored by the $\K$-dimension of $\mathrm{Z}^2\left(\g,\g\right)$. We have the equality if the kernel of the deformation $f:\O\rightarrow\A$ is contained in $\m^2\left(\O\right)$.
\end{lem}
\begin{proof} It is sufficed to remark that the local $\K$-morphism $f$ induces a surjective $\K$-linear morphism from $\overline{\m}\left(\O\right)$ to $\overline{\m}\left(f\left(\O\right)\right)$.
\end{proof}
\subsection{Parametrization of deformations}
Let $\A\in\widehat{\got{R}}$ be a Noetherian complete local ring with maximal ideal $\m=\m\left(\A\right)$ and residue field $\K=\A/\m$, and augmentation map $\mathrm{pr}:\A\rightarrow\K$. Let $\mathrm{H}^k$ be a complement of the $k$-coboundaries $\mathrm{B}^{k}(\g,\g)$ in the $k$-cocycles $\mathrm{Z}^{k}(\g,\g)$ and  $\mathrm{W}^k$ a complement of $\mathrm{Z}^k(\g,\g)$ in $\mathrm{C}^k(\g,\g)$.  It is easy to check that $\mathrm{H}^k$ is canonically isomorphic to $\mathrm{H}^k(\g,\g)$ and $d$ induces an isomorphism from $\mathrm{W}^k$ onto $\mathrm{B}^{k+1}(\g,\g)$. Thus we obtain
\begin{equation}\label{e3.1}
	\mathrm{C}^k(\g,\g)=\mathrm{Z}^k(\g,\g)\oplus\mathrm{W}^k=\mathrm{B}^k(\g,\g)\oplus\mathrm{H}^k\oplus\mathrm{W}^k.
\end{equation}
This decomposition is called a decomposition of Hodge associated with $\g$. We can check that $\phi=\phi_0+\xi$ is a deformation of $\phi_0$ if and only if $\xi$ belongs to the set of solutions of Maurer-Cartan's equation
\begin{equation}\label{e10}
\mathrm{MC}_m\left(\m\right)=\left\{\xi\in\m\otimes\mathrm{C}^2\left(\g,\g\right):d\xi-\frac{1}{2}\left[\xi,\xi\right]=0\right\},
\end{equation}
where the differential $d$ extends to $\A\otimes\mathrm{C}\left(\g,\g\right)$ by $\mathrm{id}_{\A}\otimes d$ and denote it again by $d$. We have
\begin{equation}\label{e3.3}
d(\xi)-\frac{1}{2}\left[\xi,\xi\right]=d(\xi)-\frac{1}{2}\left[\xi,\xi\right]_{\mathrm{B}^3}-\frac{1}{2}\left[\xi,\xi\right]_{\mathrm{H}^3}-\frac{1}{2}\left[\xi,\xi\right]_{\mathrm{W}^3},
\end{equation}
where $\left[\xi,\xi\right]_{\mathrm{B}^3},\left[\xi,\xi\right]_{\mathrm{H}^3},\left[\xi,\xi\right]_{\mathrm{W}^3}$ denote the projections on $\m\otimes\mathrm{B}^3(\g,\g)$, $\m\otimes\mathrm{H}^3$ and $\m\otimes\mathrm{W}^3$ respectively.
Set 
$$\mathrm{MC}'_m(\m)=\left\{\xi\in\m\otimes\mathrm{C}^2\left(\g,\g\right):d(\xi)-\frac{1}{2}\left[\xi,\xi\right]_{\mathrm{B}^3}=0\right\},$$
we can see $\mathrm{MC}_m(\m)\subset\mathrm{MC}'_m(\m)$. 
Using the equation (\ref{e3.1}), we may choose a basis  $(e^k_\alpha)_{\alpha\in\mathcal{I}^k}$ of $\mathrm{C}^k(\g,\g)$ such that the elements $e^k_\alpha$ are indexed by $\mathcal{H}^k$, $\mathcal{B}^k$ and $\mathcal{W}^k$, and contained in $\mathcal{I}^k$, and form a basis of $\mathrm{H}^k$, $\mathrm{B}^k(\g,\g)$ and $\mathrm{W}^k$ respectively. Denote by $|\mathcal{A}|$ the cardinal of a set $\mathcal{A}$, we have $|\mathcal{W}^k|=|\mathcal{B}^{k+1}|$. Set $p:=|\mathcal{W}^2|=|\mathcal{B}^3|.$

\begin{lem}\label{l2.1}There is an unique map  $g:\m\otimes\mathrm{Z}^2(\g,\g)\rightarrow\m\otimes\mathrm{W}^2(\g,\g)$ with $g(0)=0$ such that $\mathrm{MC}_m'(\m)=\left\{\xi=z+g(z):z\in\m\otimes\mathrm{Z}^2(\g,\g)\right\}$.
\end{lem}
\begin{proof} Consider a map $F:\m\otimes\mathrm{Z}^2(\g,\g)\times\m\otimes\mathrm{W}^2\rightarrow\m\otimes\mathrm{B}^3(\g,\g)$ defined by $F(z,w)=d(w)-\frac{1}{2}\left[z+w,z+w\right]_{\mathrm{B}^3}$.
Denote by  $(F^\alpha)_{\alpha\in\mathcal{W}^2}$, $(z^\alpha)_{\alpha\in\mathcal{Z}^2}$ and $(w^\alpha)_{\alpha\in\mathcal{W}^2}$ the components of $F$, $z$ and $w$ relative to above bases. Since $F$ is polynomial in variables $z^\alpha$ and $w^\alpha$, we may regard  $F$ as an element of the ring $\A[[Z,W]]^{p}$ of formal power series in variables $Z=(Z^\alpha)_{\alpha\in\mathcal{Z}^2}$ and $W=(W^\alpha)_{\alpha\in\mathcal{W}^2}$. We have $F(0)=0$. The Jacobian $\left[\partial F^\alpha/\partial{W^\beta}(0)\right]$ is invertible in $\A$, since $\mathrm{pr}(\left[\partial F^\alpha/\partial {W^\beta}(0)\right])$ is invertible in $\A/\m=\K$ (i.e. $d_{\phi_0}:\mathrm{W}^2\rightarrow\mathrm{B}^3$ is an isomorphism).
It follows from above and  the formal implicit theorem \cite{B}, the formal equation $F(Z,W)=0$ admits a solution if and only if there is an unique formal map $G=(G^\alpha)_{\alpha\in\mathcal{W}^2}\in\A[[Z]]^{p}$ such that $G(0)=0$ and $G(Z)=W$. It is clear that $G$ converges in the $\m$-adic sense. Since $\A$ is complete, it follows that we might define by substituting, a map $g^\alpha:\A^{|\mathcal{Z}^2|}\rightarrow\A^{p}$ by $g^\alpha(a_1,...,a_{|\mathcal{Z}^2|}):=G^\alpha(a_1,...,a_{|\mathcal{Z}^2|})$ which satisfies $g^\alpha(\m^{|\mathcal{Z}^2|})\subset\m^{p}$, for all $\alpha$.
\end{proof}\\
It follows from Lemma \ref{l2.1} that we can define a map
$$\Omega:\m\otimes\mathrm{Z}^2(\g,\g)\rightarrow\m\otimes\mathrm{H}^3,\quad z\mapsto\left[z+g(z),z+g(z)\right]_{\mathrm{H}^3},$$
which is called obstruction map.
\begin{thm}\label{t2.1} Then there are two maps $g$ and $\Omega$ defined as above
satisfying $g(0)=\Omega(0)=0$ such that 
$$\mathrm{MC}_m(\m)=\left\{z+g(z):\Omega(z)=0,z\in\m\otimes\mathrm{Z}^2(\g,\g)\right\}.$$
\end{thm}
\begin{proof} $\xi$ is a solution of $\mathrm{MC}_m(\m)$ iff $$d(\xi)-\frac{1}{2}\left[\xi,\xi\right]_{\mathrm{B}^3}=0,\quad\left[\xi,\xi\right]_{\mathrm{W}^3}=0,\quad\left[\xi,\xi\right]_{\mathrm{H}^3}=0.$$ 
We will show that if $\xi$ is such that $\left[\xi,\xi\right]_{\mathrm{H}^3}=0=d(\xi)-\frac{1}{2}\left[\xi,\xi\right]_{\mathrm{B}}$ then $\left[\xi,\xi\right]_{\mathrm{W}^3}=0$. We suppose that $$\left[\xi,\xi\right]_{\mathrm{H}^3}=0=d(\xi)-\frac{1}{2}\left[\xi,\xi\right]_{\mathrm{B}^3},$$
since (graded Jacobi's identity)
$$0=\left[\phi,\left[\phi,\phi\right]\right]=\left[\phi,d(\xi)-\frac{1}{2}\left[\xi,\xi\right]\right]$$ 
and by Eq. (\ref{e3.3}), it follows that $\left[\phi,\left[\xi,\xi\right]_{\mathrm{W}^3}\right]=0$. We can check that $$\left[\phi,\left[\xi,\xi\right]_{\mathrm{W}^3}\right]=(d+ad\xi)(\left[\xi,\xi\right]_{\mathrm{W}^3}),$$ 
and the injectivity of $d:\mathrm{W}^2\rightarrow\mathrm{C}^3(\g,\g)$ implies $ad \phi:\mathrm{W}^2\rightarrow\m\otimes\mathrm{C}^3(\g,\g)$ does. It follows that $\left[\xi,\xi\right]_{\mathrm{W}^3}=0$.
By Lemma \ref{l2.1}, $\xi$ is equal to $z+g(z)$ is equivalent to $d(\xi)-\frac{1}{2}\left[\xi,\xi\right]_{\mathrm{B}^3}=0$ and we deduce the result.  
\end{proof}
\begin{cor}\label{c3.1}Then there are two maps $g$ and $\Omega$ defined as above
satisfying $g(0)=\Omega(0)=0$ such that 
$$\Def(\phi_0,\A)=\left\{\phi_0+z+g(z):\Omega(z)=0,z\in\m\otimes\mathrm{Z}^2(\g,\g)\right\}.$$
\end{cor}
This decomposition of Corollary \ref{c3.1} is not canonical since it depends on the decomposition (\ref{e3.1}).
\subsection{Versal Deformation}
Assume now that $\A\in\widehat{\got{R}}$ with maximal ideal $\m=\m\left(\A\right)$, residue field $\K=\A/\m$ and augmentation map $\mathrm{pr}:\A\rightarrow\K$. A subscheme $\got{F}$ of $\L_m$ is said to be transverse at the orbit $\left[\phi_0\right]$ under the canonical action of complete linear group at a point $\phi_0\in\got{F}$ if the tangent space at $\phi_0$ satisfies
$$T_{\phi_0}\L_m=T_{\phi_0}\left[\phi_0\right]\oplus T_{\phi_0}\left(\got{F}\right),$$
where the orbit is endowed with the structure of reduced scheme. \\
We will construct $\got{F}$ satisfying these properties such that each deformation of $\phi_0$ in $\L_m$ is equivalent to a deformation of $\phi_0$ in $\got{F}$. Such a deformation is null on some parameters which describe the orbit $\left[\phi_0\right]$. These parameters are said to be a family of orbital parameters at $\phi_0$. The canonical deformation of $\got{F}$ is said to be a versal deformation of $\phi_0$ for the scheme $\L_m$. These parameters can be taken in $\m\otimes\mathrm{H}^2\left(\g,\g\right)$.
Let $\theta^\alpha$ denote the $\alpha$-th component of the map $\theta:\Gl\left(\V\right)\rightarrow\mathrm{C}^2\left(\V,\V\right)$ defined by $\theta\left(s\right)=s\ast\phi_0$, on the basis $(e_{\alpha})_{\alpha\in\mathcal{I}}$ of $\mathrm{C}^2\left(\V,\V\right)$. Then its differential mapping at unit $1$ is given by
\begin{equation}
	(D_1\theta^\alpha)\left(L\right)=\left(L\cdot\phi_0\right)^\alpha=-(dL)^\alpha,\quad L\in\mathrm{C}^1\left(\V,\V\right),
\end{equation}
since we have,\\ $\theta^\alpha\left(1+tL\right)=\left(\left(1+tL\right)\ast\phi_0\right)^\alpha=\phi_0+t\left(L\circ\phi_0-\phi_0\left(L,-\right)-\phi_0(-,L\right))^\alpha\mod t^2$.
%\begin{de} A set $\mathcal{A}$ is said to be admissible at point $\phi_0$ if it is a subset of $\mathcal{I}$ such that it indexes a basis of $\mathrm{B}^2\left(\g,\g\right)$.
%\end{de}
\begin{de} A set $\mathcal{A}$ is said to be admissible at point $\phi_0$ if it is a minimal subset of $\mathcal{I}$ such that the rank of $(D_1\theta^\alpha)_{\alpha\in\mathcal{A}}$ is maximal, i.e., $|\mathcal{A}|=\mathrm{rank}\left(d_{\phi_0}\right)=\dim\mathrm{B}^2\left(\g,\g\right)=\dim\left[\phi_0\right]$.
\end{de}
\begin{rem}\label{R2.1}Let $(e_\alpha)_{\alpha\in\mathcal{I}}$ be a basis of $\mathrm{C}^2\left(\g,\g\right)$. From the exchange basis theorem, there are parts \,$\mathcal{B}:=\mathcal{I}-\mathcal{A}\subset\mathcal{I}$ such that the $(e_\beta)_{\beta\in\mathcal{B}}$ complete a given basis of $\mathrm{B}^2(\g,\g)$ to the basis $(e_\alpha)_{\alpha\in\mathcal{I}}$. The admissible sets $\mathcal{A}$ are the complements of $\mathcal{B}$ in $\mathcal{I}$.
\end{rem}
The parameters $\left(\xi^\alpha\right)_{\alpha\in\mathcal{A}}$ are called orbital parameters at point $\phi_0$. The set $\mathcal{A}$ permits to define a subscheme of $\L_m$ (cf section $1.8$) called a slice associated with $\mathcal{A}$ by
\begin{equation}
	\got{F}:=\L^\mathcal{A}_{m,\phi_0}
\end{equation}
where the components $\xi^{\alpha}$, $\alpha\in\mathcal{A}$, are expressed as $\K$-linear combinations of $X^\alpha-(\phi_0)^{\alpha}$.
We check that the tangent space $T_{\phi_0}\left(\got{F}\right)$ is isomorphic to $\mathrm{H}^2(\g,\g)$.\\
Now we consider the set of deformations of $\phi_0$ with base $\A$ in $\got{F}$, i.e.,
\begin{equation}
	\Def_\mathcal{A}(\phi_0,\A)=\left\{\phi\in\L^\mathcal{A}_{m,\phi_0}(\A):\mathrm{pr}(\phi)=\phi_0\right\}.
\end{equation}
\begin{prop}\label{p3.1}For all $\A$ there is a bijection,
$$\Def_\mathcal{A}(\phi_0,\A)\rightarrow\mathrm{Hom}_\K(\O_{\phi_0}^\mathcal{A},\A),$$ 
where $\O^\mathcal{A}_{\phi_0}$ is the local ring at $\phi_0$ of the subscheme $\L^{\mathcal{A}}_{m,\phi_0}$.
\end{prop} 
Let $\mathcal{A}$ an admissible set of $\mathcal{I}$ at point $\phi_0$ and then  the vector subspace $\mathrm{V}^2_{\mathcal{A}}(\g,\g)$ of $\mathrm{C}^2(\g,\g)$ generated by $(e_\alpha)_{\alpha\in\mathcal{I}-\mathcal{A}}$ 
is a linear complement  of $\mathrm{B}^2(\g,\g)$ in $\mathrm{C}^2(\g,\g)$, by Remark \ref{r2.1}.
The set $\Def_{\mathcal{A}}(\phi_0,\A)$ is given by 
\begin{equation}
	\Def_\mathcal{A}(\phi_0,\A)=\left\{\phi_0+v;d v-\frac{1}{2}\left[v,v\right]=0,v\in\m\otimes\V^2_{\mathcal{A}}(\g,\g)\right\}.
\end{equation}
We can choose a Hodge decomposition (\ref{e3.1}) associated with $\phi_0$ such that 
 \begin{equation}
	\mathrm{C}^2(\g,\g)=\mathrm{V}^2_{\mathcal{A}}(\g,\g)\oplus\mathrm{B}^k(\g,\g),\quad\mathrm{and}\quad\mathrm{V}^2_{\mathcal{A}}(\g,\g)= \mathrm{H}^2\oplus\mathrm{W}^2. 
\end{equation} 
Such a decomposition is called a Hodge decomposition at $\phi_0$ compatible with $\mathcal{A}$
\begin{prop}\label{p2.2}For any admissible set $\mathcal{A}$ at $\phi_0$ there is a Hodge decomposition at $\phi_0$ compatible with $\mathcal{A}$ and two maps $g$ and $\Omega$ defined as above
satisfying $g(0)=\Omega(0)=0$ such that 
$$\Def_{\mathcal{A}}(\phi_0,\A)=\left\{\phi_0+h+g(h);\Omega(h)=0,h\in\m\otimes\mathrm{H}^2\right\}.$$
\end{prop}
\begin{proof} See Corollary \ref{c3.1}
\end{proof}\\
%Let $(L_i)_{i\in\mathcal{I}^1}$ be a basis of $\mathrm{C}^1(\g,\g)$ which completes a basis $(L_i)_{i\in\mathcal{Z}^1}$ of $\mathrm{Z}^1(\g,\g)$ and  $\mathrm{W}^1$ be the vector subspace of $\mathrm{C}^1(\g,\g)$ generated by $(L_i)_{i\in\mathcal{W}^1=\mathcal{I}^1-\mathcal{Z}^1}$. 
%Denote by $\G(\m\otimes\mathrm{W}^1)$ the subgroup of $\G_m(\A)$ generated by $\mathrm{id}+\m\otimes\mathrm{W}^1$.
\begin{lem}\label{t3.2}If $\mathcal{A}$ is an admissible subset of $\mathcal{I}$ at $\phi_0$, then the map
$$f:(w,v)\mapsto (\mathrm{id}+w)\ast(\phi_0+v)-\phi_0$$ 
is a bijection from $ (\m\otimes\mathrm{W}^1)\times(\m\otimes\mathrm{V}^2_{\mathcal{A}}(\g,\g))$ to $\m\otimes\mathrm{C}^2(\g,\g)$.
\end{lem}
\begin{proof} 
Let $(f^\alpha)_{\alpha\in\mathcal{I}}$, $(v^\alpha)_{\alpha\in\mathcal{I}-\mathcal{A}}$ and $(w^\alpha)_{\alpha\in\mathcal{W}^1}$ be the components of $f$, $v$ and $w$ relative to above bases (see Section 2.2). We may extend $f$ to an element $F$ of $\A[[V,W]]^{|\mathcal{I}|}$ of formal power series in indeterminates $V=(V^\alpha)_{\alpha\in\mathcal{I}-\mathcal{A}}$ and $W=(W^\alpha)_{\alpha\in\mathcal{W}^1}$. We have $F(0)=0$. The Jacobian $\mathrm{J}:=\left[(F^\alpha/\partial {V^{\beta}},F^{\alpha'}/\partial{W^{\beta'}})(0)\right]$ is given by
\begin{center}
$\mathrm{J}$=\begin{math}\bordermatrix{&\cr        
  & \mathrm{I}_{\mathcal{I}-\mathcal{A}}  & 0  \cr          
  & 0  & d_{{\phi_0}}  \cr  }
  \end{math}
\end{center}
where $\mathrm{I}_{\mathcal{I}-\mathcal{A}}$ is the identy matrix with order $|\mathcal{I}-\mathcal{A}|$.
It is easy to check that $\mathrm{J}$ is invertible in $\A$, since its projection $\mathrm{pr}(\mathrm{J})$ is in $\A/\m=\K$ (i.e $d_{\phi_0}$ is an isomorphism from $\mathrm{W}^1$ to $\mathrm{B}^2(\g,\g)$). It follows from above and the formal inversion theorem \cite{B} that the formal map $F$ is bijective. It is clear that $F$ converges in the $\m$-adic sense. Since $\A$ is complete, it follows that $f$ concides with the associated map $F$ by substituting the indeterminates $V_\alpha$ and $W_\alpha$ by elements of $\A$.
\end{proof}\\

\begin{thm}If $\mathcal{A}$ is an admissible subset of $\mathcal{I}$ at $\phi_0$, then the map
$$F:(w,\phi)\mapsto (\mathrm{id}+w)\ast\phi$$ 
is a bijection from $(\m\otimes\mathrm{W}^1)\times\Def_{\mathcal{A}}(\phi_0,\A)$ to $\Def(\phi_0,\A)$.
\end{thm}
\begin{proof} Let $\phi'=\phi_0+\eta'$ be an element of $\Def(\phi_0,\A))$ with $\eta'\in\m\otimes\mathrm{C}^2(\g,\g)$. By Lemma \ref{t3.2} there are unique $w\in\m\otimes\mathrm{W}^1$ and $v\in\m\otimes\mathrm{V}^2_{\mathcal{A}}(\g,\g)$ such that $\eta'=s\ast(\phi_0+v)-\phi_0$ i.e. $\phi'=s\ast(\phi_0+v)$, with $s=\mathrm{id}+w$. Hence $\phi:=\phi_0+v=s^{-1}\ast\phi'$ and $\phi\in\Def(\phi_0,\A)\cap(\phi_0+\m\otimes\mathrm{V}^2_{\mathcal{A}}(\g,\g))=\Def_{\mathcal{A}}(\phi_0,\A).$
\end{proof}
\begin{cor}\label{c3.2}Let $\mathcal{A}$ be an admissible set of $\mathcal{I}$ at $\phi_0$. 
For all deformation $\phi'\in\Def(\phi_0,\A)$ of $\phi_0$, there exists a $\phi$ in the orbit $[\phi']$ under $\G_m(\A)$ such that $\phi^\alpha=(\phi_0)^\alpha$ for all $\alpha\in\mathcal{A}$. 
\end{cor}
\begin{rem}If $\mathrm{H}^1(\g,\g)=0$, we can show that the element $\phi$ defined in Corollary \ref{c3.2} is unique by using Proposition \ref{p1.13}.
\end{rem}
The local ring $\O_{\phi_0}^\mathcal{A}$ of $\got{F}$ at $\phi_0$ can be obtained by the quotient of the local ring $\O_{\phi_0}$ of $\L_m$ at $\phi_0$ by the ideal $\sum_{\alpha\in\mathcal{A}}\O\cdot\xi^\alpha$, generated by the elements $\xi^\alpha\in\m\left(\O\right)$ where $\alpha\in\mathcal{A}$. The morphism $\mathrm{id}:\O^\mathcal{A}_{\phi_0}\rightarrow\O^\mathcal{A}_{\phi_0}$ which defines the canonical deformation of $\L^{\mathcal{A}}_{m,\phi_0}$ at $\phi_0$ by germs of coordinate functions $Y=(Y^\alpha)_{\alpha\in\mathcal{I}-\mathcal{A}}$, can be also obtained by projection of the canonical deformation of $\phi_0$ in $\L_m$ consisting of germs coordinate functions $X=(X^\alpha)_{\alpha\in\mathcal{I}}$,
\begin{equation}\label{e2.8}
	\xymatrix{
	\O_{\phi_0} 
	\ar@{->}[r]^{\mathrm{id}}
	\ar@{->}[d]_{\pi}
	& \O_{\phi_0}
		\ar@{->}[d]^{\pi}
	\\
	\O^\mathcal{A}_{\phi_0}\ar@{->}[r]_{\mathrm{id}}&	\O^\mathcal{A}_{\phi_0}}
	\end{equation}
We use the notation of Proposition \ref{p2.2}.	
\begin{cor}\label{N21}For any admissible subset $\mathcal{A}$ of $\mathcal{I}$ at $\phi_0$, then 
\begin{enumerate}
	\item all deformation $h:\O_{\phi_0}\longrightarrow\A$, with $\A\in\got{R}$ is equivalent to a deformation $h_0$ which is defined by:
$$h_0:\O_{\phi_0}\stackrel{\pi}{\longrightarrow}\O^{\mathcal{A}}_{\phi_0}\stackrel{\overline{h}_0}{\longrightarrow}\A$$
$\pi$ is said to be the versal deformation of $\phi_0$ relative to $\mathcal{A}$.
\item The canonical deformation $\mathrm{id}:\O^\mathcal{A}_{\phi_0}\rightarrow\O^\mathcal{A}_{\phi_0}$ of $\phi_0$ in $\L^{\mathcal{A}}_{m,\phi_0}$, defined by (\ref{e2.8}), may be written on the completion ring $\widehat{\O^\mathcal{A}_{\phi_0}}$ by
$$\mathrm{id}(\phi_0+Y)=\phi_0+(Y_\alpha)+g(Y_\alpha)\quad\mathrm{with}\quad\Omega((Y_\alpha))=0,$$
where $\alpha$ runs through the set of indices $\mathcal{H}^2$.
The parameters $(Y_\alpha)_{\alpha\in\mathcal{H}^2}$ are said to be essential.
\end{enumerate}
\end{cor}
\begin{proof} 1. By Corollary \ref{c3.2}, hence each deformation $h:\O_{\phi_0}\rightarrow\A$ is equivalent to a deformation $h_0=\phi_0+\eta$ where $h_0\left(X^\alpha\right)=\eta^\alpha=0$ for all $\alpha\in\mathcal{A}$. $h_0$ vanishes on the ideal generated by the elements $X^\alpha$, $\alpha\in\mathcal{A}$, consequently, it factorizes through $\O^\mathcal{A}_{\phi_0}$. Then there is $\overline{h}_0$ such that $h$ is equivalent to  $h_0=\overline{h}_0\circ\pi.$ The statement 2. holds by Proposition \ref{p2.2}
\end{proof}\\
%In the formalism of local rings, $\G_m(\A)$ acts on $\mathrm{Hom}_\K(\O_{\phi_0},\A)$ and, by Corallary \ref{N21} gives the
%\begin{cor}For any admissible subset $\mathcal{A}$ of $\mathcal{I}$ at $\phi_0$. There is one to one correspondind between $\mathrm{Hom}_\K(\O_{\phi_0},\A)/\G_m(\A)$ and $\mathrm{Hom}_\K(\O^{\mathcal{A}}_{\phi_0},\A)$.
%\end{cor}
\begin{rem}\label{r2.1}
\begin{enumerate}
	\item It is clear that the results in this section remain true in $\L_n^{\mathrm{R}}$ by replacing $\G_m(\A)$, $\phi_0$, $\O_{\phi_0}$, $\O^{\mathcal{A}}_{\phi_0}$ and $\mathrm{C}^k(\g,\g)$ to $\G_n(\A)^{\mathrm{R}}$, $\varphi_0$, $\O^\mathrm{R}_{\varphi_0}$, $\O^{\mathrm{R},\mathcal{A'}}_{\varphi_0}$, $\mathrm{C}^k(\n,\n)^{\mathrm{R}}$ respectively. We denote by $\Omega'$, $g'$, $F'$, $\mathcal{I'}$ and $\mathcal{A'}$ the analog correspondings.
	\item One finds in \cite{C9} the local study in each point of the scheme $\L_3$ given by this method. 
\end{enumerate}
 
\end{rem}
\section{Reduction Theorem}
\subsection{Reduction Theorem}
Let $(\g,\phi_0)=\mathrm{R}\ltimes(\n,\varphi_0)$ be an algebraic Lie algebra with $\n$ the nilpotent radical and $\mathrm{R}=\mathrm{U}\oplus\mathrm{S}$ a maximal reductive Lie subalgebra and $n$ the dimension of $\n$. The comparison of local studies of $\varphi_0$ and $\phi_0$ must be done on the essential parameters (those of a versal deformation), i.e. '' modulo'' the actions of groups. This means that we need to give some conditions for which those local studies are equivalent. The best favorable case is to compare the deformation classes with base the local ring  \footnote{\tt The local morphisms $\O_{\phi_0}^{\mathcal{A}}\rightarrow\K[\overline{t}]$, identify with the vetors of the Zariski tangent which is equal to $\mathrm{H}^2(\g,\g)$.} $\K[[t]]/(t^2)=\K+\K\overline{t}$; it is equivalent to compare the two cohomology spaces via the morphism $\overline{\mathrm{i}}_2$. Here the identification of the local problems is equivalent to the bijectivity of $\overline{\mathrm{i}}_2$, it is that one calls the \texttt{minimal condition}, preamble any equivalence of classes of deformations with base a local ring $\A$. Indeed, the reduction with base $\A$ implies the reduction with base its different quotients and in particular  with base the ring $\K[\overline{t}]$. The minimal condition is a structural property which characterizes a class of algebras.\\
The graded Lie algebra morphism $\mathrm{i}:\mathrm{C}\left(\n,\n\right)^\mathrm{R}\rightarrow\mathrm{C}\left(\g,\g\right)$ is defined by the composition of graded differential complex morphisms $\alpha$ and $\beta$ defined by
$$0\rightarrow\mathrm{C}\left(\n,\n\right)^\mathrm{R}\stackrel{\alpha}{\rightarrow}\mathrm{C}\left(\n,\g\right)^\mathrm{R}\stackrel{\pi}{\rightarrow}\mathrm{C}\left(\n,\g/\n\right)^\mathrm{R}\rightarrow 0,$$
and 
$$\mathrm{C}\left(\n,\g\right)^\mathrm{R}\stackrel{\beta}{\rightarrow}\mathrm{C}\left(\g,\g\right).$$
It is defined  by $$\mathrm{i}\left(f\right)|_{\mathrm{R}\times\g^{q-1}}=0\quad\mathrm{and}\quad\mathrm{i}\left(f\right)|_{\n^{q}}=f,$$
for each element $f\in\mathrm{C}^q\left(\n,\n\right)^\mathrm{R}$.
The reader can verify the following lemma:
\begin{lem}\label{n31} For all elements $f$ and $g$ of $\mathrm{C}\left(\n,\n\right)^\mathrm{R}$, we then have
\begin{enumerate}
	\item $i\left(\left[f,g\right]\right)=\left[i\left(f\right),i\left(g\right)\right]$
	\item $i\left(\left[\varphi_0,f\right]\right)=\left[\phi_0,i\left(f\right)\right]$.
\end{enumerate}
\end{lem}
Note that the property $(2.)$ means that the map $\mathrm{i}$ commutes with the differentials $d_{\phi_0}$ and $d_{\varphi_0}$. We deduce the linear morphism of cohomologies: $$\overline{\mathrm{i}}=\oplus_q\overline{\mathrm{i}_q}:\mathrm{H}\left(\n,\n\right)^\mathrm{R}\rightarrow\mathrm{H}\left(\g,\g\right).$$
The groups $\Gl\left(\n\right)_0^\mathrm{R}$ and $\Gl\left(\g\right)$ canonically act on $\mathrm{C}\left(\n,\n\right)^\mathrm{R}$ and $\mathrm{C}\left(\g,\g\right)$  respectively by:
$$\left(s\ast f\right)\left(x_1,...,x_q\right)=\left(s\circ f\right)\left(s^{-1}x_1,...,s^{-1}x_q\right)$$
with $f$ a $q$-cochain and $s$ an element of the group. We define an injection, again denoted by $\mathfrak{I}$ from the group $\Gl\left(\n\right)_0^\mathrm{R}$ into $\Gl\left(\g\right)$ which sends $s\in\Gl\left(\n\right)_0^\mathrm{R}$ into $\mathfrak{I}(s)$ with $\mathfrak{I}(s)|_\n=s$ and $\mathfrak{I}(s)|_\mathrm{R}=\mathrm{id}_{\mathrm{R}}$. We get:
\begin{lem}\label{n32}For all $s\in\Gl\left(\n\right)_0^\mathrm{R}$ and $f\in\mathrm{C}\left(\n,\n\right)^\mathrm{R}$. Let $\varphi_0\in\L_n^\mathrm{R}\left(\K\right)$ and set $\phi_0=\mathfrak{I}(\varphi_0)$, then
\begin{enumerate}
	\item $i\left(s\ast f\right)=\mathfrak{I}(s)\ast i\left(f\right)$,
	\item the following diagram is commutative
	\begin{equation}
	\xymatrix{
	\varphi_0+f 
	\ar@{->}[r]^{\mathfrak{I}}
	\ar@{->}[d]_{s\ast}
	& \phi_0+\mathrm{i}(f)
		\ar@{->}[d]^{\mathfrak{I}(s)\ast}
	\\
	s\ast(\varphi_0+f)\ar@{->}[r]_{\mathfrak{I}}&	\mathfrak{I}(s)\ast(\phi_0+\mathrm{i}(f))}
	\end{equation}
\end{enumerate}
\end{lem}
\begin{proof} The equality $\mathfrak{I}(s)\ast\phi_0=\mathfrak{I}(s\ast\varphi_0)$ comes from the commutativity of $\mathfrak{I}(s)$ with $\mathrm{R}$. The statement 2. is deduced from 1.
\end{proof}\\
The map $\mathrm{i}$ extends  to a graded Lie algebra morphism by   
$$ \mathrm{id}\otimes\mathrm{i}:\A\otimes\mathrm{C}\left(\n,\n\right)^\mathrm{R}\rightarrow\A\otimes\mathrm{C}\left(\g,\g\right)$$
satisfying the similar Lemma \ref{n31} and Lemma \ref{n32}, we denote again it by $\mathrm{i}$. From similar Lemma \ref{n31}, then if $\varphi=\varphi_0+\xi\in\Def_n^\mathrm{R}(\varphi_0,\A)$, i.e. 
$$d\xi-\frac{1}{2}\left[\xi,\xi\right]=0,$$ 
then $$\mathrm{i}(d\xi-\frac{1}{2}\left[\xi,\xi\right])=d\mathrm{i}(\xi)-\frac{1}{2}\left[\mathrm{i}(\xi),\mathrm{i}(\xi)\right]=0,$$ 
i.e., $\phi_0+\mathrm{i}(\xi)\in\Def_m(\phi_0,\A)$.\\
From similar Lemma \ref{n32}, the scheme morphism $\mathfrak{I}$ maps $\L_n^\mathrm{R}$ to $\L_m$ and induces a map denoted again by $\mathfrak{I}$ and  defined by 
$$\mathfrak{I}:\Def_n^\mathrm{R}(\varphi_0,\A)\rightarrow\Def_m(\phi_0,\A)$$
$$\varphi_0+\xi\mapsto\phi_0+\mathrm{i}(\xi)$$
passes to the quotient
$$\overline{\mathfrak{I}}:\overline{\Def}_n^\mathrm{R}\left(\varphi_0,\A\right)\longrightarrow\overline{\Def}_m\left(\phi_0,\A\right)$$
modulo the actions of groups
$$\G_n\left(\A\right)^\mathrm{R}=\mathrm{id}\oplus(\m\otimes\mathrm{C}^1(\n,\n)^{\mathrm{R}})\quad\mathrm{and}\quad \G_m\left(\A\right)=\mathrm{id}\oplus(\m\otimes\mathrm{C}^1(\g,\g)),$$ 
respectively.
\begin{lem}\label{n33}We then have
\begin{enumerate} 
	\item $\mathrm{Z}\left(\g,\g\right)\cap 	\mathrm{i}(\mathrm{C}\left(\n,\n\right)^\mathrm{R})=\mathrm{i}(\mathrm{Z}\left(\n,\n\right)^\mathrm{R})$.
	\item $\mathrm{B}\left(\g,\g\right)\cap \mathrm{i}(\mathrm{C}^{p+1}\left(\n,\n\right)^\mathrm{R})=\mathrm{i}(\mathrm{B}^{p+1}(\n,\n)^{\mathrm{R}}+\delta\mathrm{Z}^p(\n,\g/\n)^{\mathrm{R}})$, 
\end{enumerate}
where $\delta$ defines the connection $\partial$ in the long sequence of the cohomology .
\end{lem}
\begin{proof} 1. The statement $d f=0$ for all $f\in\mathrm{C}\left(\n,\n\right)^\mathrm{R}$ is equivalent to $\mathrm{i}\left(d f\right)=d i\left(f\right)=0$, i.e., $\mathrm{i}\left(f\right)\in\mathrm{Z}\left(\g,\g\right)$ since $i$ is injective.\\
2. Take $d h\in\mathrm{B}^{p+1}\left(\g,\g\right)$ of the form $\mathrm{i}\left(f\right)$ where $f\in\mathrm{C}^{p+1}\left(\n,\n\right)^\mathrm{R}$. The $\mathrm{R}$-invariance implies
 $$\mathrm{ad}_\g\mathrm{R}\cdot d h=d(\mathrm{ad}_\g\mathrm{R}\cdot h)=0,$$ 
thus $\mathrm{ad}_\g\mathrm{R}\cdot h\subset\mathrm{Z}^{p}\left(\g,\g\right)$. It follows that $h$ may be written as $h_0+h_1$ with $h_0\in\mathrm{C}^{p}\left(\g,\g\right)^\mathrm{R}$ and $h_1\in\mathrm{Z}^{p}\left(\g,\g\right)$, and $d h_0=\mathrm{i}\left(f\right)$. Denote by $\rho\left(x\right)$ the operator defined by 
\begin{equation}
	\left(\rho\left(x\right)h_1\right)\left(x_1,...,x_p\right):=h_1\left(x,x_1...,x_p\right),
\end{equation}
where $x,x_1,...,x_p\in\g$ and $\theta$ is the adjoint representation of $\g$ in $\mathrm{C}\left(\g,\g\right)$. It is well-known that $\theta$ and $\rho$ satisfy the following formula 
\begin{equation}\label{e31}
	d\circ\rho\left(x\right)+\rho\left(x\right)\circ d=\theta\left(x\right),\quad\forall x\in\g.
\end{equation}
Applying the formula (\ref{e31}) for $h_0$ and $x\in\mathrm{R}$, hence we obtain $d\rho\left(x\right)h_0=0$. We define a map $\widetilde{h}_0\in\mathrm{C}^{p}\left(\g,\g\right)^\mathrm{R}$ by 
$$\widetilde{h}_0|_{\n^p}=0\quad\mathrm{and}\quad\rho\left(x\right)\widetilde{h}_0=\rho\left(x\right)h_0,$$ 
for all $x\in\mathrm{R}$. Hence $\widetilde{h}_0$ is a cocycle, indeed, we have $d\widetilde{h}_0=0$ on $\n^{p+1}$ by construction and for all $x\in\mathrm{R}$, we have 
$$\rho\left(x\right)(d\widetilde{h}_0)=-d\rho\left(x\right)\widetilde{h}_0+\theta\left(x\right)\widetilde{h}_0=-d(\rho\left(x\right) \widetilde{h}_0)=0.$$
It is clear that $\rho(x)(h_0-\widetilde{h}_0)(x)=0$ for all $x\in\mathrm{R}$ and then 
\begin{equation}\label{e3.4}
	h_0-\widetilde{h}_0=\mathrm{i}(k)+l,
\end{equation}
with $k\in\mathrm{C}^p(\n,\n)^{\mathrm{R}}$, $l=\mathrm{pr}\circ(	h_0-\widetilde{h}_0)$ and $\mathrm{pr}:\g\rightarrow\mathrm{R}$ the projection. It follows from (\ref{e3.4}) that
\begin{equation}
	d h=d(h_0-\widetilde{h}_0)=\mathrm{i}(d k)+d l,
\end{equation}
since $\mathrm{i}\circ d=d\circ\mathrm{i}$. We deduce that $d h$ belongs to the image of $\mathrm{i}$ if and only if is $d l$ i.e. 
\begin{equation}\label{e3.6}
	\mathrm{pr}\circ d l=0.
\end{equation}
Since $l$ is $\mathrm{R}$-invariant (see (\ref{e3.4})), it is easy to prove that $\rho(x) l=0$ for all $x\in\mathrm{R}$. Hence
\begin{equation}
	\rho(x) d l=-d \rho(x) l+\theta(x)d l=0,
\end{equation}
for all $x\in\mathrm{R}$. For all $(x_1,...,x_{p+1})\in\n^{p+1}$, we have
\begin{equation}\label{e3.8}
	(\mathrm{pr}\circ d l)(x_1,...,x_{p+1})=\sum_{i<j}(-1)^{i+j}l([x_i,x_j],...,\widehat{x}_i,...,\widehat{x}_j,..x_{p+1}),
\end{equation}
which is the expression of the differential of the restriction $l|_{\n^p}$ on $\n^p$ in the complex $\mathrm{C}(\n,\mathrm{R})^{\mathrm{R}}$ with value in the trivial $\n$-module $\mathrm{R}$ which is identified with the adjoint action $\n$ in $\g/\n$. It follows from (\ref{e3.6}) and (\ref{e3.8}) that 
\begin{equation}
	l|_{\n^p}\in\mathrm{Z}^p(\n,\g/\n)^{\mathrm{R}}
\end{equation}
such that $(	d l)|_{\n^{p+1}}=\delta (l|_{\n^p})\in\mathrm{Z}^{p+1}(\n,\n)^{\mathrm{R}}$, where $\delta$ defines the connection $\partial$ in the long sequence of cohomology. Then we have
\begin{equation}
	d h=\mathrm{i}(d k+\delta(l|_{\n^p}))\in\mathrm{i}(\mathrm{B}^{p+1}(\n,\n)^{\mathrm{R}}+\delta\mathrm{Z}^p(\n,\g/\n)^{\mathrm{R}}).
\end{equation}
\end{proof}\\
\begin{cor}\label{C3.1}
\begin{enumerate}
	\item $\overline{\mathrm{i}}_k$ is a monomorphism iff $$\mathrm{B}^k(\g,\g)\cap\mathrm{i}_k(\mathrm{C}^k(\n,\n)^{\mathrm{R}})=\mathrm{i}_k(\mathrm{B}^k(\n,\n)^{\mathrm{R}}).$$
	\item  $\overline{\mathrm{i}}_k$ is an epimorphism iff 
	$$\mathrm{Z}^k(\g,\g)=\mathrm{i}_k(\mathrm{Z}^k(\n,\n)^{\mathrm{R}})+\mathrm{B}^k(\g,\g).$$
\end{enumerate}
\end{cor}
We identify the vector space $\mathrm{C}^2(\n,\n)$ with its image by $\mathrm{i}_k$. We choose $\mathcal{I}$ (res. $\mathcal{I}'$) a set which indexes a basis of the vector space $\mathrm{C}^2(\g,\g)$ (resp. $\mathrm{C}^2(\n,\n)^{\mathrm{R}}$) and $\mathcal{A}$ (resp. $\mathcal{A}'$) an admissible set of $\mathcal{I}$ (resp. $\mathcal{I}'$ at $\phi_0$ (resp. $\varphi_0$) such that $\mathcal{I}'\subset\mathcal{I}$.\\
\textbf{The Local Epimorphism $\overline{\eta}$}:\\
If $\overline{\mathrm{i}}_2$ is injective, by Corollary \ref{C3.1}, we then have the equality, $$\mathrm{B}^2(\g,\g)\cap\mathrm{i}_2(\mathrm{C}^2(\n,\n)^{\mathrm{R}})=\mathrm{i}_2(\mathrm{B}^2(\n,\n)^{\mathrm{R}}).$$ 
Let $\mathrm{E}$ be the linear complement of $\mathrm{C}^2(\n,\n)^{\mathrm{R}}$ in $\mathrm{C}^2(\g,\g)$ generated by the $(e^k_{ij})$, with $i>n$, or $j>n$, or $k>n$, and the $\delta.e^k_{ij}$, with $\delta\in\mathrm{ad}\mathrm{R}|_\n$, and $i,j,k\leq n$, cf (\ref{E1.18}). We then have 
$$\mathrm{B}^2(\g,\g)+\mathrm{E}=\mathrm{i}_2(\mathrm{B}^2(\n,\n)^{\mathrm{R}})+\mathrm{E}.$$
If we choose a basis $(e_\alpha)_{\alpha\in\mathcal{I'}}$ of $\mathrm{C}^2(\n,\n)^{\mathrm{R}}$ which is completed by a basis $B$ of $\mathrm{C}^2(\g,\g)$ indexed by $\mathcal{I}$, containing all the $(e^k_{ij})$, with $i>n$, or $j>n$, or $k>n$ and some of $\delta.e^k_{ij}$, with $i$, $j$ and $k$ $\leq n$, then each admissible set $\mathcal{A}'\subset\mathcal{I}'$ at $\varphi_0$ will be completed by an admissible set $\mathcal{A}\subset\mathcal{I}$ at $\phi_0$, which will be contained in $\mathcal{A}'\cup(\mathcal{I}-\mathcal{I}')$. The scheme $\L_{n,\varphi_0}^{\mathrm{R},\mathcal{A}'}$ is the spectrum of the quotient of $\mathrm{P}_m=\K\left[X^k_{ij}:1\leq i<j\leq m,1\leq k\leq m\right]$ by the ideal 
$$\J_{\mathcal{A}'}:=\J_m+\Delta_m+\left\langle X^k_{ij};i>n,\mathrm{or}\,j>n, \mathrm{or}\,k>n\right\rangle+\left\langle \xi^{\alpha}(X)-(\varphi_0)^\alpha;\alpha\in\mathcal{A}'\right\rangle,$$
where $\Delta_m$ is defined by (\ref{E1.14}), and the $\xi^{\alpha}(X)$ are $\K$-linear combination of $X^k_{ij}$ associated with the changing basis $(e^k_{ij})\rightarrow B$. The scheme $\L_{m,\phi_0}^{\mathcal{A}}$ is the spectrum of the quotient of $\mathrm{P}_m$ by the ideal
$$\J_{\mathcal{A}}:=\J_m+\left\langle \xi^{\alpha}(X)-(\phi_0)^\alpha; \alpha\in\mathcal{A}\right\rangle.$$
The identity map with the inclusion $\J_{\mathcal{A}}\hookrightarrow\J_{\mathcal{A}'}$, induces an epimorphism on the quotient algebras from $\mathrm{P}_m/\J_{\mathcal{A}}$ to $\mathrm{P}_m/\J_{\mathcal{A}'}$, an embedding scheme from $\L_{n,\varphi_0}^{\mathrm{R},\mathcal{A}'}$ to $\L_{m,\phi_0}^{\mathcal{A}}$ and finally, the local epimorphism  $\overline{\eta}:\O^{\mathcal{A}}_{\phi_0}\rightarrow\O_{\varphi_0}^{\mathrm{R},\mathcal{A}'}$.
 %Using the notation and the statements of Theorem \ref{c2.2} and Remark \ref{r2.1}, we have the following diagram:
% \begin{equation}
%\xymatrix{
%\G(\m(\widehat{\O_{\phi_0}})\otimes\mathrm{W}^1)\times\widehat{O^{\mathcal{A}}_{\phi_0}} \ar@{->}[r]& \widehat{\O_{\phi_0}} \ar@{->}[d]^{\eta}& \\
%\G(\m(\widehat{\O^{\mathrm{R}}_{\varphi_0}})\otimes\mathrm{W'}^1)\times\widehat{\O^{\mathcal{A}'}_{\varphi_0}}\ar@{->}[r]& \widehat{\O^{\mathrm{R}}_{\varphi_0}} }
%\end{equation}
%Then we can build a local ring morphism $\overline{\eta}$ from $\O^{\mathcal{A}}_{\phi_0}$ to $\O^{\mathcal{A}'}_{\varphi_0}$, given by the composite of local morphism rings, 
%\begin{equation}
%	\overline{\eta}=\mathrm{P}\circ (F')^{-1}\circ\eta\circ F\circ\mathrm{I}
%\end{equation}
%where $\mathrm{I}:\widehat{\O^{\mathcal{A}}_{\phi_0}}\rightarrow\G(\m(\widehat{\O_{\phi_0}})\otimes\mathrm{W}^1)\times\widehat{\O^{\mathcal{A}}_{\phi_0}}$, defined by $\mathrm{I}(\pi(X))=(\mathrm{id},\pi(X))$, and $\mathrm{P}:\G(\m(\widehat{\O_{\varphi_0}^{\mathrm{R}}})\otimes\mathrm{W'}^1)\times\widehat{\O^{\mathcal{A}'}_{\varphi_0}}\rightarrow\widehat{\O^{\mathcal{A}'}_{\varphi_0}}$, defined by $\mathrm{P}(s',\pi'(X))=\pi'(X)$.

\begin{thm}\label{N31}(\textbf{Reduction Theorem}). Let be $\g:=\mathrm{R}\ltimes\n$ be an algebraic Lie algebra such that 
\begin{enumerate}
	\item $\overline{\mathrm{i}}_1:\mathrm{H}^1\left(\n,\n\right)^\mathrm{R}\rightarrow\mathrm{H}^1\left(\g,\g\right)$ is an epimorphism,
	\item $\overline{\mathrm{i}}_2:\mathrm{H}^2\left(\n,\n\right)^\mathrm{R}\rightarrow\mathrm{H}^2\left(\g,\g\right)$  is an isomorphism, and
	\item $\overline{\mathrm{i}}_3:\mathrm{H}^3\left(\n,\n\right)^\mathrm{R}\rightarrow\mathrm{H}^3\left(\g,\g\right)$ is a monomorphism.
\end{enumerate}
Then for every $\A\in\widehat{\got{R}}$, $\mathcal{A}\subset\mathcal{I}$ and $\mathcal{A'}\subset\mathcal{I'}$ admissible sets as above, the map $\mathfrak{I}$ from $\L_n^\mathrm{R}$ to $\L_m$ induces
\begin{enumerate}
	\item[i)]a bijection from $\overline{\Def}\left(\varphi_0,\A\right)^\mathrm{R}$ to $\overline{\Def}\left(\phi_0,\A\right)$, defined by $\left[\varphi\right]\longmapsto\left[\phi\right]$.
	\item [ii)] a local ring morphism, $\eta:\O_{\phi_0}\rightarrow\O^{\mathrm{R}}_{\varphi_0}$ which induces a local ring isomorphism $\bar{\eta}:\O^{\mathcal{A}}_{\phi_0}\rightarrow\O^{\mathrm{R},\mathcal{A}'}_{\varphi_0}$ for $\K=\C$, such that the following diagram  is commutative
	\begin{equation}\label{e3.2}
	\xymatrix{
	\O_{\phi_0} 
	\ar@{->}[r]^{\eta}
	\ar@{->}[d]_{\pi}
	& \O^\mathrm{R}_{\varphi_0}
		\ar@{->}[d]^{\pi'}
	\\
\O^{\mathcal{A}}_{\phi_0}\ar@{->}[r]_{\bar{\eta}}&\O^{\mathrm{R},\mathcal{A}'}_{\varphi_0}}
	\end{equation}
\end{enumerate}
\end{thm}
\begin{proof} 
\textbf{Surjectivity on the classes.} \\
We will show that for all deformation $\phi\in\Def\left(\phi_0,\A\right)$ there is $\varphi\in\Def\left(\varphi_0,\A\right)^\mathrm{R}$ such that $\phi$ and $\mathfrak{I}(\varphi)$ are equivalent under $\G_m(\A)$. 
We will prove for $\phi=\sum_{\alpha}t^{\alpha}\phi_{\alpha}$ by induction on the integer $p\in\N$ the following property:\\
there are $s_{p}=\mathrm{id}+\sum_{|\alpha|=p}t^{\alpha}s_{\alpha}\in\G_m\left(A\right)$ and $\Phi=\sum_{\alpha}t^{\alpha}\Phi_{\alpha}\in\Def\left(\phi_0,\A\right)$ such that
\begin{equation}
	\phi=s_p\ast\Phi,\quad \Phi_{\alpha}=\mathrm{i}_2(\varphi_{\alpha}),
\end{equation}
where $\varphi_{\alpha}\in\mathrm{C}^2(\n,\n)^{\mathrm{R}}$ for $|\alpha|\leq p$.\\
It is obvious if $p=0$. We assume that $\Phi$ satisfies the induction hypothesis for $p$, the deformation equation can be expressed by:
\begin{equation}\label{e37}
\left[\Phi,\Phi\right]=\sum_{\alpha_1}\sum_{\alpha_2}t^{\alpha_1}t^{\alpha_2}\left[\Phi_{\alpha_1},\Phi_{\alpha_2}\right]=0.
\end{equation}
We get:
\begin{equation}\label{e38}
	\left[\Phi,\Phi\right]=2\sum_{\alpha}t^\alpha (d\phi_\alpha-\omega_\alpha)=0,
\end{equation}
with $\omega_\alpha:=\frac{1}{2}\sum^{\alpha_1+\alpha_2=\alpha}_{\alpha_1\neq 0,\alpha_2\neq 0}[\phi_{\alpha_1},\phi_{\alpha_2}]$. 
It follows from the equation (\ref{e38}) that,  
\begin{equation}\label{e39}
\sum_{|\alpha|\leq p+1} t^\alpha d\Phi_{\alpha}\equiv\sum_{|\alpha|\leq p+1} t^\alpha\omega_\alpha\mod t^{p+2}.
\end{equation}
From the equation (\ref{e39}), and the induction hypothesis, hence 
\begin{equation}
	\sum_{|\alpha|=p+1} \overline{t}^{\alpha}d\Phi_{\alpha}=\sum_{|\alpha|= p+1}\overline{t}^{\alpha}\mathrm{i}_3\left(\frac{1}{2}\sum^{\alpha_1+\alpha_2=\alpha}_{\alpha_1\neq 0,\alpha_2\neq 0}\left[\varphi_{\alpha_1},\varphi_{\alpha_2}\right]\right)\in\m/\m^{p+2}\otimes\mathrm{i}_3(\mathrm{C}^3\left(\n,\n\right)^\mathrm{R}).
\end{equation} 
We deduce from Lemma \ref{n33}(1) that $d_{\phi_0}\Phi_{\alpha}\in\mathrm{i}_3(\mathrm{Z}^3\left(\n,\n\right)^\mathrm{R})$  for all $|\alpha|=p+1$. 
Since $\overline{\mathrm{i}}_3$ is injective then  $d_{\phi_0}\phi_\alpha\in\mathrm{i}_3(\mathrm{B}^3(\n,\n)^{\mathrm{R}})$ by Corollary \ref{C3.1}(1). Hence there are $\varphi'_{\alpha}\in\mathrm{C}^2(\n,\n)^{\mathrm{R}}$ with $|\alpha|=p+1$ such that 
$$d\Phi_{\alpha}=d\mathrm{i}\left(\varphi'_{\alpha}\right),$$ 
and 
$$\Phi_{\alpha}\in\mathrm{i}\left(\varphi'_{\alpha}\right)+\mathrm{Z}^2\left(\g,\g\right).$$ 
Since $\overline{\mathrm{i}}_2$ is sujective, it follows from Corollary \ref{C3.1}(2) that there are $\varphi''_{\alpha}\in\mathrm{C}^2(\n,\n)^{\mathrm{R}}$ with $|\alpha|=p+1$ such that 
\begin{equation}
	\Phi_{\alpha}=\mathrm{i}\left(\varphi'_{\alpha}+\varphi''_{\alpha}\right)+ds_{\alpha},
\end{equation}
where $s_{\alpha}\in\mathrm{C}^1\left(\g,\g\right)$. we set:
\begin{equation}
	s_{p+1}=\mathrm{id}+\sum_{|\alpha|=p+1}t^{\alpha}s_{\alpha},
\end{equation}
and 
\begin{equation}
	\varphi_\alpha=\varphi'_\alpha+\varphi''_\alpha,
\end{equation}
for all $|\alpha|=p+1$.
We have
\begin{eqnarray*}
	s_{p+1}\ast\Phi &&=\sum_{|\alpha|\leq p}t^{\alpha}\Phi_{\alpha_{}}+\sum_{|\alpha|=p+1}t^{\alpha_{}}\left(\Phi_{\alpha_{}}-d s_{\alpha_{}}\right)\\&&
	+\mathrm{(term\,of\,degrees\,>p+1)},
\end{eqnarray*}
hence,
\begin{eqnarray*}
	s_{p+1}\ast\Phi=\sum_{|\alpha|\leq p+1}t^{\alpha_{}}\mathrm{i}\left(\varphi_{\alpha_{}}\right)+\mathrm{(term\,of\,degrees\,>p+1)}.
\end{eqnarray*}
This deformation satisfies the property $(p+1)$ and the sequence of deformations $$\left(s_p\circ\cdots\circ s_1\right)\ast\Phi$$ 
converges in the sense of Krull to a limit under the form $\mathfrak{I}(\varphi)$ which is equivalent to $\phi$. Then $\varphi$ is belong to $\Def(\phi_0,\A)^{\mathrm{R}}$, since $i$ is injective.\\
\textbf{Injectivity on the classes.}\\
We suppose that $\varphi_1,\varphi_2\in\Def(\varphi_0,\A)^{\mathrm{R}}$ are two deformations of $\varphi_0$ such that their images $\mathfrak{I}(\varphi_1),\mathfrak{I}(\varphi_2)\in\Def(\phi_0,\A)^{\mathrm{R}}$ are equivalent under $\G_m(\A)$. We will show that there is $\sigma\in\mathfrak{I}(\mathrm{G}^\mathrm{R}_n(\A))$ such that $\mathfrak{I}(\varphi_2)=\sigma\ast\mathfrak{I}(\varphi_1)$ .\\
We reason by induction on the integer $p\in\N$,
there exists $\sigma=\sum_{\alpha}t^{\alpha }\sigma_{\alpha}\in\G_m(\A)$ such that $\sigma_p:=\sum_{|\alpha|\leq p}t^{\alpha} \sigma_{\alpha}\in\mathfrak{I}(\G^{\mathrm{R}}_n(\m))$ and $\sigma\ast\mathfrak{I}(\varphi_1)=\mathfrak{I}(\varphi_2)$. This is obvious if $p=0$. Then 
$$\sigma_p^{-1}\cdot\sigma=\mathrm{id}+\sum_{|\alpha|=p+1}t^{\alpha}\sigma_{\alpha}+\mathrm{(term\,of\,degrees\,>p+1)}.$$
Since $\sigma_p\in\mathfrak{I}(\G_n(\A)^{\mathrm{R}})$, $\mathfrak{I}(\varphi_1),\mathfrak{I}(\varphi_2)\in\mathfrak{I}(\Def(\varphi_0,\A)^{\mathrm{R}})$ imply that the following expression,
\begin{eqnarray*} 
  \sigma_p^{-1}\ast(\sigma\ast\mathfrak{I}(\varphi_1))=(\sigma_p^{-1}\cdot\sigma)\ast\mathfrak{I}(\varphi_1) &&=(\mathrm{id}+\sum_{|\alpha|=p+1}t^{\alpha}\sigma_{\alpha}+\cdots)\ast\mathfrak{I}(\varphi_1)\\
 &&=\sum_{|\alpha|\leq p}t^{\alpha}\left(\mathfrak{I}(\varphi_1)\right)_{\alpha}+\sum_{|\alpha|
=p+1}t^{\alpha}\left(\left(\mathfrak{I}(\varphi_1)\right)_{\alpha}-d_{\phi_0}\sigma_{\alpha}\right)\\&&+\mathrm{(term\,of\,degrees\,>p+1)}\\&&  =\sigma_p^{-1}\ast\mathfrak{I}(\varphi_2)
\end{eqnarray*}
belongs to $\mathfrak{I}(\Def(\varphi_0,\A)^{\mathrm{R}})$. From Lemma \ref{n33}(1), then we have   
$$ d_{\phi_0}\sigma_{\alpha}\in\mathrm{i}_2(\mathrm{Z}^2(\n,\n)^{\mathrm{R}}),$$ 
for all $|\alpha|=p+1$. Since $\overline{\mathrm{i}}_2$ is injective then  $d_{\phi_0}\sigma_{\alpha}\in\mathrm{i}_2(\mathrm{B}^2(\n,\n)^{\mathrm{R}})$ by Corollary \ref{C3.1}(1). Hence there are $s_{\alpha}\in\mathrm{C}^1(\n,\n)^{\mathrm{R}}$ such that
$$d_{\phi_0}\sigma_{\alpha}=d_{\phi_0}\mathrm{i}_1(s_{\alpha})$$
and 
\begin{equation}
	\sigma_{\alpha}\in\mathrm{i}_1(s_{\alpha})+\mathrm{Z}^1(\g,\g),
\end{equation}
for all $|\alpha|=p+1$. Since $\overline{\mathrm{i}}_1$ is surjective, it follows from Corollary \ref{C3.1}(2) that there are $\delta^1_{\alpha}\in\mathrm{ad}\g=\mathrm{B}^1(\g,\g)$ and $\delta^2_{\alpha}\in\mathrm{C}^1(\n,\n)^{\mathrm{R}}$ such that
$$\sigma_{\alpha}=\mathrm{i}_1(s_{\alpha})+\delta^1_{\alpha}+\mathrm{i}_1(\delta^2_{\alpha})$$
for all $|\alpha|=p+1$. The derivations $\delta^1_{\alpha}$ being inner, it follows that they may be lifted to a derivation of $\mathfrak{I}(\varphi_1)$ by Lemma \ref{n6}. Let $\delta^1_{\alpha}(t)$ denote the lifted derivation of $\delta^1_{\alpha}$. Since $\Pi:=\exp(-\sum_{|\alpha|=p+1}t^{\alpha}\delta^1_{\alpha}(t))$ is an automorphism of $\mathfrak{I}(\varphi_1)$ by Lemma \ref{n7}, and using the induction hypothesis, hence we get
$$\sigma_{p+1}:=\sigma_{p}\circ\Pi=\sigma_p+\sum_{|\alpha|=p+1}t^{\alpha}\mathrm{i}_1(s_{\alpha}+\delta^2_{\alpha})+\mathrm{(term\,of\,degrees\,>p+1)}.$$
Then $\sigma\circ\Pi$ verifies the property for $p+1$. The sequence $(\sigma_p)$ converges in the sense of Krull to an automorphism $\sigma_\infty$ which belongs to $\mathfrak{I}(\G_n^{\mathrm{R}}(\A))$ satisfying $\sigma_\infty\ast\mathfrak{I}(\varphi_1)=\mathfrak{I}(\varphi_2)$.\\
2. \textbf{Local Isomorphism $\overline{\eta}$}\\
Let $f:\mathrm{B}\rightarrow\mathrm{C}$ be a morphism of local rings, denote by $\widehat{f}:\widehat{\mathrm{B}}\rightarrow\widehat{\mathrm{C}}$ its unique extension on the completion rings.\\
Set $\A=\widehat{\O^{\mathcal{A}}_{\phi_0}}$ the completion ring of $\O^{\mathcal{A}}_{\phi_0}$. From the surjectivity on the classes, there are $h:\O^{\mathrm{R}}_{\varphi_0}\rightarrow\widehat{\O^{\mathcal{A}}_{\phi_0}}$ an element of $\Def^{\mathrm{R}}(\varphi_0,\widehat{\O^{\mathcal{A}}_{\phi_0}})$ and $s$ an element of $\G_m(\widehat{\O^{\mathcal{A}}_{\phi_0}})$ such that 
\begin{equation}
	h\circ\eta=\mathfrak{I}(h)=s\ast\pi.
\end{equation}
Hence 
\begin{equation}\label{E3.15}
	\widehat{h}\circ\widehat{\eta}=s\ast\widehat{\pi}.
\end{equation}
By Corollary \ref{N21} and Remark \ref{r2.1}, there is a local morphism $\overline{h}:\O^{\mathrm{R},\mathcal{A}'}_{\varphi_0}\rightarrow\widehat{\O^{\mathcal{A}}_{\phi_0}}$ and $s'\in\G_n^{\mathrm{R}}(\widehat{\O^{\mathcal{A}}_{\phi_0}})$ such that 
\begin{equation}\label{E3.16}
\widehat{\overline{h}}\circ\widehat{\pi'}=s'\ast \widehat{h}.
\end{equation}
Since $\widehat{\pi}$ is an epimorphism (because is $\pi$), it follows from the equation (\ref{E3.15}) that $\widehat{h}$ is an epimorphism. From the equation (\ref{E3.16}), we deduce that $\widehat{\overline{h}}$ is an epimorphism. Consequently, $\widehat{\overline{h}}\circ\widehat{\overline{\eta}}:\widehat{\O^{\mathcal{A}}_{\phi_0}}\rightarrow\widehat{\O^{\mathcal{A}}_{\phi_0}}$ is an epimorphism (since $\widehat{\overline{\eta}}$ is an epimorphism, since $\mathrm{i}_2$ is injective). Since $\widehat{\O^{\mathcal{A}}_{\phi_0}}$ is the inverse limite of $(\O^{\mathcal{A}}_{\phi_0}/\m^n(\O^{\mathcal{A}}_{\phi_0}))_n$, it follows that the homomorphisms $(\widehat{\overline{h}}\circ\widehat{\overline{\eta}})_n:\O^{\mathcal{A}}_{\phi_0}/\m^n(\O^{\mathcal{A}}_{\phi_0})\rightarrow\O^{\mathcal{A}}_{\phi_0}/\m^n(\O^{\mathcal{A}}_{\phi_0})$ induced by $\widehat{\overline{h}}\circ\widehat{\overline{\eta}}$ are isomophisms $(n=1,2,\cdots)$, since $\widehat{\overline{h}}\circ\widehat{\overline{\eta}}$ is surjective and $\O^{\mathcal{A}}_{\phi_0}/\m^n(\O^{\mathcal{A}}_{\phi_0})$ is a $\O^{\mathcal{A}}_{\phi_0}/\m(\O^{\mathcal{A}}_{\phi_0})$-module of finite dimension (since $\O^{\mathcal{A}}_{\phi_0}$ is Noetherian ring). Hence $\widehat{\overline{h}}\circ\widehat{\overline{\eta}}$ becomes an isomorphism. 
It follows that $\widehat{\overline{\eta}}$ is a monomorphism and then is an isomorphism (since $\widehat{\overline{\eta}}$ is an epimorphism). We deduce that $\widehat{\overline{\eta}}$ is an isomorphism if and only if is $\overline{\eta}$, by Corollary 1.6, p.282 in \cite{A}, since $\O^{\mathcal{A}}_{\phi_0}$ and $\O^{\mathrm{R},\mathcal{A}'}_{\varphi_0}$ are analytic local rings and $\K=\C$.
\end{proof}
%\begin{cor}If $\g$ is a Lie algebra satisfying the hypotheses of Reduction Theorem then $\g$ is rigid in $\L_m$ iff is $\n$ in $\L^{\mathrm{R}}_n$.
%\end{cor}
%\begin{cor}\label{c4.2}If $\g$ is a formal rigid  Lie algebra in $\L_m$ such that $\overline{\mathrm{i}}_1$ is an epimorphism and  $\overline{\mathrm{i}}_2$ is a monomorphism then $\n$ is formal rigid in $\L_n^{\mathrm{R}}$. In particular any complete Lie algebra such that $\mathrm{H}^1(\n,\g/\n)^{\mathrm{T}}=0$, satisfies those conditions.
%\end{cor}
%\begin{proof} We assume that $\n$ is not formal rigid in $\L_n^\mathrm{R}$, hence there exists a formal deformation, $\n_t\in\mathrm{Def}(\varphi_0,\K[[t]])^{\mathrm{R}}$ such that $\n_t$ and $\n$ are not conjugated under $\G_n(\K[[t]])^{\mathrm{R}}$. It follows that $\mathrm{R}\ltimes\n$ and $\mathrm{R}\ltimes\n_t$ are not conjugated under $\G_m(\K[[t]])$ by the injectivity of \quad$\widehat{}$ \quad on the classes, in Theorem \ref{N31}.
%\end{proof}
\begin{rem}Note that the proof of Theorem \ref{N31}.
\begin{enumerate} 
	\item If $\overline{\mathrm{i}}_2$ is bijective and $\overline{\mathrm{i}}_3$ is injective, the local rings $\O^{\mathcal{A}}_{\phi_0}$ and $\O^{\mathcal{A}',\mathrm{R}}_{\varphi_0}$ are isomorphic, and then the local study of  $\phi_0$ in the slice $\L_m^\mathcal{A}$ is equivalent to that of $\varphi_0$ in the slice $\L^{\mathcal{A}',\mathrm{R}}_n$. Moreover, if $\mathcal{B}$ is another admissible set at $\phi_0$ satisfying the hypothesis of Theorem \ref{N31}, then the local rings $\O^{\mathcal{B}}_{\phi_0}$ and $\O^{\mathcal{A}',\mathrm{R}}_{\varphi_0}$ will be isomorphic. Consequently, the local rings $\O^{\mathcal{A}}_{\phi_0}$ and $\O^{\mathcal{B}}_{\phi_0}$ will be isomorphic too.
	\item The map $\overline{\mathrm{i}}_1:\mathrm{H}^1\left(\n,\n\right)^\mathrm{R}\rightarrow\mathrm{H}^1\left(\g,\g\right)$ can be non injective.
\end{enumerate}
\end{rem}
\subsection{Lie algebras satisfying the hypotheses of Reduction Theorem}
\begin{thm}\label{n3.2}An algebraic Lie algebra $\g=\mathrm{R}\ltimes\n$ with $\mathrm{R}=\mathrm{U}\oplus\mathrm{S}$, $\mathrm{S}$ the semi-simple part of $\mathrm{R}$ and $\mathrm{U}$ the torus part, satisfies the hypotheses of Theorem \ref{N31} iff it belongs to one of the following cases
\begin{enumerate}
	\item $\mathrm{U}\neq 0$, $\g$ is complete, $\mathrm{H}^1\left(\n,\g/\n\right)^\mathrm{R}=0=\mathrm{H}^2(\n,\g/\n)^{\mathrm{R}}$.
	\item $\mathrm{U}=0$, $\mathrm{H}^1\left(\n,\g/\n\right)^\mathrm{R}=0=\mathrm{H}^2(\n,\g/\n)^{\mathrm{R}}$.
\end{enumerate}
\end{thm}
\begin{proof} The canonical morphism $\overline{\mathrm{i}}_p$ defined by the composition $\overline{\beta}_p\circ\overline{\alpha}_p$ with the injection $\overline{\beta}_p:\mathrm{H}^p\left(\n,\g\right)^\mathrm{R}\rightarrow\mathrm{H}^p\left(\g,\g\right)$ given in the Hochschild-Serre spectral sequence and the morphism $\overline{\alpha}_p$ given in the long sequence:
\begin{equation}\label{lon1}
	0\rightarrow\mathrm{U}\rightarrow\mathrm{H}^1\left(\n,\n\right)^\mathrm{R}\stackrel{\overline{\alpha}_1}{\rightarrow}\mathrm{H}^1\left(\n,\g\right)^\mathrm{R}\rightarrow\mathrm{H}^1\left(\n,\g/\n\right)^\mathrm{R}
\end{equation}
\begin{equation}\label{lon2}
\stackrel{}{\rightarrow}\mathrm{H}^2\left(\n,\n\right)^\mathrm{R}\stackrel{\overline{\alpha}_2}{\rightarrow}\mathrm{H}^2\left(\n,\g\right)^\mathrm{R}\rightarrow\mathrm{H}^2\left(\n,\g/\n\right)^{\mathrm{R}}\rightarrow\mathrm{H}^3\left(\n,\n\right)^{\mathrm{R}}
\end{equation}
\begin{equation}\label{lon3}
\cdots\stackrel{}{\rightarrow}\mathrm{H}^p\left(\n,\n\right)^\mathrm{R}\stackrel{\overline{\alpha}_p}{\rightarrow}\mathrm{H}^p\left(\n,\g\right)^\mathrm{R}\rightarrow\mathrm{H}^p\left(\n,\g/\n\right)^{\mathrm{R}}\rightarrow\mathrm{H}^{p+1}\left(\n,\n\right)^{\mathrm{R}}
\end{equation}
The surjectivity of $\overline{\mathrm{i}}_p$ induces that of $\overline{\beta}_p$. If $\overline{\beta}_p$ is surjective then $$\mathrm{H}^p\left(\n,\g\right)^{\mathrm{R}}\simeq\mathrm{H}^p\left(\g,\g\right).$$ 
The Hochschild-Serre spectral sequence gives
\begin{equation}\label{E32}
\mathrm{H}^p\left(\n,\g\right)^{\mathrm{R}}\simeq\mathrm{H}^p\left(\g,\g\right)\simeq\bigoplus_{i+j=p}\mathrm{H}^{i}\left(\mathrm{R},\K\right)\otimes\mathrm{H}^j\left(\n,\g\right)^\mathrm{R}.
\end{equation}
In particular for $p=1,2$, we obtain by using Withehead Lemmas (i.e., $\mathrm{H}^1(\mathrm{S},\K)=\mathrm{H}^2(\mathrm{S},\K)=0$) and the fact that $\overline{\beta}_p$ are assumed surjective:
$$\mathrm{H}^1(\n,\g)^{\mathrm{R}}=(\mathrm{U}^*\otimes\mathrm{Z}(\g))\oplus\mathrm{H}^1(\n,\g)^{\mathrm{R}}$$
and
$$\mathrm{H}^2(\n,\g)^{\mathrm{R}}\simeq(\wedge^2\mathrm{U})^*\otimes\mathrm{Z}(\g)\oplus(\mathrm{U}^*\otimes\mathrm{H}^1(\n,\g)^{\mathrm{R}})\oplus\mathrm{H}^2(\n,\g)^{\mathrm{R}}.$$
Hence
\begin{equation}\label{e4.18}
\mathrm{U}^*\otimes\mathrm{Z}(\g)=(\wedge^2\mathrm{U})^*\otimes\mathrm{Z}(\g)=\mathrm{U}^*\otimes\mathrm{H}^1(\n,\g)^{\mathrm{R}}=0.
\end{equation}
If $\mathrm{U}\neq 0$, then the equation (\ref{e4.18}) is verified iff $\mathrm{Z}(\g)=0=\mathrm{H}^1(\n,\g)^{\mathrm{R}}$. It follows that $\mathrm{H}^1(\g,\g)=0$.\\
$\overline{\mathrm{i}}_2$ is injective iff $\overline{\alpha}_2$ since is $\overline{\beta}_2$. If $\overline{\alpha}_2$ is injective  iff $\mathrm{H}^1(\n,\g/\n)^{\mathrm{R}}=0$ since $\mathrm{H}^1(\n,\g)^{\mathrm{R}}=0$, and thus $\overline{\mathrm{i}}_1$ is surjective, see (\ref{lon1}).\\
If $\mathrm{U}=0$, then the equation (\ref{e4.18}) is trivial and $\overline{\beta}_i$ is an isomorphism (i.e. surjective) with $i=1,2$. Consequently,
$\overline{\mathrm{i}}_1$, $\overline{\mathrm{i}}_2$ and $\overline{\mathrm{i}}_3$ are surjective, bijective and injective respectively iff  $\mathrm{U}=0$, $\mathrm{H}^1\left(\n,\g/\n\right)^\mathrm{R}=0=\mathrm{H}^2(\n,\g/\n)^{\mathrm{R}}$.
\end{proof}
\begin{rem}\label{r4.2} 
\begin{enumerate}
	\item The condition $2$ of Theorem \ref{n3.2} is trivially satisfied for $\mathrm{R}=0$.
	\item If $\mathrm{H}^1(\n,\g/\n)^{\mathrm{R}}=0$ then \quad$\overline{\alpha}_1$ is a surjection and $$\mathrm{H}^1(\g,\g)\simeq(\mathrm{U}^*\otimes\mathrm{Z}(\g))\oplus\frac{\mathrm{H}^1(\n,\n)^{\mathrm{R}}}{\mathrm{U}}.$$
\end{enumerate}
\end{rem}
\section{Formal Rigidity}
\subsection{Formal Rigidity and Dimension}
Let $\phi_0\in\L_m(\K)$ be a point and $\mathcal{A}$ an admissible set of $\mathcal{I}$ at $\phi_0$. The Krull dimension $d$ of  the completion ring $\widehat{{\O}^\mathcal{A}_{\phi_0}}$ of $\O^\mathcal{A}_{\phi_0}$ for the topology $\m(\O^\mathcal{A}_{\phi_0})$-adic, is the maximal number of elements $t_1,...,t_d$ of $\m(\O^\mathcal{A}_{\phi_0})$ such that the subring consisting of formal power series in $t_1,...,t_d$, is isomorphic to $\K[[T_1,...,T_d]]$, the formal power series ring with $d$ indeterminates on $\K$. This gives the dimension of $\L^{\mathcal{A}}_{m,\phi_0}$ at point $\phi_0$. The formal rigidity is the rigidity with base $\A=\K\left[\left[T\right]\right]$.\\
If we suppose that there is a valuation on $\K$, for example $\C$, we can define a separated strong topology on $\K^m$ and use the notion of convergence of series (analyticity). In this case the formal rigidity is equivalent to the analytic rigidity, by M. Artin's theorem, \cite{A},\cite{C7}. This rigidity is also equivalent to the orbit is an open set in the sens of Zariski and in the sens of strong topology. We call it the geometric rigidity at point $\phi_0$.
\begin{thm}\label{t4.1}For all admissible set $\mathcal{A}$ of $\mathcal{I}$ at $\phi_0$, then the following conditions are equivalent
\begin{enumerate}
\item $\dim_\K\O_{\phi_0}^\mathcal{A}<\infty$;
\item the Krull dimension $d$ of $\widehat{\O_{\phi_0}^\mathcal{A}}$ is null; $\phi_0$ is an isolated point of $\L_m^{\phi_0,\mathcal{A}}(\K)$
	\item the elements of the maximal ideal of $\O_{\phi_0}^\mathcal{A}$ are nilpotent;
	\item the ideal generated by the elements $(\xi^{\alpha})_{\alpha\in\mathcal{A}}$ of $\O_{\phi_0}$ contains $\m(\O)^p$ for some $p\in\N^*$;
\item $\phi_0$ is formally rigid in $\L_m$;
\item Moreover if $\K$ is a valued field, the orbit $\left[\phi_0\right]$ is a Zariski's open set in $\L_m(\K).$
\end{enumerate}
\end{thm}
\begin{proof}$1\Rightarrow 2.$ If $t$ is an element of $\m(\O^{\mathcal{A}}_{\phi_0})$ then the sequence $(t^i)_{i\in\N}$ linearly generates a vector subspace of finite dimension. Consequently the subring $\K[[t]]$ of $\widehat{\O^{\mathcal{A}}_{\phi_0}}$ cannot be a ring of series power, so $d=0$.\\
$2\Rightarrow 3.$ If $t\in\m(\O^{\mathcal{A}}_{\phi_0})$ then the ring $\K[[t]]$ is isomorphic to the ring $\K[T]/(T^p)$ where $p\in\N^*$, and $t^p=0$.\\
$3\Rightarrow 4.$ The ring $\O^\mathcal{A}_{\phi_0}$ being the quotient of $\O_{\phi_0}$ by the ideal generated by the elements $(\xi^\alpha)_{\alpha\in\mathcal{A}}$ and of finite type, the condition 3 means there is $p$ such that $\m(\O_{\phi_0}^\mathcal{A})^p=0$, i.e., $\m(\O_{\phi_0})^p\subset\sum_{\alpha\in\mathcal{A}}\xi^\alpha\O_{\phi_0}$.\\
$4\Rightarrow 5.$ By Corollary \ref{N21}, any deformation with formal base, $h:\O_{\phi_0}\rightarrow\K[[t]]$  is equivalent to $h_0=\overline{h}_0\circ\pi:\O_{\phi_0}\stackrel{\pi}{\longrightarrow}\O^{\mathcal{A}}_{\phi_0}\stackrel{\overline{h}_0}{\longrightarrow}\K[[t]].$ The hypothesis $4$ is equivalent to $\m(\O^{\mathcal{A}}_{\phi_0})^p=0$ where $p\in\N^*$. Then any element of  $\m(\O^{\mathcal{A}}_{\phi_0})$ is nilpotent and the morphism $\overline{h}_0$ is trivial since $\K[[t]]$ is integral. It follows that $h_0$ is also trivial on $\m(\O_{\phi_0})$ and the deformation $h$ is trivial. Then $\phi_0$ is $\K[[t]]$-rigid.\\
$5\Rightarrow 1.$ If the $\K$-dimension of $\O^{\mathcal{A}}_{\phi_0}$ is infinite then there is a variable $t_i$ among the generators $(t_j)_{j=1...p}$ of $\m(\O^{\mathcal{A}}_{\phi_0})$ such that the ideal $J_i$ generated by the $t_j\neq t_i$ gives a quotient $\O^{\mathcal{A}}_{\phi_0}/J_i$ of infinite dimension over $\K$. We will obtain a non trivial subring of $\K[\overline{t}_i]$ and this quotient defines a non trivial deformation $\overline{h}_0$, as well as $h_0=\overline{h}_0\circ\pi$, this contradicts the hypothesis 5.\\
$6\Leftrightarrow5.$ is classical since $\K$ is valued, see above.
\end{proof}

It is shown in \cite{C1} if $\g$ is formal rigid in $\L_m$ that $\g$ is algebraic. It follows that 
it admits a Chevelley's decomposition, i.e., $\g=\mathrm{R}\ltimes\n$ where $\n$ is the nilpotent radical and  $\mathrm{R}=\mathrm{S}\oplus\mathrm{U}$ is consisted of a Lie subalgebra of Levi $\mathrm{S}$, and a torus $\mathrm{U}$ on $\n$.
\begin{cor}\label{c4.1}\
\begin{enumerate}
	\item If $\g=\mathrm{R}\ltimes\n$ such that $[\g,\n]=\n$ and $\mathrm{H}^2(\n,\g/\n)^2=0$, then $\g$ is formal rigid in $\L_m$ if and only if is $\n$ in $\L_n^{\mathrm{R}}$.
	\item If $\g=\mathrm{R}\ltimes\n$ such that $[\g,\n]=\n$ and $\g$ is formal rigid in $\L_m$, then is $\n$ in $\L_n^{\mathrm{R}}$.
\end{enumerate}
\end{cor}
\begin{proof} 1. Let $\phi_0$ (resp. $\varphi_0$) denote the law of $\g$ (resp. $\n$). For all admissible set $\mathcal{A}\subset\mathcal{I}$ (resp. $\mathcal{A'}\subset\mathrm{I'}$) at $\phi_0$ (resp. $\varphi_0$), defined as in Theorem \ref{N31}, the local morphism $\overline{\eta}:\O^{\mathcal{A}}_{\phi_0}\rightarrow\O^{\mathrm{R},\mathcal{A'}}_{\varphi_0}$ is an isomorphism since $[\g,\n]=\n$ and $\mathrm{H}^2(\n,\g/\n)^2=0$, it follows that $\O^{\mathcal{A}}_{\phi_0}$ is a finite-dimensional space if and only if is the dimension of $\O^{\mathrm{R},\mathcal{A'}}_{\phi_0}$. Then we deduce the statement from Theorem \ref{t4.1}.\\
2. $\mathrm{\eta}$ is surjective since $[\g,\n]=\n$, then we deduce the statement from Theorem \ref{t4.1}.
\end{proof}
\subsection{Criterion of formal Rigidity}
We denote by $\mathfrak{R}_m$ the Zariski closed of $\L_m$ consisting of solvable Lie algebras. The formal rigidity of $\g$ in $\mathfrak{R}_m$ means its orbit is open in $\mathfrak{R}_m$. 
Let $\mathrm{C}^\infty\g$ denote the intersection of the terms $\mathrm{C}^p\g$ of the descending central series with $\mathfrak{s}\ltimes\overline{\n}$ its Levi decomposition. It follows from \cite{C4} that there is a bijection between the isomorphism classes of complete Lie algebras $\g$ and the isomorphism classes $\mathcal{C}^\infty$ consisting of the parts $\mathfrak{s}\ltimes\overline{\n}$ satisfying the following conditions
\begin{equation}\label{e4.1}
	\overline{\n}^{\mathrm{S}\oplus\mathrm{U}}\cap\mathrm{Z}(\overline{\n})=0,
\end{equation}
	and
\begin{equation}\label{e4.2}
	(\mathrm{S}\oplus\mathrm{U})\cdot\overline{\n}+\left[\overline{\n},\overline{\n}\right]=\overline{\n},
\end{equation}
where $\mathrm{S}=\mathrm{ad}\mathfrak{s}$, $\mathrm{U}$ is a maximal torus of $\mathrm{Der}(\overline{\n})^{\mathrm{S}}$. Moreover there is an isomorphism of Lie algebras  
\begin{equation}\label{E6.4}
	\g\simeq\mathrm{R}\ltimes\mathrm{N}^{\mathrm{U}}|_{\mathrm{m}}\oplus\mathfrak{m},
\end{equation}
where $\mathrm{R}=\mathrm{S}\oplus\mathrm{U}$, $\mathfrak{m}=\mathrm{R}\cdot\overline{\n}$ and $\mathrm{N}$ is the  largest nilpotency ideal of $(\mathrm{Der}(\overline{\n}))^{\mathrm{S}}$.
In particular there is a bijection between the isomorphism classes of complete solvable Lie algebras and the isomorphism classes of nilpotent Lie algebras $\overline{\n}\in\mathcal{C}^\infty_{N}$ satisfying the following conditions
\begin{equation}\label{E6.5}
	\overline{\n}^{\mathrm{T}}\cap\mathrm{Z}(\overline{\n})=0,
\end{equation}
and
\begin{equation}\label{E6.6}
	\mathrm{T}\cdot\overline{\n}+\left[\overline{\n},\overline{\n}\right]=\overline{\n}
\end{equation}where $\mathrm{T}$ is a maximal torus on $\overline{\n}$ and $\g\simeq\mathrm{T}\ltimes\mathrm{N}^{\mathrm{T}}|_{\mathrm{m}}\oplus\mathfrak{m}$.\\
All direct product of nilpotent complex Lie algebras up to dimension $6$ belong to the class $\mathcal{C}^\infty_{N}$, see \cite{C5}. These nilpotent Lie algebras $\overline{\n}$ give complete solvable Lie algebras $\g$ satisfying $\n(\g)=\overline{\n}$, this is equivalent to $\left[\g,\n\right]=\n$, where $\n(\g)$ is the nilpotent radical of $\g$. \\
If $\dim_\K\overline{\n}$ $\geq 7$, there are nilpotent Lie algebras do not belong to $\mathcal{C}^\infty_{N}$, in particular the characteristically nilpotent Lie algebras, \cite{Di},\cite{L}. There are also complete solvable Lie algebras of dimension $9$ which do not satisfy the property $\n(\g)=\overline{\n}$, \cite{C4}. \\
Let $\g=\mathrm{R}\ltimes\mathfrak{n}$ be an algebraic Lie algebra of finite dimension Lie algebra with $\mathrm{R}=\mathrm{U}\oplus\mathrm{S}$ reductive Lie subalgebra and $\mathfrak{n}$ the largest nilpotent ideal. It is easy to see that the group $\mathrm{H}^1(\n,\g/\n)^{\mathrm{R}}$ is isomorphic to $\mathrm{Hom}_{\mathrm{R}}(\n/\left[\n,\n\right],\g/\n)$ where $\g/\n$ is identified with the $\mathrm{R}$-module adjoint $\mathrm{R}$. We then have
\begin{equation}
\mathrm{H}^1(\n,\g/\n)^{\mathrm{R}}=\mathrm{Hom}_{\mathrm{S}}((\n/\left[\n,\n\right])^{\mathrm{U}},\mathrm{S})\oplus\mathrm{Hom}_{\K}((\n/\left[\n,\n\right])^{\mathrm{R}},\mathrm{U}).
\end{equation}
\begin{lem}\label{n3.5}The following properties are equivalent
\begin{enumerate}
\item $\mathrm{H}^1(\n,\r/\n)^{\mathrm{R}}=\mathrm{Hom}_{\K}((\n/\left[\n,\n\right])^{\mathrm{R}},\mathrm{U})=0.$
	\item $\r=\n$ or $\n=\left[\n,\g\right]$ where $\r$ is the radical of $\g$.
\end{enumerate}
\end{lem}
\begin{proof} 1. is equivalent to $\mathrm{U}=0$ or $\mathrm{R}.\n+\left[\n,\n\right]=\n$, i.e., $\left[\g,\n\right]=\n$. 
\end{proof}\\
The image of $\mathrm{H}^1(\n,\g/\n)^{\mathrm{R}}$ by $\partial_1$ in the long sequence of cohomologies (\ref{lon1}) verifies the following:
\begin{lem} \label{n3.6}Let $f:\n\rightarrow\g/\n$ be a morphism of $\mathrm{R}$-modules such that $f(\left[\n,\n\right])=0$ which represents a class of \quad$\mathrm{H}^1(\n,\g/\n)^{\mathrm{R}}$. The differential of $f$ in $\mathrm{C}^1(\n,\g)^{\mathrm{R}}$ gives a representative $\xi\in\mathrm{Z}^2(\n,\n)^{\mathrm{R}}$ of the image by $\partial_1$ of the class of $f$ which verifies for all $(x,y,z)\in\n^3$:
$$(\xi\circ\xi)(x,y,z)=\oint_{(x,y,z)}\left[\left[f(x),f(y)\right],z\right].$$
\end{lem}
\begin{proof} The condition of cocycle in $\mathrm{C}^1(\n,\g/\n)^{\mathrm{R}}$ implies that $f(\left[x,y\right])=0$ and we obtain for $x,y\in\n$:
$$\xi(x,y)=d f(x,y)=\left[f(x),y\right]+\left[x,f(y)\right]\in\n,$$
where $d$ is the operator differential of the complex $\mathrm{C}(\n,\g)^{\mathrm{R}}$. It follows that
$$\xi(\xi(x,y),z)=ad f(\xi(x,y))(z)-ad f(z)(\xi(x,y)).$$
It follows from this and by using the $\mathrm{R}$-invariance for $f$ (i.e $[f(x),y]=f[x,y]$) that
\begin{equation}\label{e4.7}	\xi(\xi(x,y),z)=2\left[\left[f(x),f(y)\right],z\right]+\left[\left[f(x),y\right],f(z)\right]-\left[f(z),\left[x,f(y)\right]\right].
\end{equation}
Since 
$$\xi\circ\xi(x,y,z)=\oint_{(x,y,z)}\xi(\xi(x,y),z),$$ 
thus we obtain the result by summing the expression (\ref{e4.7}) and by using the Jacobi identity.
\end{proof} \\

%\begin{prop}If a no perfect formal rigid  Lie algebra $\g=\mathrm{R}\ltimes\n$ in $\L_m$, then $\g$ is complete.
%\end{prop}
%\begin{proof} We suppose that $\mathrm{Z}(\g)\neq 0$. The adjoint representation of $\mathrm{R}$ in $\mathrm{Der}(\n,\g)$ being semi-simple, then
%\begin{equation}
%	\mathrm{H}^1(\n,\g)=\mathrm{H}^1(\n,\g)^{\mathrm{R}}\oplus\mathrm{R}\cdot\mathrm{H}^1(\n,\g).
%\end{equation}
%We can see $\mathrm{R}\cdot\mathrm{H}^1(\n,\g)=0$, since $\mathrm{R}\cdot\mathrm{Der}(\n,\g)\subset\mathrm{ad}\g$. It follows from Dixmier's theorem,\cite{D}, that $\mathrm{H}^1(\n,\g)\neq 0$, since $\n$ is nilpotent and there is a trivial $\n$-submodule $\mathrm{Z}(\g)$ in $\g$. It follows from above that $\mathrm{H}^1(\n,\g)^{\mathrm{R}}\neq 0$, this contradicts the fact any no perfect Lie algebra $\g$ in $\mathfrak{L}_m$ satisfies that $\mathrm{H}^1(\n,\g)^{\mathrm{R}}=0$, \cite{C1}. 
%\end{proof}
\begin{prop}\label{N3.2}Any algebraic Lie algebra $\g=\mathrm{R}\ltimes\n$ such that $\n$ is rigid in $\L_n^{\mathrm{R}}$ and $\mathrm{H}^1(\g,\g)=\mathrm{H}^0(\g,\g)$, satisfies that $\r=\n$ or $\n=\left[\g,\n\right]$.
\end{prop}
\begin{proof} 
%We suppose that the radical of $\g$ is not nilpotent, i.e., $\r\neq\n$. Then the rigid Lie algebra satisfies $\dim_\K\mathrm{U}=1$ or $\mathrm{Z}(\g)=0$, cf \cite{C1}.\\
%If $\mathrm{Z}(\g)\neq 0$, it follows from above that the codimension of $\left[\g,\g\right]$ is equal to $1$. Since $\r\neq\n$ hence $\left[\g,\g\right]\subset\mathrm{S}\oplus\n\subsetneq\mathrm{S}\oplus\r$. It follows that $\left[\g,\g\right]=\mathrm{S}\oplus\n$, i.e., $\left[\g,\n\right]=\n$.\\
%Since $\g$ is a complete Lie algebra then $\g$ satisfies the statement, cf (\ref{e4.2}). 
%It follows from Corollary \ref{c4.2} if $\g$ is complete and rigid that $\n$ is rigid in $\L_n^{\mathrm{R}}$. 
The hypothesis $\mathrm{H}^1(\g,\g)=\mathrm{H}^0(\g,\g)$ implies that $\mathrm{H}^1(\n,\g)^\mathrm{R}=0$ and the long sequence of cohomology gives an injective map $\partial_1$, $$0\rightarrow\mathrm{H}^1(\n,\g/\n)^{\mathrm{R}}\stackrel{\partial_1}{\rightarrow}\mathrm{H}^2(\n,\n)^{\mathrm{R}}\rightarrow\cdots$$
The map $\partial_1$ maps a non zero class $\overline{f}$ of $\mathrm{H}^1(\n,\r/\n)^{\mathrm{R}}$ onto a class $\overline{\xi}\neq 0$ of $\mathrm{H}^2(\n,\n)^{\mathrm{R}}$ satisfying Lemma \ref{n3.6}. Hence we have $\xi\circ\xi=0$ since the image of $f$ is in $\mathrm{U}$. We obtain a non trivial deformation $\varphi_0+t\xi$ of $\n$ in $\L^{\mathrm{R}}_n$, so this contradicts the rigidity of $\n$. Hence we have $\mathrm{H}^1(\n,\r/\n)^{\mathrm{R}}=0$ this expresses the condition $2.$ of Lemma \ref{n3.5}.
\end{proof}\\
%If $\g$ is a rigid complete solvable Lie algebra in $\L_m$, then $\n(\g)=\overline{\n}$ by Proposition \ref{N3.2}
%and $\g\simeq\mathrm{T}\ltimes\mathrm{N}^\mathrm{T}\oplus\mathfrak{m}$ where $\mathfrak{m}=\mathrm{T}\cdot\n$, $\mathrm{N}$ is the largest nilpotency ideal of $\Der(\n)$ and $\mathrm{T}$ is a maximal torus on $\n$.
%Hence $\n_0\simeq\mathrm{N}^{\mathrm{T}}$ and  $\mathrm{N}^{\mathrm{T}}=\mathrm{ad}\n_0$ where $\n_0$ is the weightspace of zero of $\n$. The last condition may be expressed by
%\begin{equation}\label{E6.9}
%	\mathrm{Der}_\K(\n)^{\mathrm{T}}=\mathrm{T}\oplus\mathrm{ad}\n_0.
%\end{equation}
The space of multilinear maps $\mathrm{C}^i(\V_1,\V_2)$ may be identified with the space of linear maps $\mathrm{End}(\wedge^i\V_1,\V_2)$ by associating to any element $h\in\mathrm{C}^i(\V_1,\V_2)$ an element  $\widetilde{h}\in\mathrm{End}(\wedge^i\V_1,\V_2)$ defined by $\widetilde{h}(x_1\wedge...\wedge x_i)=h(x_1,...,x_i)$ for all $x_i\in\V_1$, where $\V_1$ and $\V_2$ are $\K$-module. \\
The rigidity of $\g$ in $\mathfrak{R}_m$ in the case $\mathrm{H}^1(\g,\g)=\mathrm{H}^0(\g,\g)$ can be characterized by the following 
\begin{thm}\label{t6.2}We suppose that there is a valuation on $\K$. If $\g=\mathrm{T}\ltimes\n$ is a solvable Lie algebra satisfying $\mathrm{H}^1(\g,\g)=\mathrm{H}^0(\g,\g)$. Then the two following properties are equivalent:
\begin{enumerate}
	\item [i)]$\n$ is formal rigid in $\L_n^{\mathrm{T}}$;
	%\item $\mathrm{Der}(\n)^{\mathrm{T}}=\mathrm{T}\oplus\mathrm{ad}\n_0$, 
	\item [ii)] we have $[\g,\n]=\n$ and $\g$ is formal rigid in $\mathfrak{R}_m$.
\end{enumerate}
\end{thm} 
\begin{proof} $i)\Rightarrow ii).$  The rigidity of $\n$ in $\L_n^{\mathrm{T}}$ and the hypothesis $\mathrm{H}^1(\g,\g)=\mathrm{H}^0(\g,\g)$ imply that $[\n,\g]=\n$ by Proposition \ref{N3.2}. Let $\phi=\sum_{k\geq 0}t^k\phi_k$ be a deformation of $\phi_0$ with $\phi_0$ the Lie multiplication of $\g$ and $\varphi_0$ its restriction on $\n$. We may assume that $\phi$ is $\mathrm{T}$-invariant (up to conjugation) under the adjoint action of torus $\mathrm{T}$ on $\g$, see \cite{C7}. The $\mathrm{T}$-invariance respects the different weightspaces, $\n$ admits a decomposition $$\n=\mathfrak{m}\oplus\n_0$$ 
where $\mathfrak{m}=\oplus_{\alpha\neq 0}\n_\alpha$. It follows from $[\n,\g]=\g$ that the map $\widetilde{\phi_0}$ defines a surjection from $(\wedge^2\n)_0$ to $\n_0$.
Let $p$ denote the smallest index such that $\phi_p(\n\times\n)$ is not contained in $\n$. The restriction $\phi_{p}$ to $\n\times\n$ admits a decomposition $\varphi_p+f$ satisfying
$$\varphi_p(\n\times\n)\subset\n\quad\mathrm{and}\quad f(\n\times\n)\subset\mathrm{T},$$ 
where $\varphi_k$ (for $k\in\N$) is the projection onto $\n$, relatively to the direct sum $\mathrm{T}\oplus\n$, of the restriction $\phi_k|_{\n\times\n}$.
The Jacobi identity of $\phi$ modulo $t^{p+1}$ implies that 
$$f\in\mathrm{Z}^2(\n,\g/\n)^{\mathrm{T}}.$$ 
Let $\mathrm{F}$ be a complement of  $\mathrm{F}_0:=\ker(\widetilde{\phi}_0|_{({\wedge^2\n})_0})$ in $({\wedge^2\n})_0$, the kernel of the restriction of the map $\widetilde{\phi_0}$ on $(\wedge^2\n)_0$. The linear map 
$$\widetilde{f}|_{(\wedge^2\n)_0}:(\wedge^2\n)_0\rightarrow\mathrm{T},$$ 
may be decomposed as $\widetilde{f}_1+\widetilde{f}_2$ such that $$\widetilde{f}_1(\mathrm{F})=\widetilde{f}_2(\mathrm{F_0})=0.$$ 
We remark that the maps $\widetilde{f}_i$ and $\widetilde{f}$ may be extended to $\wedge^2\n$ by vanishing on a $\mathrm{T}$-invariant complement of $(\wedge^2\n)_0$ in $\wedge^2\n$ since they are $\mathrm{T}$-invariant with $i=1,2$. The bijection $\widetilde{\phi}_0|_{\mathrm{F}}:\mathrm{F}\rightarrow\n_0$ permits to define a linear map $L:\g\rightarrow\mathrm{T}$ given by 
$$L|_{\n_0}=(\widetilde{f}_2\circ\widetilde{\varphi}^{-1}_0)|_{\n_0}\quad\mathrm{and} \quad L|_{\mathrm{T}\oplus\m}=0.$$ 
It follows that the element $\widetilde{f}_2$ is equal to $(L\circ\widetilde{\varphi}_0)|_{\wedge^2\n}$ and $f_2$ belongs to $\mathrm{B}^2(\n,\g/\n)^{\mathrm{T}}$. The truncated of order $p$ of the deformation $\Psi=s\ast\phi$ with $s=\mathrm{id}-t^p L$ is given by 
$$\Psi^{(p)}=\phi^{(p)}+t^p d(L)$$ 
which its restriction on $\n\times\n$ may be written as 
$$\varphi^{(p)}+t^p(f-L\circ\varphi_0),\quad\mathrm{with}\quad \varphi^{(p)}=\sum_{k\leq p}t^k\varphi_k,$$
which is equal to 
$$\varphi^{(p)}+t^p f_1$$ 
where $f_1$ is considered as an element of $\mathrm{H}^2(\n,\g/\n)^{\mathrm{T}}$. We then have $$\widetilde{\Psi}^{(p)}(\mathrm{F})=\widetilde{\varphi}^{(p)}(\mathrm{F}),$$
since $\widetilde{f}_1$ is vanished on $\mathrm{F}$. Hence the restriction of $\widetilde{\Psi}^{(p)}$ to $\mathrm{F}$ defines a deformation consisting of linear maps from $\mathrm{F}$ to $\n_0$ which coincides with $\widetilde{\varphi}_0$ for $t=0$. The equality 
$$\widetilde{\varphi}_0(\mathrm{F})=\n_0,$$ 
(i.e., $[\n,\g]=\n$) implies this deformation is contained in the open set consisting of linear maps of maximal rank and it follows that 
$$\widetilde{\varphi}^{(p)}(\mathrm{F})=\n_0\quad\mathrm{for}\quad |t| \quad\mathrm{small}.$$
It is clear that 
$$\widetilde{\Psi}^{(p)}(\mathrm{F}_0)\subset\mathrm{T}\oplus\n_0,$$ 
and if $\widetilde{f}_1\neq 0$ then $\widetilde{\Psi}^{(p)}(\mathrm{F}_0)$ is not contained in $\n_0$.  
It follows if $\widetilde{f}_1$ was different from zero that $\dim_\K\widetilde{\Psi}^{(p)}((\wedge^2\n)_0)$ would be greater than $\dim_\K\n_0$ and then  $\dim_\K\widetilde{\Psi}((\wedge^2\n)_0)$ would be greater than $\dim_\K\n_0$. The deformation $\Psi$ being staying in $\mathfrak{R}_m$, its derived ideal $\Psi(\g\times\g)\subset\n$ which contains $\m$ and a subspace of $\mathrm{T}\oplus\n_0$ of dimension greater than $\dim_\K\n_0$. Consequently, the radical nilpotent of $\Psi$ would be of dimension greater than $\dim_\K\n$, that is forbidden because the dimension of this ideal cannot increase at neighbourhood of this point in $\mathfrak{R}_m$, cf \cite{C1}. The hypothesis $f_1\neq 0$ is not considered. By induction on the integer $p$, we may construct a deformation $\Psi$ which is equivalent to $\phi$ such that $\Psi(\g\times\g)=\n$. The restriction $\varphi$ of $\Psi|_{\n\times\n}$ to $\n\times\n$ is a $\mathrm{T}$-invariant deformation of $\varphi_0$. It follows from the hypothesis i) that $\varphi$ is trivial, i.e., there is an automorphism $s=\mathrm{id}+\sum_{k>0}t^kL_k$ of $\n$ such that $\varphi_0=s\ast\varphi$. Let $s'$ denote the automorphism of $\g$ such that $s'|_{\n}=s$ and $s'|_{\mathrm{T}}=\mathrm{id}$, then $s'\ast\Psi=\phi_0$ and so $\g$ is rigid.\\
$ii)\Rightarrow i)$. 
$\K$ being valued, we assume that $\n$ is not analytical rigid in $\L_n^\mathrm{R}$. There is then a analytical deformation $\n_t$ of $\n$ with law $\varphi_t$ and vector space $\n=\K^n$ such that for all $\epsilon>0$, there is $0<t<\epsilon$ satisfying that $\varphi_t$ does not belong to the $\mathrm{Gl}(\n)^{\mathrm{T}}_0$-orbit of $\varphi_0$ in $\L_n^{\mathrm{T}}$. If $\epsilon$ is enough small, then $\varphi_t$ will not belong to the $\mathrm{Gl}(\n)$-orbit of $\varphi_0$ in $\L_n$ since its trace on $\L_n^{\mathrm{T}}$ is finite union of irreducible $\mathrm{Gl}(\n)^{\mathrm{T}}_0$-orbits, by Proposition \ref{p}. Hence $\n_t$ will be not isomorphic to $\n$.
The deformed Lie algebra $\g_t:=\mathrm{T}\oplus\n_t$ belongs to $\L_m$ and the dimension of the derived ideal cannot decrease in neighbourhood of $t=0$, we then have $[\g_t,\g_t]=\n_t$, since $\left[\g,\g\right]=\n$. If the ideal $\n_t$ is not nilpotent, there will exist a Lie subalgebra of Levi $\mathfrak{s}$ of $\n_t$ such that for t fixed, $\g_t$ admits a Chevalley's decomposition $\mathrm{T}\oplus\mathfrak{s}\oplus\n'_t$ where $\n'_t$ is the radical nilpotent of $\n_t$ and $\g_t$. The action of $\mathrm{T}$ being fixed, we have $\m:=[\mathrm{T},\g_t]=[\mathrm{T},\n'_t]\subset\n'_t$ and the vector space $\m+\varphi_t(\m\times\m)\subset\n'_t$ is the image of a linear map which depends on $t$ and whose the rank cannot decrease at neighbourhood of $t=0$. It follows that $\n'_t=\n$ and $\mathfrak{s}=0$ such that $\g_t$ belongs to $\mathfrak{R}_m$. Since $\g$ is rigid in $\mathfrak{R}_m$, then $\g_t$ is isomorphic to $\g$ for all $|t|<\epsilon$. Consequently, the maximal nilpotent ideals $\n_t$ and $\n$ are isomorphic, hence a contradiction.
\end{proof}
\begin{cor}If $\mathrm{H}^1(\g,\g)=\mathrm{H}^0(\g,\g)\neq 0$, we have:\\
$\n$ is rigid in $\L^{\mathrm{T}}_n$ if and only if is $\g$ in $\mathfrak{R}_m$.
\end{cor}
\begin{rem}
\begin{enumerate}
	\item A complete solvable Lie algebra $\g$ with $\left[\g,\n\right]=\n$ is equivalent to the semi-direct product $\mathrm{T}\ltimes\n$ with $\mathrm{T}$ a maximal torus on $\n$ and $\n\in\mathcal{C}^\infty_{N}$ satisfying $\mathrm{Der}(\n)^{\mathrm{T}}=\mathrm{T}\oplus\mathrm{ad}\n_0$.
	\item Note the following observations might be induced from Theorem \ref{t6.2}.
If we assume that the deformation $\phi$ of $\phi_0$ does not belong to $\mathfrak{R}_m$ and $\phi_0$ is complete. Let $f_1\in\mathrm{Z}^2(\n,\g/\n)^{\mathrm{T}}$ a representative of a nonzero class which is defined on $\mathrm{F}_0$. In particular $f_1$ will not be vanished on a pair $(x_\alpha,x_{-\alpha})\in\n_\alpha\times\n_{-\alpha}$ satisfying $\left[x_\alpha,x_{-\alpha}\right]=0$, with $\alpha\neq 0$. It is possible that there is a $S$-triplet $(x_\alpha,x_{-\alpha},h_\alpha)$, i.e., $$\phi(x_\alpha,x_{-\alpha})=f_1(x_{\alpha},x_{-\alpha})=h_\alpha\in\mathrm{T}-\left\{0\right\},$$
 hence there exists a nontrivial Levi subalgebra of the deformed Lie algebra. In particular we would have $$f_1(\n_0\times\n_0)\subset\sum_{\alpha\neq 0} f_1(\n_\alpha\times\n_{-\alpha}).$$ 
Because the condition Eq. (\ref{E6.6}) of $1.$ implies that $\n_0\subset\sum_{\alpha\neq 0}\left[\n_\alpha,\n_{-\alpha}\right]$ and each element of $f_1(\n_0,\n_0)$ is sum of vectors $f_1(x_0,[x_\alpha,x_{-\alpha}])$, for all $(x_\alpha,x_{-\alpha})\in\n_\alpha\times\n_{-\alpha}$, which can be written
$$-f_1(x_\alpha,[x_{-\alpha},x_0])-f_1(x_{-\alpha},[x_0,x_\alpha])\in f_1(\n_\alpha\times\n_{-\alpha}),$$
thanks to the condition of cocycle, with $x_0\in\mathrm{T}$.\\
This situation imposes on the deformed algebra to be of the form $\tau\oplus\mathfrak{s}\ltimes\n'_t$, where the nilpotent radical $\n'_t$ is of dimension less than $n$, depends on the parameter $t$ which is nonzero and small. This Lie algebra would have the same orbital dimension as $\g$ which is maximal since $\g$ is complete, and $t$ would be a true parameter. One would obtain a family of laws parameterized which are isomorphic where $\g$ is rigid in $\got{R}_m$ but is not in $\L_m$. This same phenomenon can more generally exist in $\L_m$ and corresponds to an increase in the subalgebra of Levi, what will be characterized by the group $\mathrm{H}^2(\n,\g/\n)^{\mathrm{R}}$. We will obtain the existence of continuous families. Up to now, these phenomena were not observed in the known examples.
\end{enumerate}
\end{rem}
\textbf{Particular cases}. Let $\n$ be a Lie algebra, a derivation $\delta$ of $\n$ has a strict positive spectrum, denoted by $\delta>0$ if all its eigenvalues are in $\Q^+-\left\{0\right\}$. The algebraic Lie algebra $\Der_\K(\n)$ admits a Chevalley's decomposition $\mathrm{S}\oplus\mathrm{A}\oplus\mathrm{N}$ where $\mathrm{S}$ is a Levi subalgebra, $\mathrm{A}$ is an Abelian subalgebra such that $\left[\mathrm{S},\mathrm{A}\right]=0$, $\mathrm{N}$ is the largest ideal of nilpotence and $\mathrm{R}=\mathrm{S}\oplus\mathrm{A}$ the reductive part.
\begin{prop}\label{n3.3}If there is a derivation of $\n$ such $\delta>0$ then $\n$ is nilpotent and we have:
\begin{enumerate}
	\item $\mathrm{A}$ is nonzero and contains a derivation with strict positive spectrum.
	\item Any complete Lie algebra with the radical nilpotent $\n$ is given by $\g_1=\mathrm{R}_1\ltimes\n$ where $\mathrm{R}_1=\mathrm{S}_1\oplus\mathrm{T}_1\subset \mathrm{R}$ such that $\mathrm{S}_1$ is a semisimple subalgebra of $\mathrm{S}$ and $\mathrm{T}_1$ is a maximal torus in $\mathrm{Der}(\n)^{\mathrm{S}_1}=\mathrm{S}^{\mathrm{S}_1}+\mathrm{A}+ \mathrm{N}^{\mathrm{S}_1}$ such that $\mathrm{N}^{\mathrm{R}_1}=0$. There is a solvable Lie algebra of this type $\mathrm{T}\ltimes\n$ if $\mathrm{N}^{\mathrm{T}}=0$ for a maximal torus $\mathrm{T}$ on $\n$.
	\item Each Lie algebra $\g$ in $2.$ satisfies the hypotheses of the theorem of reduction.
\end{enumerate}
\end{prop}
\begin{proof} The nilpotence of $\n$ is classical. 1. The semisimple part in the Jordan's decomposition of a derivation $\delta$ with strict positive spectrum is still a derivation with positive spectrum. It follows from Theorem of Mostow that we may suppose that it is contained in $\mathrm{S}\oplus\mathrm{A}$. Up to a conjugation it admits a decomposition $\delta_1+\delta_0$ with $\delta_0\in\mathrm{S}$ and $\delta_1\in\mathrm{A}$. If $\n=\oplus_i\n_i$ is the decomposition in eigenvector spaces of $\delta_1$ with $\lambda_i$ the eigenvalues on $\n_i$, the fact $\delta_1$ commutes with $\mathrm{S}$ implies that $\mathrm{S}.\n_i\subset\n_i$ for all $i$ and the trace of $\delta_0\in\mathrm{S}$ is equal to zero on each $\n_i$. Then the trace over each $\n_i$ verifies:
$$\mathrm{tr}_{\n_i}(\delta_0)+\mathrm{tr}_{\n_i}(\delta_1)=\lambda_i\dim\n_i=\mathrm{tr}_{\n_i}(\delta)>0,$$
so $\lambda_i\in\Q^+-\left\{0\right\}$ and in particular $\delta_1\neq 0$ and $\mathrm{A}\neq 0$.\\
2. Let $\g_1$ denote the Lie algebra $\mathrm{R}_1\ltimes\n$ with $\mathrm{R}_1$ a subalgebra of $\mathrm{R}$. We suppose that $\n$ possesses a derivation $\delta$ with strict positive spectrum. It follows from $1.$ that there is a derivation $\delta_1$ with strict positive spectrum in $\mathrm{A}$. Since the algebra $\mathrm{A}$ is contained in the subalgebra $\mathrm{Der}(\n)^{\mathrm{R}_1}$ then each element $l$ of $\mathrm{A}$ can be extended to element of $\mathrm{Der}(\g_1)$ setting $l(\mathrm{R}_1)=0$. If $\g_1$ is complete we have $\mathrm{A}\subset\mathrm{R}_1$ and the subspace $\m=\mathrm{R}_1\cdot\n$ is equal to $\overline{\n}$ and so to $\n$, see (\ref{e4.1}), (\ref{e4.2}) and (\ref{E6.4}). We may write
$$(\mathrm{Der}(\overline{\n}))^{\mathrm{S_1}}=\mathrm{Der}(\n)^{\mathrm{S_1}}=\mathrm{S}^{\mathrm{S}_1}+\mathrm{A}+ \mathrm{N}^{\mathrm{S}_1},$$
with $\mathrm{R}_1=\mathrm{S}_1\oplus\mathrm{T}_1$. It follows in the case $\mathfrak{m}=\n$ that $\g_1\cong\mathrm{S}_1\oplus\mathrm{T}_1\ltimes\n$ where $\mathrm{T}_1$ is a maximal torus of $\mathrm{Der}(\n)^{\mathrm{S}_1}$ and $\mathrm{N}^{\mathrm{R}_1}=0$ since $\delta_1>0$ and $\delta_1\in\mathrm{R}_1$.\\
3. The hypothesis $\delta>0$ shows that for all complete Lie algebra $\g_1=\mathrm{R}_1\ltimes\n$ there is $\delta_1>0$ in $\mathrm{A}$. Hence $\mathrm{H}_k(\n)^{\delta_1}=\mathrm{H}_k(\n)^{\mathrm{A}}=0$ for $k\geq 1$ and so $\mathrm{H}^k(\n,\g_1/\n)^{\mathrm{R}_1}=0$ for $k\geq 1$. It follows from the long sequence of cohomology that 
\begin{equation}
	\mathrm{H}^{k+1}(\n,\n)^{\mathrm{R}_1}=\mathrm{H}^{k+1}(\n,\g_1)^{\mathrm{R}_1}\quad k\geq 1,
\end{equation}
and in particular $\g_1$ satisfies the hypotheses of the theorem of reduction.
\end{proof}\\
\begin{rem}
\begin{enumerate}
	\item The complete Lie algebras up to dimension $8$, cf \cite{C5}, satisfy the Proposition \ref{n3.3}.
	\item We can apply Proposition \ref{n3.3} in the following cases:\\
if there is a derivation $\delta>0$ of $\n$ such that $\mathrm{N}^\delta=0$, then for all semisimple subalgebra $\mathrm{S}_1$ of $\mathrm{S}$ and maximal torus $\mathrm{T}_1$ in $(\mathrm{Der}(\n))^{\mathrm{S}_1}$, the Lie algebra $\mathrm{S}_1\oplus\mathrm{T}_1\ltimes\n$ is complete and satisfies the property $3$ of Proposition \ref{n3.3}.
\end{enumerate}
\end{rem}
%\begin{de}A nilpotent Lie algebra $\n$ is called metabelien if $\left[\n,\left[\n,\n\right]\right]=0$.
%\end{de}
%For more informations about metabelian Lie algebras, see \cite{L}. Set $\mathfrak{w}=\left[\n,\n\right]$ and $\mathfrak{v}$ a complement of $\mathfrak{w}$ in $\n$. 
%\begin{lem}\label{L5.3}If $\n$ is a metabelian Lie algebra, then the map $\epsilon:\n\rightarrow\n$ defined by $\epsilon(v)=1$ for all $v\in\mathfrak{v}$ and $\epsilon(w)=2$ for all $w\in\mathfrak{w}$, is a derivation of $\n$. 
%\end{lem}

%\begin{lem}If $\n$ is a metabelian Lie algebra then  $$\Der(\n)=(\Der(\n))^{\epsilon}\oplus\mathrm{i}(\mathrm{Hom}_\K(\mathfrak{v},\mathfrak{w})),$$
%where $\mathrm{i}(f)|_\mathfrak{v}=f|_\mathfrak{v}$ and $\mathrm{i}(f)|_\mathfrak{w}=0$, for all $f\in\mathrm{Hom}_\K(\mathfrak{v},\mathfrak{w})$. 
%\end{lem}
%It follows from Lemma \ref{L5.3} that a metabelian Lie algebra then belongs to $\mathcal{C}^\infty_{N}$.  
%\begin{prop}Let $\n$ be a metabelian Lie algebra and $\g$ its associated complete Lie algebra. Then $\g$ is isomorphic to $\mathrm{T}\ltimes\mathrm{N}^\mathrm{T}\oplus\mathfrak{m}$, where $\mathfrak{m}=\mathrm{T}\cdot\n$, $\mathrm{N}$ is the largest nilpotency ideal of $\Der(\n)$, and $\mathrm{T}$ is a maximal torus on $\overline{\n}$ containing $\epsilon$. If $\mathrm{N}^{\mathrm{T}}\neq 0$ then $\g$ is not rigid.
%\end{prop}
%\begin{proof} If $\mathrm{N}^{\mathrm{T}}\neq 0$ then $\n(\g)\neq\n$, i.e., $\left[\g,\n\right]\neq\n$, we deduce the result from Theorem \ref{N3.2}
%\end{proof}
\section{Study of Schemes $\L_n^{\mathrm{T}}$ of nilpotent laws}
The study of rigid complete solvable Lie algebras in \cite{C2} had led quite naturally to that a sequence of algebraic solvable Lie algebras $\g^{(n)}:=\mathrm{T}\ltimes\varphi_n$ where the passage from $\varphi_n$ to $\varphi_{n+1}$ is made by central extension, the maximal torus $\mathrm{T}$ on $\varphi_n$ being extended to $\varphi_{n+1}$ with an additional weight $\alpha_{n+1}$ and $n$ is the dimension of $\varphi_n$. We denote again by $\mathrm{T}$ this new torus. 
One will develop next this construction under the hypotheses of simple path of weights $(\alpha_1,...,\alpha_n,...)$ related to $\mathrm{T}$; this defines a sequence of schemes $\L_n^{\mathrm{T}}$ for $n\geq n_0$. Under the hypotheses chosen on the weights, the local study of the Lie algebras $\g^{(n)}$ will profit from a double method:\\
1) the theorem of reduction which is applied here and reduces the local study of $\g^{(n)}$ to that of $\varphi_n$ in $\L_n^{\mathrm{T}}$;\\
2) the local study of the laws $\varphi_n$ is made by induction on $n$ starting from an initialization $n_0$ where the torus appears.\\
This work proposed here, in an affine geometry context (affine schemes), has its parallel in projective geometry which expresses the Lie algebras defined by generators and relations, particularly adapted to the study of the nilpotent Lie algebras. This second point of view initiated by G. Favre for finite dimension \cite{Fa} and A. Fialowski for infinite dimension \cite{F} will be separately treated in a forthcoming work, \cite{CC}.
\subsection{Notation}
Let $\mathrm{Diag_n}\left(\K\right)$ denote the vector space of diagonal matrices $d=\mathrm{diag}\left(a_1,...,a_n\right)$ defined by $d\cdot e_i=a_ie_i$ $\left(1\leq i\leq n\right)$. Since $\K$ being algebraically closed, then each torus on $\K^n$ is (up conjugation) a vector subspace $\mathrm{T}$ of $\mathrm{Diag_n}\left(\K\right)$. The torus $\mathrm{T}$ is defined by its weights $\alpha_i\in\mathrm{T}^*$ which are the restrictions of the dual linear forms $(e^i_i)^*$ on $\mathrm{T}$ where $(e^i_i)_{1\leq i\leq n}$ is the canonical basis of $\mathrm{Diag_n}\left(\K\right)$ satisfying $e^i_i(e_k)=\delta_{ik}e_k$ $\left(1\leq i\leq n\right)$. Thus 
$$t\cdot e_i=\alpha_i(t)e_i, \quad\left(1\leq i\leq n\right), \quad t\in\mathrm{T}.$$
The family $\mathrm{P}_n=\left(\alpha_i\right)_{1\leq i\leq n}$ linearly generates $\mathrm{T}^*$ and we have $s:=n-\dim\mathrm{T}$ independent linear equations linking the elements $\alpha_i$
\begin{equation}\label{e41}
	\sum^n_{i=1}\lambda_k^i\alpha_i=0\quad 1\leq k\leq s,
\end{equation}
where $\lambda^i_k\in\K$. The relations (\ref{e41}) define the torus $\mathrm{T}$ as the intersection of the kernels of $\sum^n_{i=1}\lambda_k^i (e^i_i)^*$ in $\mathrm{Diag_n}\left(\K\right)$ where $1\leq k\leq s$ .
If $\mathrm{T}$ is a torus on $\K^n$, we denote again by $\mathrm{T}$ the torus on $\K^{n+1}$ defined by adding the weight $\alpha_{n+1}\in\mathrm{T}^*$:
$$t\cdot e_i=\alpha_i(t)e_i, \quad\left(1\leq i\leq n\right),\quad t\cdot e_{n+1}=\alpha_{n+1}(t)e_{n+1},\quad t\in\mathrm{T},$$
where $(e_1,...,e_n,e_{n+1})$ is the canonical basis of $\K^{n+1}$.
The scheme $\L_n^{\mathrm{T}}$ is the set of Lie multiplications $\varphi$ defined by $\varphi\left(e_i,e_j\right)=\sum_{k=1}^n\varphi^k_{ij}e_k$ such that $t\cdot\varphi=0$, i.e.,
\begin{equation}\label{e42}
	(\alpha_k-\alpha_i-\alpha_j)(t)\varphi^k_{ij}=0\quad\forall t\in\mathrm{T}, i<j,k;
\end{equation}
this is equivalent to $\varphi^k_{ij}=0$ for all $i<j,k$ such that $\alpha_k\neq\alpha_i+\alpha_j$.\\
Denote by $\mathcal{J}$ the set the multiindices $\left(^k_{ij}\right)$ such that $i<j$ and $\alpha_i+\alpha_j=\alpha_k$. The coordinates in $\L_n^{\mathrm{T}}$ are variables $X^k_{ij}$ indexed by $\mathcal{J}$, and the Jacobi's relations are given by
\begin{equation}\label{e43}
	\mathrm{J}^h_{ijk}=\oint_{(ijk)}\sum_{l=1}^nX^l_{ij}X^h_{lk}=0,\quad 1\leq i<j<k\leq n,\quad 1\leq h\leq n,
\end{equation}
for $\alpha_i+\alpha_j+\alpha_k=\alpha_h\in\mathrm{P}_n.$ \\
Let $\sum_n\left(\mathrm{T}\right)$ be the subset of $\L_n^{\mathrm{T}}\left(\K\right)$ consisting of laws such that $\mathrm{T}$ is exactly a maximal torus of derivations. If we choose $\mathrm{T}$ such that $\sum_n\left(\mathrm{T}\right)$ is nonempty, then $\mathrm{T}$ is algebraic on $\K$, i.e., the relations (\ref{e41}) satisfied by the weights $\alpha_i$ are defined over $\Q$.
\begin{ex}
\begin{enumerate}
	\item $\sum_n\left(\mathrm{0}\right)$ is the set of characteristically nilpotent Lie multiplications on $\K^n$ (i.e., all derivations are nilpotent). It is nonempty for $n\geq 7$ and not Zariski open for $n=7$, cf \cite{C}.
	\item If $\mathrm{T}_n=\mathrm{Diag}_n\left(\K\right)$ then $\L^{\mathrm{T}_n}_n=\sum_n\left(\mathrm{T}_n\right)=\left\{0\right\}$.
\end{enumerate}
\end{ex}
\begin{rem} The set $\sum_n\left(\mathrm{T}\right)$ is constructible for the Zariski topology, (i.e., it is a finite union of local closed sets) and stable under the canonical group $\mathrm{G}_n:=\mathrm{Gl}_n\left(\K\right)^{\mathrm{T}}_0$.
\end{rem}
\begin{prop}\label{p5.1}The isomorphic classes in $\sum_n(\mathrm{T})$ are the orbits of the normalizer group $\mathrm{H}$ of $\mathrm{T}$ in $\Gl_n(\K)$ under the canonical action.
\end{prop}
\begin{proof} We can see that $\mathrm{H}$ stabilizes $\L_n^{\mathrm{T}}(\K)$ and $\sum_n(\mathrm{T})$. Conversely, if $\varphi_1$ and $\varphi_2$ are isomorphic, then there is $s\in\Gl_n(\K)$ such that $s\ast\varphi_1$ is equal to $\varphi_2$ with maximal torus $s\cdot\mathrm{T}\cdot s^{-1}$. The field $\K$ being algebraically closed, there is an automorphism $s'$ of $\varphi_2$ which conjugates $s\cdot\mathrm{T}\cdot s^{-1}$ and $\mathrm{T}$ according to Mostow's theorem; we have $s'\cdot s\in\mathrm{H}$.
\end{proof}
\subsection{Induction on schemes: $\L_n^{\mathrm{T}}\rightarrow\L_{n+1}^{\mathrm{T}}$}
Let $\mathrm{T}$ be a torus on $\K^{n+1}$ as the above section and its set of weights is given by $\mathrm{P_{n+1}}=\mathrm{P_{n}}\cup\left\{\alpha_{n+1}\right\}$. We suppose the choice of $\alpha_{n+1}$ such that $\sum_{n+1}\left(\mathrm{T}\right)$ is nonempty. The weight of nonzero weightvector of the form $\left[e_i,e_{n+1}\right]$ for $i\leq n$, is $\alpha_i+\alpha_{n+1}\in\mathrm{P_{n+1}}$. If $0\notin\mathrm{P_{n}}$ and $\alpha_{n+1}\notin\mathrm{P_{n}}-\mathrm{P_{n}}$ then $\left[e_i,e_{n+1}\right]=0$ for $i\leq n$ and $\K e_{n+1}$ is central. This means that any law $\varphi_{n+1}\in\L^{\mathrm{T}}_{n+1}\left(\K\right)$ is a central extension of a law $\varphi_n\in\L^{\mathrm{T}}_{n}\left(\K\right)$ of kernel $\K e_{n+1}$, i.e., we have the exact sequence of Lie algebras
\begin{equation}\label{e44}
	\left\{0\right\}\longrightarrow\K e_{n+1}\longrightarrow\left(\K^{n+1},\varphi_{n+1}\right)\longrightarrow\left(\K^{n},\varphi_{n}\right).
\end{equation}
If $\mathrm{T}$ is maximal on $\left(\K^{n+1},\varphi_{n+1}\right)$, i.e.,  $\varphi_{n+1}\in\sum_{n+1}\left(\mathrm{T}\right)$, then the extension is not trivial. This extension corresponds to a nonzero class of $\mathrm{H}^2\left(\varphi_n,\K e_{n+1}\right)^{\mathrm{T}}$, for the trivial action. The group of homology $\mathrm{H}_2\left(\varphi_n\right)$ is a $\mathrm{T}$-module and may be decomposed as
\begin{equation}\label{e45}
\mathrm{H}_2\left(\varphi_n\right)=\bigoplus_{\alpha\in\mathrm{T}^*}\mathrm{H}_2\left(\varphi_n\right)_{\alpha}.
\end{equation}
Thus
\begin{equation}\label{e46}
	\mathrm{H}^2\left(\varphi_n,\K e_{n+1}\right)^{\mathrm{T}}=\mathrm{Hom}_{\mathrm{T}}\left(\mathrm{H}_2\left(\varphi_n\right),\K e_{n+1}\right)\cong\left(\mathrm{H}_2(\varphi_n\right)_{\alpha_{n+1}})^*.
\end{equation}
It follows from (\ref{e46}) that a nonzero cohomology class  corresponds to a nonzero homology class such that the weight $\alpha_{n+1}$ appears in the decomposition (\ref{e45}) and implies that 
\begin{equation}\label{e47}
	\alpha_{n+1}=\alpha_i+\alpha_j\quad \mathrm{for}\quad i<j\leq n.
\end{equation}
We suppose that $\mathrm{P}_n$ consists of strict positive weights, i.e., there is $t\in\mathrm{T}$ such that each weight  $\alpha_i$ satisfies $\alpha_i(t)\in\Q^{*+}$ for $1\leq i\leq n$. We deduce the following properties:
\begin{enumerate}
		\item $\L_n^{T}\left(\K\right)$ is formed of nilpotent Lie multiplications,
		\item $\mathrm{P}_{n+1}$ is consisted of strict positive weights by (\ref{e47}),
		\item A complete law $\mathrm{T}\ltimes\varphi_n$ where $\varphi_n\in\L^{\mathrm{T}}_n$ satisfies the reduction theorem. 
\end{enumerate}\
We will introduce the notion of path of weights which permits to find the previous properties like that of central extension.
\begin{de}A sequence of weights, $\mathrm{Q}=\left(\alpha_p\right)_{p\in\N}\subset\mathrm{T}^*$, is said to be a path of weights if
\begin{enumerate}
	\item there is $n_0\in\N^*$ such that the family $\left(\alpha_1,...,\alpha_{n_0}\right)$ generates $\mathrm{T}^*$ over $\K$ and $\alpha_i>0$ for all $i=1,...,n_0$,
	\item $\sum_{n}\left(\mathrm{T}\right)$ is a nonempty set for all $n\geq n_0$,
	\item $\alpha_{n+1}$ does not belong to $(\mathrm{P_n}-\mathrm{P_n})$ for all $n\geq n_0$.
\end{enumerate}
$n_0$ is called the initialization of $\mathrm{Q}$.
\end{de}
The initialization $n_0$ of $\mathrm{Q}$ corresponds to the value of $n$ where the torus $\mathrm{T}$ appears.\\ 
The condition of the central extension will give the existence of a morphism of schemes which generalizes the quotient map (\ref{e44}). Let $\mathrm{A}_n$ denote the polynomial ring $\K[X^k_{ij}:\binom{k}{ij}\in\mathcal{J}]$ and $\mathrm{J}_n$ the ideal generated by the Jacobi's polynomials $\mathrm{J}^h_{ijk}$ with $1\leq i<j<k,h\leq n$, see (\ref{e43}). 
\begin{prop}\label{n41}If $\mathrm{Q}=\left(\alpha_p\right)_{p\in\N}\subset\mathrm{T}^*$ is a path of weights with $n_0$ initialization, then 
\begin{enumerate}
	\item the ideal $\mathrm{J}_{n+1}$ is generated by the ideal $\mathrm{J}_n$ and the polynomials $\mathrm{J}^{n+1}_{ijk}$ where $1\leq i<j<k\leq n$;
	\item the canonical ring monomorphism $\mathrm{i}_n:\mathrm{A}_n\rightarrow\mathrm{A}_{n+1}$ induces a scheme morphism $$\pi_{n+1}:\L^\mathrm{T}_{n+1}\rightarrow\L^\mathrm{T}_{n}$$ 
defined on the space of rational points $\L^{\mathrm{T}}_{n+1}(\K)$ by the quotient map of (\ref{e44}).
\end{enumerate}
\end{prop}
\begin{proof} 1. Consider $\mathrm{J}^h_{ijk}$ a nontrivial Jacobi's polynomial where $h\leq n$, we will show that it belongs to $\mathrm{J}_n$. It is clear for $k\leq n$. If $k=n+1$ then $\alpha_i+\alpha_j+\alpha_{n+1}\in\mathrm{P}_{n+1}$. It follows by assumptions that one of the following elements $\alpha_{n+1}+\alpha_i$, $\alpha_{n+1}+\alpha_j$, $\alpha_{n+1}+\alpha_i+\alpha_j$ belongs to $\mathrm{P}_{n+1}$ with $\alpha_i+\alpha_j\in\mathrm{P}_{n+1}$. However  if  $\alpha_{i}+\alpha_j\in\mathrm{P}_n$ then the elements $\alpha_{n+1}+\alpha_i$, $\alpha_{n+1}+\alpha_j$ and $\alpha_{n+1}+\alpha_i+\alpha_j$  do not belong to $\mathrm{P}_{n+1}$ by hypothesis of $\mathrm{P}_{n+1}$. If $\alpha_{i}+\alpha_j\notin\mathrm{P}_n$ then $\alpha_{i}+\alpha_j=\alpha_{n+1}$ and its multiplicity is $1$; The weight $2\alpha_{n+1}$ corresponds to a trivial bracket. By  similar arguments we can show that for $i<j<k$ the polynomial $\mathrm{J}^{n+1}_{ijk}$ satisfies that $k\leq n$. Hence we deduce that $\mathrm{J}_{n+1}$ is generated by $\mathrm{J}_n$ and the polynomials $\mathrm{J}^{n+1}_{ijk}$. \\
2. The correspondence $X^k_{ij}\rightarrow X^k_{ij}$ for $i<j\leq n$ and $k\leq n$ induces an injective map $\mathrm{i}_n:\mathrm{A}_n\rightarrow\mathrm{A}_{n+1}$ which sends $\mathrm{J}_n$ to $\mathrm{J}_{n+1}$. Then it induces a quotient morphism $\overline{\mathrm{i}}_n:\mathrm{A}_{n}/\mathrm{J}_{n}\rightarrow \mathrm{A}_{n+}/\mathrm{J}_{n+1}$ and $\pi_{n+1}:=\mathrm{Spec}(\overline{\mathrm{i}}_n):\L^\mathrm{T}_{n+1}\rightarrow\L^\mathrm{T}_{n}$.
\end{proof}\\
Under the hypotheses of Proposition \ref{n41}, there is a scheme morphism induced from the quotient ring morphism $\mathrm{pr}_n:\mathrm{A}_{n+1}\rightarrow\mathrm{A}_n$ by the ideal generated by the new coordinates $X^{n+1}_{ij}$. We have $\mathrm{pr}_n(\mathrm{J}_{n+1})=\mathrm{J}_n$, it follows that there is a quotient morphism $\mathrm{\overline{pr}}_n:\mathrm{A}_{n+1}/\mathrm{J}_{n+1}\rightarrow \mathrm{A}_n/\mathrm{J}_{n}$. We deduce a morphism of scheme
\begin{equation}
	\mathrm{i}_n=\mathrm{Spec}(\mathrm{\overline{pr}}_n):\L_{n}^{\mathrm{T}}\rightarrow\L_{n+1}^{\mathrm{T}}
\end{equation}
which sends a law $\mathrm{\varphi}_n\in\L^{\mathrm{T}}_n$ to $(\varphi_n,\K^n)\times\K e_{n+1}$ the direct product of $(\varphi_n,\K^n)$ with the abelian ideal $\K$. We can see that $\mathrm{i}_n(\L_{n}^{\mathrm{T}})\cap\sum_{n+1}(\mathrm{T})=\emptyset$.
\subsection{Continuous families of Lie algebras}
\begin{de} \label{d5.2}A path of weights is said to be simple if all weights are distinct.
\end{de}
From Definition \ref{d5.2}, the coordinates $X^k_{ij}$ may be indexed by the weights themselves. One then can write $X_{ij}$ instead of $X^k_{ij}$ since the index $k$ is fixed by the weight $\alpha_i+\alpha_j$. The set of pairs $i<j$ such that $\alpha_i+\alpha_j\in\mathrm{P}_n$ will be denoted again by $\mathcal{J}$.

\begin{prop}Let $\mathrm{T}$ be a torus on $\K^n$. One supposes that its set of weights $\mathrm{P}_n$ is a simple path of weights. Then $\sum_{n}\left(\mathrm{T}\right)$ is a Zariski open set which is equal to the set of elements $\varphi_n\in\L_n^{\mathrm{T}}$ such that $\mathrm{T}\ltimes\varphi_n$ is complete. It is the union of the $\G_n$-orbits of maximal dimesion, $n-\dim_\K\mathrm{T}$.
\end{prop}
\begin{proof} If $\varphi$ is an element of $\L_n^{\mathrm{T}}$, then $\Der\left(\varphi\right)^\mathrm{T}$ is a torus $\tau$ which is contained in $\mathrm{Diag}_n\left(\K\right)$. Then the law $\varphi$ belongs to $\mathcal{C}^\infty_N$ and $\mathrm{\tau}\ltimes\varphi$ is a complete Lie algebra, see Section 5. Then the set $\sum_{n}\left(\mathrm{T}\right)$ consists of elements $\varphi$ such that $\tau=\mathrm{T}$, i.e., $\dim\tau=\dim\mathrm{T}$ since $\mathrm{T}\subset\tau$. It follows that it is the Zariski open set consisting of orbits of maximal dimension.
\end{proof}
\begin{rem} Under the above hypotheses, any element $\varphi\in\L^{\mathrm{T}}_n$ belongs to the Zariski open $\sum_{n}\left(\tau\right)$ of the subscheme $\L^{\tau}_n$ of $\L^{\mathrm{T}}_n$ in which we can apply the reduction theorem for $\tau\ltimes\varphi$. The different $\Sigma_n(\mathrm{\tau})$ form a stratification of $\L_n^{\mathrm{T}}$.
\end{rem}
The biggest stratum $\Sigma_n(\mathrm{T})$ gives the quotient variety $\Sigma_n(\mathrm{T})/\G_n$ under the action of $\G_n$.
\begin{prop}The isomorphic classes in the variety $\Sigma_n(\mathrm{T})/\G_n$ are the orbits of the finite group $\Gamma:=\mathrm{H}/\mathrm{H}_0$, with $\mathrm{H}$ the normalizer group of $\mathrm{T}$ in $\Gl_n(\K)$ and $\mathrm{H}_0$ its identity component. This variety is called continuous family associated with the maximal torus $\mathrm{T}$.
\end{prop}
\begin{proof} It is a direct consequence of Propostion \ref{p5.1}.
\end{proof}\\
%\textbf{Group\quad$\Gamma$.}\\
Next, we will study the quotient variety $\Sigma_n(\mathrm{T})/\G_n$ with the help of slices $\L^{\mathrm{T},\mathcal{A}}_{n,\varphi_0}$, $\varphi_0\in\Sigma_n(\mathrm{T})$, which are local affine charts.\\
The weights being distinct, then $\mathrm{G}_n$ is the diagonal group identified with $(\K^*)^n$. The canonical action of an element $s=\left(s_1,...,s_n\right)$ of $\mathrm{G}_n$ on $X$ is defined by
\begin{equation}
	Z^k_{ij}=(s\ast X)^k_{ij}=\frac{s_k}{s_is_j}X^k_{ij},
\end{equation}
where $X$ is a law defined by its coordinates $X(e_i,e_j)=\sum_kX^k_{ij}e_k$. It is particularly easy to characterize an admissible set of $\mathcal{J}$ at a law $\varphi$. It is just a subset $\mathcal{A}$ of $\mathcal{J}$ such that the following system of equations 
\begin{equation}
	\varphi^k_{ij}=\frac{s_k}{s_is_j}\varphi^k_{ij}\quad \mathrm{for}\quad ({i<j})\in\mathcal{A},
\end{equation}
is equivalent to the system of equations
\begin{equation}
	s\ast\varphi=\varphi,
\end{equation}
and which is minimal for this property. 
The minimality of $\mathcal{A}$ implies that $\varphi_{ij}^k\neq 0$ for all $({i<j})\in\mathcal{A}$. It follows that $\mathcal{A}$ is contained in the set $\mathcal{J}_{\varphi}$ consisting of pairs $(i<j)$ such that $\varphi_{ij}^k\neq 0$.
\begin{prop}An subset $\mathcal{A}$ of $\mathcal{J}$ is an admissible set at $\varphi$ if and only if $\mathcal{A}$ is a minimal set of $\mathcal{J}_{\varphi}$ such that the following system
\begin{equation}\label{e53}
	s_k=s_is_j,\quad (i<j)\in\mathcal{A},
\end{equation}
defines the subgroup $\mathrm{Aut(\varphi)}_0^{\mathrm{T}}$ of $\mathrm{G}_n$.
\end{prop}  
The laws $\varphi$ of $\L_n^{\mathrm{T}}(\K)$ which have the same subset $\mathcal{J}_{\varphi}\subset\mathcal{J}$ thus have the same admissible sets $\mathcal{A}\subset\mathcal{J}_{\varphi}$, the same group $\mathrm{Aut(\varphi)}_0^{\mathrm{T}}$ and the same maximal torus $\tau$. 
Two equivalent laws under  $\mathrm{G}_n$ have the same subset $\mathcal{J}_{\varphi}$ of $\mathcal{J}$ since Eq (\ref{e43}) implies that a nonzero coordinate remains nonzero under $\mathrm{G}_n$.
\begin{thm}\label{n51}Under the hypothesis of a simple path of weights. Each law $\varphi$ of $\L_{n,\varphi_0}^{\mathrm{T},\mathcal{A}}(\K)$ admits an admissible set which contains $\mathcal{A}$. If the cardinal of $\mathcal{A}$ is maximal (i.e., $|\mathcal{A}|=n-\dim_\K\mathrm{T}$) then 
\begin{enumerate}
	\item $\varphi$ admits $\mathcal{A}$ as admissible set iff $\varphi^k_{ij}\neq 0$ for all $(^k_{ij})\in\mathcal{A}$ where $i<j$.
	\item All laws of $\L_{n,\varphi_0}^{\mathrm{T},\mathcal{A}}(\K)$ admit $\mathcal{A}$ as admissible set.
	\item $\L_{n,\varphi_0}^{\mathrm{T},\mathcal{A}}(\K)$ is contained in $\sum_{n}\left(\mathrm{T}\right)$ and its isomorphic classes are the traces of the $\Gamma$-orbits.
	\item $\varphi$ admits $\mathcal{A}$ as admissible set iff there is an element $s\in\mathrm{G}_n$ such that $s\ast\varphi\in\L_{n,\varphi_0}^{\mathrm{T},\mathcal{A}}(\K)$.
\end{enumerate}
\end{thm}
\begin{proof} The set of the laws $\varphi$ which admits $\mathcal{A}$ as admissible set satisfy $\varphi^\alpha\neq 0$ for all $\alpha\in\mathcal{A}$. This condition is sufficient if the cardinal of $\mathcal{A}$ is maximal. This proves the statements $1,2$ and $3$. 
The orbit $\left[\varphi_0\right]$ of $\varphi_0$ under $(\K^*)^n$ is given by
$$\left[\varphi_0\right]=\left\{(\frac{s_k}{s_is_j}(\varphi_0)_{ij}^k):s=(s_1,...,s_n)\in\mathrm{G}_n\right\}.$$
The components indexed by an admissible set $\mathcal{A}$ are non null and the corresponding terms $(\frac{s_k}{s_is_j})$ are free and can take arbitrary nonnull values; there thus  exits in the orbit of $\varphi_0$ a representative $\varphi$ for any arbitrary choice of non null components $(\varphi^\alpha)_{\alpha\in\mathcal{A}}.$ If $|\mathcal{A}|$ is maximal, we obtain $4$.
\end{proof}\\
\textbf{Consequence.} If the cardinal of $\mathcal{A}$ is maximal, the subschemes defined in Theorem \ref{n51} are such  that the laws $\varphi$ play a symmetric role between them since they admit $\mathcal{A}$ as admissible set and have the same values $\varphi^\alpha$ for all $\alpha\in\mathcal{A}$. By definition of an admissible set, such a subscheme contains only a finite number of points of a same orbit and thus represents a continuous family of Lie structures.\\ 
\textbf{Convention.} Let $\mathcal{A}$ be a maximal admissible set at $\varphi_1$ and $\varphi_2$. It follows from Theorem \ref{n51}, 4. that the both schemes $\L_{n,\varphi_1}^{\mathrm{T},\mathcal{A}}$ and $\L_{n,\varphi_2}^{\mathrm{T},\mathcal{A}}$ are conjugated by $s\in\mathrm{G}_n$ and we may fix arbitrary values of the components $(\varphi^{\mathcal{A}})$ of a law $\varphi$ in $\K$. One adopts the convention below: $X^k_{ij}=1$ for all $(^k_{ij})\in\mathcal{A}$. The corresponding subscheme will be denoted by $\L_n^{\mathrm{T},\mathcal{A}}$ which is isomorphic to $\L_{n,\varphi_0}^{\mathrm{T},\mathcal{A}}$. \\
This also gives a convention for a bracket of a nilpotent Lie algebra verifying the hypotheses above which consists to set $\varphi^k_{ij}=1$ for all $(i<j)\in\mathcal{A}$ for all admissible set. 
\subsection{Study of direct filiations $\mathcal{A}_n\hookrightarrow\mathcal{A}_{n+1}$}
The laws of the open set $\sum_{n+1}(\mathrm{T})$ dot not always come from $\sum_{n}(\mathrm{T})$ and it also happens there are laws $\varphi_n$ of $\sum_{n}(\mathrm{T})$ do not have an extension in $\sum_{n+1}(\mathrm{T})$; this is verified when $\mathrm{H}_2(\varphi_n)_{\alpha_{n+1}}$ is zero. However an extension $\varphi_{n+1}$ of $\varphi_n\in\sum_{n}(\mathrm{T})$, cf (\ref{e44}), belongs to $\sum_{n+1}(\mathrm{T})$ iff it is not trivial. We will say that $\varphi_{n+1}$ is obtained by a direct filiation of $\varphi_n$. This means that there are two admissible sets $\mathcal{A}_n$ and $\mathcal{A}_{n+1}$ associated to $\varphi_{n+1}$ and $\varphi_n$ respectively such that $$|\mathcal{A}_{n+1}|=n+1-r,\quad|\mathcal{A}_{n}|=n-r\quad\mathrm{and}\quad \mathcal{A}_{n+1}=\mathcal{A}_{n}\cup\left\{(^{n+1}_{pq})\right\}$$ 
for some $p<q$ satisfying $(\varphi_{n+1})^{n+1}_{pq}\neq 0$ where $r=\dim\mathrm{T}$. \\
The inclusion $\mathcal{A}_{n}\hookrightarrow\mathcal{A}_{n+1}$ thus induces a process of construction by extension (\ref{e44}), denoted by 
$$\L_n^{\mathrm{T},\mathcal{A}_n}\longrightarrow\L_{n+1}^{\mathrm{T},\mathcal{A}_{n+1}}.$$
It is said to be a direct filiation; one will develop it.
We consider $\pi^{-1}_{n+1}(\varphi_n)\cap\L_{n+1}^{\mathrm{T},\mathcal{A}_{n+1}}$ the trace of the fiber of $$\pi_{n+1}:\L_{n+1}^{\mathrm{T}}\rightarrow\L_n^{\mathrm{T}},$$ 
of each point $\varphi_n$ of $\L_n^{\mathrm{T},\mathcal{A}_n}$ on $\L_{n+1}^{\mathrm{T},\mathcal{A}_{n+1}}$. It is the set of laws $\varphi_{n+1}$ of $\L_{n+1}^{\mathrm{T},\mathcal{A}_{n+1}}$ which satisfy the equation (\ref{e44}) and the following condition $(\varphi_{n+1})^{n+1}_{pq}=1$. As subscheme, it is the set of $(X^{n+1}_{ij})$ which verify the following properties:
\begin{equation}\label{e51}
	\oint_{ijk}\sum_l(\varphi_n)^l_{ij} X_{lk}^{n+1}=0, 
\end{equation}
i.e. Jacobi's identity, and
\begin{equation}\label{e52}
	X_{pq}^{n+1}=1.
\end{equation}
The injection of tangent of Zariski of the subscheme $\pi^{-1}_{n+1}(\varphi_n)\cap\L_{n+1}^{\mathrm{T},\mathcal{A}_{n+1}}$ in that of $\L_{n+1}^{\mathrm{T},\mathcal{A}_{n+1}}$ at the point $\varphi_{n+1}$ is interpreted thanks to a morphism of a long exact sequence of cohomology related to (\ref{e44}).
\subsection{Study of a long exact sequence of cohomology} 
We use the notation of section 5. More generally, let us take $\n$ a nilpotent Lie algebra, and $\mathfrak{a}$ a central ideal contained in $\left[\n,\n\right]$. Let $\mathrm{R}$ be a completely reducible Lie subalgebra of $\mathrm{Der}(\n)$ which leaves stable $\mathfrak{a}$. We denote again by $\mathrm{R}$ the subalgebra of derivations on the quotient $\n/\mathfrak{a}$.
The exact sequence of $\n$-adjoint modules
\begin{equation}
	0\rightarrow\mathfrak{a}\rightarrow\n\rightarrow\n/\mathfrak{a}\rightarrow 0,
\end{equation}
gives an exact sequence of differential graded Chevalley-Eilenberg's complex 
\begin{equation}
0\rightarrow\mathrm{C}(\n,\mathfrak{a})\rightarrow\mathrm{C}(\n,\mathfrak{n})\rightarrow\mathrm{C}(\n,\n/\mathfrak{a})\rightarrow 0,
\end{equation}
and an exact sequence for the $\mathrm{R}$-invariant for the canonical action of $\mathrm{R}$ deduced from the action on $\n$:
\begin{equation}
0\rightarrow\mathrm{C}(\n,\mathfrak{a})^{\mathrm{R}}\stackrel{i}\rightarrow\mathrm{C}(\n,\mathfrak{n})^{\mathrm{R}}\stackrel{p}{\rightarrow}\mathrm{C}(\n,\n/\mathfrak{a})^{\mathrm{R}}\rightarrow 0.
\end{equation}
Then its long exact sequence of cohomology is given by 
\begin{equation}
0\rightarrow\mathrm{H}^0\left(\n,\mathfrak{a}\right)^\mathrm{R}\stackrel{\overline{i}_0}{\rightarrow}\mathrm{H}^0\left(\n,\n\right)^\mathrm{R} \stackrel{\overline{p}_0}{\rightarrow}\mathrm{H}^0\left(\n,\n/\mathfrak{a}\right)^\mathrm{R}\stackrel{\partial_0}{\rightarrow}\mathrm{H}^1\left(\n,\mathfrak{a}\right)^\mathrm{R}\cdots
\end{equation}
\begin{equation} \rightarrow\mathrm{H}^k\left(\n,\mathfrak{a}\right)^\mathrm{R}\stackrel{\overline{i}_k}{\rightarrow}\mathrm{H}^k\left(\n,\n\right)^\mathrm{R}\stackrel{\overline{p_k}}{\rightarrow}\mathrm{H}^k\left(\n,\n/\mathfrak{a}\right)^\mathrm{R}\stackrel{\partial_k}{\rightarrow}\mathrm{H}^{k+1}\left(\n,\mathfrak{a}\right)^\mathrm{R}\rightarrow
\end{equation}
and satisfies the following properties:
\begin{prop}\label{nn51}The homomorphism $\overline{i}_2$ verifies the following conditions
\begin{enumerate}
	\item if $\quad\mathrm{Hom}_{\mathrm{R}}(\n\wedge\mathfrak{a},\mathfrak{a})=\left\{0\right\}$, then the image of $\mathrm{H}^2(\n,\mathfrak{a})^{\mathrm{R}}$ by $\overline{i}_2$ is consisted of integrable classes in the scheme $\L_n^{\mathrm{R}}$, i.e., 
	$$\forall\xi\in\mathrm{Z}^2(\n,\mathfrak{a})^{\mathrm{R}}\Rightarrow\left[i_2(\xi),i_2(\xi)\right]=0;$$
	\item if $\mathrm{Hom}_{\mathrm{R}}(\mathfrak{a},\n/\mathfrak{a})=\left\{0\right\}$ and $\mathrm{H}^1(\n,\n)^{\mathrm{R}}\simeq\mathrm{H}^1(\n/\mathfrak{a},\n/\mathfrak{a})^{\mathrm{R}}$ then $\overline{i}_2$ is injective.
\end{enumerate}
\end{prop}
\begin{proof} 1. Let $\varphi_0$ denote the law of $\n$, for each $\xi\in\mathrm{Z}^2(\n,\mathfrak{a})^{\mathrm{R}}$ one associates
$$\psi:=i_2(\xi)\in\mathrm{Z}^2(\n,\n)^{\mathrm{R}}\quad\mathrm{and}\quad\varphi_0+t\psi\in\mathrm{C}^2(\n,\n)^\mathrm{R}.$$
Since 
$$\left[\varphi_0,\varphi_0\right]=\left[\varphi_0,\psi\right]=0$$ 
it follows that  
$$\left[\varphi_0+t\psi,\varphi_0+t\psi\right]=t^2\left[\psi,\psi\right].$$
Hence,
$$\psi(\psi(x,y),z)\in\psi(\mathfrak{a}\times\n)=\widetilde{\psi}(\mathfrak{a}\wedge\n)\subset\mathfrak{a},\quad\forall x,y,z\in\n,$$ 
where $\widetilde{\psi}:\wedge^2\n\rightarrow\mathfrak{a}$ is defined by $\widetilde{\psi}(x\wedge y)=\psi(x,y)$. \\
Since the restriction of $\widetilde{\psi}$ to $\n\wedge\mathfrak{a}$ belongs to $\mathrm{Hom}_{\mathrm{R}}(\n\wedge\mathfrak{a},\mathfrak{a})=\left\{0\right\}$, it follows that 
$$\psi(\psi(x,y),z)=0,\quad\forall x,y,z\in\n.$$
Consequently, $\left[\psi,\psi\right]=0$ and $\varphi_0+t\psi$ is a deformation, i.e. the class of $\psi$ is integrable.\\
2. The group $\mathrm{H}^1(\n,\mathfrak{a})^{\mathrm{R}}$ is isomorphic to $\mathrm{Hom}_{\mathrm{R}}(\n/\left[\n,\n\right],\mathfrak{a})$ and is zero since $\mathfrak{a}\subset\left[\n,\n\right]$ and $\mathrm{Hom}_{\mathrm{R}}(\mathfrak{a},\n/\mathfrak{a})=\left\{ 0\right\}$. Hence the morphism $\overline{p}_1$ is injective. Under the same hypothesis we can show that $\delta\in\mathrm{Z}^1(\n,\n/\mathfrak{a})^{\mathrm{R}}$ is zero on $\mathfrak{a}$ such that $\mathrm{Z}^1(\n,\n/\mathfrak{a})^{\mathrm{R}}$ may be identified with $\mathrm{Z}^1(\n/\mathfrak{a},\n/\mathfrak{a})^{\mathrm{R}}$ and $\mathrm{H}^1(\n,\n/\mathfrak{a})^{\mathrm{R}}$ with $\mathrm{H}^1(\n/\mathfrak{a},\n/\mathfrak{a})^{\mathrm{R}}$. Since $\mathrm{H}^1(\n/\mathfrak{a},\n/\mathfrak{a})^{\mathrm{R}}\simeq\mathrm{H}^1(\n,\n)^{\mathrm{R}}$ it follows that $\overline{p}_1$ is bijective and so $\overline{i}_2$ is injective.
\end{proof}\\
This Proposition can be applied under the hypotheses of a simple path of weights with $\n:=(\K^{n+1},\varphi_{n+1})$, $ \mathfrak{a}:=\K e_{n+1}$, $\n/\mathfrak{a}=(\K^{n},\varphi_{n})$ and a torus $\mathrm{T}$ for $\mathrm{R}$. We will compute $\mathrm{H}^2(\varphi_{n+1},\K e_{n+1})^{\mathrm{T}}$ using the proposition which follows.

More generally, let $\n$ be a nilpotent Lie algebra with bracket $\varphi$, and $\mathfrak{a}$ a central ideal of dimension one which is stable under a torus $\mathrm{T}$ on $\n$. Let $\mathfrak{b}$ be a complement direct of $\mathfrak{a}$ in $\n$ stable under $\mathrm{T}$, then it can be identified with $\n/\mathfrak{a}$ as $\mathrm{T}$-module. Hence $\varphi$ may be written as $\varphi_0+\psi$ where $\varphi_0$ is the Lie bracket defined on $\n/\mathfrak{a}$ transferred on $\mathfrak{b}$ and $\psi$ is a cocycle of $\mathrm{Z}^2(\n/\mathfrak{a},\mathfrak{a})^{\mathrm{T}}$. \\
Let $\Omega(\n)$ be the $\mathrm{T}$-submodule of $\wedge^3\n$ generated by the vectors $\oint_{(xyz)}\varphi(x,y)\wedge z$. \\
If $\beta$ is the weight of $\mathrm{T}$ on $\mathfrak{a}$ then the group $\mathrm{H}^2(\n,\mathfrak{a})^{\mathrm{T}}$ is isomorphic to 
$$\mathrm{Hom}_{\mathrm{T}}(\mathrm{H}_2(\n),\mathfrak{a})=(\mathrm{H}_2(\n)_{\beta})^*.$$
It follows from Proposition 3.2 in \cite{C2}, that $\mathrm{H}_2(\n)_{\beta}$ is isomorphic to 
\begin{equation}
	\frac{(\ker\widetilde{\varphi})_{\beta}}{\Omega(\n)_{\beta}}.
\end{equation}
It can be computed by the following
\begin{prop}\label{nn52} If $\mathrm{T}$ does not have a null weight on $\n/\mathfrak{a}$, then the vector space $\mathrm{H}_2(\n)_{\beta}$ is isomorphic to $$\frac{(\ker\widetilde{\varphi}_0)_{\beta}\cap(\ker\widetilde{\psi})_{\beta}}{\Omega(\n/\mathfrak{a})_{\beta}}.$$ 
If the class of $\psi$ is nonzero, then it is a hyperplan of $\mathrm{H}_2(\n/\mathfrak{a})_{\beta}$.
\end{prop}
\begin{proof} We have
$$
	\widetilde{\varphi}((b_1+a_1)\wedge (b_2+a_2))=\varphi(b_1+a_1,b_2+a_2)=\varphi(b_1,b_2)=\varphi_0(b_1,b_2)+\psi(b_1,b_2),$$

for all $ (b_1,b_2,a_1,a_2)\in\mathfrak{b}^2\times\mathfrak{a}^2$. It follows that
\begin{equation}\label{e5.1}
(\ker\widetilde{\varphi})_{\beta}=(\ker\widetilde{\varphi}_0)_{\beta}\cap(\ker\widetilde{\psi})_{\beta}\oplus (\n\wedge\mathfrak{a})_{\beta}.
\end{equation}
Hence $(\n\wedge\mathfrak{a})_{\beta}=0$ since $\mathfrak{a}\wedge\mathfrak{a}=0$ and $\mathfrak{b}\simeq\n/\mathfrak{a}$ does not admit a null weight by hypothesis. By Eq. (\ref{e5.1}), we have $(\ker\widetilde{\varphi})_{\beta}=(\ker\widetilde{\varphi}_0)_{\beta}\cap(\ker\widetilde{\psi})_{\beta}$.
The space $\Omega(\n)_{\beta}$ is generated by the cyclic sums of tensors,  we have
$$\varphi(b_1+a_1,b_2+a_2)\wedge(b_3+a_3)=\varphi_0(b_1,b_2)\wedge b_3+\varphi_0(b_1,b_2)\wedge a_3+\psi(b_1,b_2)\wedge b_3,$$
for all $(b_1,b_2,b_3,a_1,a_2,a_3)\in\mathfrak{b}^3\times\mathfrak{a}^3$.
According to the hypothesis, the last two tensors have a null projection on the weightspace of $\beta$ and hence there are isomorphisms
$$\Omega(\n)_{\beta}\simeq\Omega(\n/\mathfrak{a})_{\beta}\quad\mathrm{and}\quad \mathrm{H}_2(\n)_{\beta}\simeq((\ker\widetilde{\varphi}_0)_{\beta}\cap(\ker\widetilde{\psi})_{\beta})/\Omega(\n/\mathfrak{a})_{\beta}.$$
The homomorphism $\widetilde{\psi}:\wedge^2\n/\mathfrak{a}\rightarrow\mathfrak{a}$ is zero on a complement of $(\wedge^2\n/\mathfrak{a})_{\beta}$. Suppose it is zero on $(\ker\widetilde{\varphi}_0)_\beta$ i.e., on $\ker\widetilde{\varphi}_0$, then it factorizes through the quotient $\wedge^2(\n/\mathfrak{a})/(\ker\widetilde{\varphi}_0)$. Hence $\wedge^2(\n/\mathfrak{a})/(\ker\widetilde{\varphi}_0)$ and $\left[\n/\mathfrak{a},\n/\mathfrak{a}\right]$ are isomorphic and $\widetilde{\psi}$ may be written as $h\circ\widetilde{\varphi}_0$ for some $h\in\mathrm{Hom}_\K(\left[\n/\mathfrak{a},\n/\mathfrak{a}\right],\mathfrak{a})$ , thus the class of $\psi$ is zero. It follows that if the class of $\psi$ is nonzero then $\widetilde{\psi}$ can be identified with a linear form which is nonzero on $(\ker\widetilde{\varphi}_0)_\beta$.
\end{proof}

In the case of a simple path $(\alpha_1,...,\alpha_{n+1})$, the fiber of $\pi_{n+1}$ of $\varphi_n$ is controlled by the vector space $\mathrm{H}_2(\varphi_n)_{\alpha_{n+1}}$. Proposition \ref{nn51} and Proposition \ref{nn52} can be applied to the simple path of weights and give the following theorem which specifies the Zariski tangent spaces in subschemes defined by the fibers of $\pi_{n+1}$.
\begin{thm}\label{p6.1}Let $(\alpha_1,...,\alpha_{n+1})\subset\mathrm{T^*}$ be a simple path of weights with $n_0$ initialization and $\mathcal{A}_n\hookrightarrow\mathcal{A}_{n+1}$ a direct filiation for $n\geq n_0$. If $\varphi_n$ is a point of $\L_{n}^{\mathrm{T},\mathcal{A}_{n}}$ and $\varphi_{n+1}$ a point of $\pi_{n+1}^{-1}(\varphi_n)\cap\L_{n+1}^{\mathrm{T},\mathcal{A}_{n+1}}$ (assumed nonempty), then
\begin{enumerate}
	\item  the Zariski tangent at each point of $\pi_{n+1}^{-1}(\varphi_n)$ is given by
	$$\mathrm{Z}^2(\varphi_{n+1},\K e_{n+1})^{\mathrm{T}}\simeq(\mathrm{H}_2(\varphi_n)_{\alpha_{n+1}})^*$$
	\item the Zariski tangent of $\pi_{n+1}^{-1}(\varphi_{n})\cap\L_{n+1}^{\mathrm{T},\mathcal{A}_{n+1}}$ at point $\varphi_{n+1}$ is given by 
$$\mathrm{H}^2(\varphi_{n+1},\K e_{n+1})^{\mathrm{T}}\simeq\mathrm{Z}^2(\varphi_{n+1},\K e_{n+1})^{\mathrm{T}}/\K\psi',$$
where $\psi'$ is the extension of $\psi$ which is null on $\K^{n+1}\times\K e_{n+1}$, and  $\psi$ defines the central extension $\varphi_{n+1}$ of $\varphi_{n}$.\\
It injects in the Zariski tangent of $\L_{n+1}^{\mathrm{T},\mathcal{A}_{n+1}}$ at $\varphi_{n+1}$ via $\overline{i_2}$, 
	\item $\mathrm{H}^2(\varphi_{n+1},\K e_{n+1})^{\mathrm{T}}$ may be identified with the dual of the hyperplan $$\frac{(\ker\widetilde{\varphi}_0)_{\alpha_{n+1}}\cap(\ker\widetilde{\psi})}{\Omega(\n/\mathfrak{a})_{\alpha_{n+1}}}$$ 
of $\mathrm{H}_2(\varphi_n)_{\alpha_{n+1}}$.
\end{enumerate}
\end{thm}
\begin{proof} 1. If we solve the system of equations Eq. (\ref{e51}) in $\K$ we then obtain $\mathrm{Z}^2(\varphi_n,\K e_{n+1})^{\mathrm{T}}$; here $\mathrm{B}^2(\varphi_n,\K e_{n+1})^{\mathrm{T}}=0$.\\
2. From Eq. (\ref{e52}), we then obtain the kernel of the linear form which maps\\ $\xi\in\mathrm{Z}^2(\varphi_n,\K,e_{n+1})^{\mathrm{T}}$ onto $\xi^{n+1}_{pq}$. It is then a complement of $\K\psi'$. It follows that $\mathrm{B}^2(\varphi_{n+1},\K e_{n+1})^{\mathrm{T}}=\K\psi'$ and we deduce the result. The interpretation of $\overline{i_2}$ is obvious.\\
3. It follows from Proposition \ref{nn52} that the class of $\psi$ is nonzero, what is equivalent to $\psi\neq 0$.
\end{proof}
\subsection{Induction on local rings}
If $\mathcal{A}_n\hookrightarrow\mathcal{A}_{n+1}$ is a direct filiation, then the construction of the scheme $\L_{n+1}^{\mathrm{T},\mathcal{A}_{n+1}}$ is obtained by the scheme $\L_{n}^{\mathrm{T},\mathcal{A}_{n}}$ by vanishing the following new Jacobi's polynomials
\begin{equation}\label{e57}
	\J^{n+1}_{ijk}=\oint_{ijk}\sum_{l}Y^l_{ij}X^{n+1}_{lk},\quad 1<j<k\leq n,
\end{equation}
and the polynomials
\begin{equation}\label{e58}
	X^{n+1}_{pq}-1,
\end{equation}
where the variables $Y^l_{ij}$, $1<j\leq n,l\leq n$ are the ancient coordinates and $X^{n+1}_{lk}$, $(l,k\leq n)$ are the new coordinates. This remarkable distribution comes from the hypotheses on the weights, cf Proposition \ref{n41}. Denote by $X=(X_{ij}^{n+1})_{i<j\leq n}$ and $Y=(Y_{ij}^k)_{i<j\leq n,k\leq n}$. If $\O_n=\O_n(Y)$  is the local ring of $\L_{n}^{\mathrm{T},\mathcal{A}_{n}}$ at $\varphi_n$  then the local ring $\O_{n+1}=\O_{n+1}(X,Y)$ of $\L_{n+1}^{\mathrm{T},\mathcal{A}_{n+1}}$ at $\varphi_{n+1}$  is obtained by localizing the quotient ring
\begin{equation}\label{e6.10}
	\mathrm{A}=\O_n\left[X^{n+1}_{ij}\right]/\J,
\end{equation}
by the maximal ideal generated by the classes of representatives $f(X,Y)$ which vanish at point $\varphi_{n+1}$,
where $\J$ is the ideal generated by the polynomials (\ref{e57}) and (\ref{e58}). We obtain the following theorem with $\nu_n:=\dim\mathrm{H}_2(\varphi_n)_{\alpha_{n+1}}$.
\begin{thm}\label{t6.1}Let $\mathcal{A}_{n}\hookrightarrow\mathcal{A}_{n+1}$ be a direct filiation. Each point $\varphi_{n+1}$ of the fiber of $\varphi_n$ in $\L_{n+1}^{\mathrm{T},\mathcal{A}_{n+1}}(\K)$ verifies the following properties.
\begin{enumerate}
	\item If $\nu_n=0$ then there is not the fiber
	\item If $\nu_n=1$ then $\varphi_{n+1}$ is unique in the fiber and $\O_{n+1}$ is a quotient of $\O_n$. Moreover if $\varphi_n$ is rigid then $\varphi_{n+1}$ is rigid.
	\item If $\nu_n>1$ then $\O_{n+1}$ is the localization ring of the quotient $\O_n[T_1,...,T_{\nu_n-1}]$ by the ideal generated by the expressions
	$$c_0^{\lambda}(Y)+\sum_{0<k<\nu_n}c_k^{\lambda}(Y)T_k,\lambda\in\Lambda-\Lambda_1,\,\,c_l^{\lambda}(Y)\in\m(\O_n),\,0\leq l\leq \nu_n,$$
	at the point $(T_1,...,T_{\nu_n-1})=(0,...,0)$, where $\Lambda$ is the set of multi-indices $(ijk)$ satisfying $i<j<k$, which indexes the Jacobi's polynomials associated with $\alpha_i+\alpha_j+\alpha_k=\alpha_{n+1}$ and $\Lambda_1$ is a subset of $\Lambda$ defined in the proof.
	The dimension of Krull $d$ of the scheme $\L_{n+1}^{\mathrm{T},\mathcal{A}_{n+1}}$ at point $\varphi_{n+1}$ is minored by $\nu_n-1$. In particular $\varphi_{n+1}$ is not rigid.
\end{enumerate}
\end{thm}
\begin{proof} The case $\nu_n=0$ is trivial. If $\nu_n>0$ the Jacobi's polynomials indexed by the set $\Lambda$ of triples $i<j<k$ such that $\alpha_i+\alpha_j+\alpha_k=\alpha_{n+1}$ and the condition $X^{n+1}_{pq}=1$ may be written as
$$a_0^{\lambda}(Y)+\sum_{\mu\neq (p,q)}a_{\mu}^{\lambda}(Y)X_\mu,\quad\lambda\in\Lambda,$$
where $\mu$ runs through the set $M$ of couples $i<j$ such that $\alpha_i+\alpha_j=\alpha_{n+1}$, $(i,j)\neq (p,q)$ and $a_0^{\lambda},a_\mu^{\lambda}\in\O_n$. Let $(A)$ be the system obtained by vanishing those polynomials. It we take $Y=\varphi_n$, the system of equations $(A)$ becomes a system of equations $(B)$ which the solutions in $X$ are formed by $\Psi\in{Z}^2(\varphi_n,\K e_{n+1})^{\mathrm{T}}$ such that $\Psi^{n+1}_{pq}=1$. It gives an affine space of dimension $\nu_n-1$, see Proposition \ref{p6.1}. Let $\Lambda_1$ be a subset of $\Lambda$ such that the subsystem of equations of $(B)$ indexed by $\Lambda_1$ is equivalent to $(B)$ and $\Lambda_1$ is minimal for this property. There is a subset $M_1$ of $M$ of cardinal $|\Lambda_1|$ such that the submatrix $(a^\lambda_\mu(\varphi_n))_{(\lambda,\mu)\in\Lambda_1\times M_1}$ is invertible. This matrix remains invertible if we take its values in $\O_n$ since its projection by $\mathrm{pr}:\O_n\rightarrow\O_n/\m(\O_n)$ is. It follows that the subsystem of equations of $(A)$ indexed by $\Lambda_1$ permits to express the variables $(X^\mu)_{\mu\in M_1}$ , in function of variables $(X^\rho)_{\rho\in M-M_1}$, i.e.
\begin{equation}
	X^\mu=A^\mu_0(Y)+\sum_\rho A^\mu_\rho(Y)X^\rho\quad \mu\in M_1,
\end{equation}
where $A^\mu_0(Y),A^\mu_\rho(Y)\in\O_n$.
We substitute these expressions in the remaining expressions of $(A)$ for multi-indices $\Lambda-\Lambda_1$, which gives the new system of equations by vanishing the following expressions:
\begin{equation}\label{e6.12}
	b_0^\lambda(Y)+\sum_{\rho\in M-M_1} b^\lambda_\rho(Y)X^\rho\quad \lambda\in\Lambda-\Lambda_1.
\end{equation}
If we take $Y=\varphi_n$, then the system of equations (\ref{e6.12}) corresponds to the system of equations $(B)$ indexed by $\Lambda-\Lambda_1$. Theses equations give anything more on $\K$ and consequently are identically null, i.e $b^\lambda_0(Y)=b^\lambda_\rho(Y)=0$ for $\lambda\in\Lambda-\Lambda_1$ and $\rho\in M-M_1$. In other words, the elements $b^\lambda_0(Y)$ and $b^\lambda_\rho(Y)$ belong to the maximal ideal of $\O_n$. The ring $\mathrm{A}$, see (\ref{e6.10}) is equal to the quotient of $\O_n[X^\rho,\rho\in M-M_1]$ by the ideal generated by the terms (\ref{e6.12}). If we want to localize this ring at point $\varphi_{n+1}$, we will take new adapted variables $X^\rho-(\varphi_{n+1})^\rho$ which we write $T_k$ with indexation $k$ on $\left\{1,...,\nu_n-1\right\}$. Then the Jacobi's polynomials are written as $c_0^{\lambda}(Y)+\sum_{0<k<\nu_n}c_k^{\lambda}(Y)T_k,\quad\lambda\in\Lambda-\Lambda_1,$
where $c_l^{\lambda}(Y)\in\m(\O_n)$, $(0\leq l< \nu_n)$ and $\mathrm{A}$ is written as in $3.$ We localize at point $\varphi_{n+1}$ by considering the maximal ideal obtained as quotient of the maximal ideal of $\O_n[T]$ given by $\m(\O_n)[T]+\sum_{1\leq k<\nu_n}T_k.\O_n[T]$.\\
If $\nu_n=1$, we see that the ring $\A$ is reduced to the quotient of $\O_n$ by the ideal generated by the $b_0^\lambda$, $\lambda\in\Lambda-\Lambda_1$. We obtain the result.\\
If $\nu_n>1$, there are $\nu_n-1$ new parameters which are algebraically independent over $\K$, the dimension of Krull then satisfies $d\geq\nu_n-1$.
\end{proof}
\begin{rem}We observe the following properties:\
\begin{enumerate}
	\item a nilpotent parameter is a free ancient parameter appeared in lower dimension. In next section we will give examples which illustrate this phenomenon. The order of nilpotence of a such parameter can decrease strictly for an extension in higher dimension, cf \cite{CC}.
	\item if $\nu_n>1$ we find the same interpretation of genesis of the parameters as in the formalism by generators and relations, in \cite{Fa}. We have here the Grassmanian $\mathrm{Gr}_r(\mathrm{H}_2(\varphi_n)_\beta)$ of dimension $r(\nu_n-r)$ for $r=1$ and it is known that $\mathrm{H}_2(\varphi_n)_\beta$ is also interpreted in the free Lie algebra formalism, \cite{C8}.
\end{enumerate}
\end{rem}
\section{Application}
A Lie algebra may be written as $\left[e_i,e_j\right]=c_{ij}e_{i+j}$, $(i+j\leq n)$ on the basis $(e_i)_{i=1,...,n}$ where $c_{ij}\in\K$ if and only if there is a derivation $\delta$ defined by $\delta.e_i=ie_i$ for $1\leq i\leq n$. This means that there exists a maximal torus $\mathrm{T}$ containing the torus $\mathrm{T}_1=\K\delta$. The weights of $\mathrm{T}$ on the basis $(e_i)$  satisfy $\alpha_i|_{\mathrm{T_1}}=i\alpha_1|_{\mathrm{T_1}}$, $1\leq i\leq n$. The coordinates $X_{ij}$ of the scheme $\L_n^{\mathrm{T}}$ are given by the pairs $(i<j)$ such that $\alpha_i+\alpha_j=\alpha_{i+j}$ belongs to $\mathrm{P}_n$. The Jacobi's polynomials may be written as:
\begin{equation}
	\mathrm{J}_{ijk}=X_{i,j}X_{i+j,k}+X_{j,k}X_{j+k,i}+X_{k,i}X_{k+i,j}
\end{equation}
for $i<j<k$ and $\alpha_i+\alpha_j+\alpha_k=\alpha_{i+j+k}\in\mathrm{P}_n$. \\
One can look for examples of maximal torus for Lie algebras which the structure constants $c_{ij}$ $(i<j)$ are equal $0$ or $1$. One finds such a Lie algebra $\varphi$ by the set $\mathcal{J}_\varphi$ consisting of pairs $(i<j)$ such that $i+j\leq n$ and $c_{ij}\neq 0$. If one takes for $3\leq p\leq n-2$ and considers the following Lie algebra $\mathfrak{a}_{p,n}$ defined by 
$$\left[e_1,e_i\right]=e_{i+1},\quad i=2,3,...,\widehat{p-1},...,n-1$$
$$\left[e_2,e_i\right]=e_{i+2},\quad i=3,4,...,\widehat{p-2},\widehat{p-1},...,n-2$$
where the symbol $\widehat{x}$ indicates that $x$ must be omitted. The maximal torus $\mathrm{T}$ is of dimension $2$ and is defined by the following weights
\begin{equation}\label{e6.1}
	\alpha_i=i\alpha_1,(1\leq i\leq p-1),\quad\alpha_i=(i-p)\alpha_1+\alpha_{p},\quad(p\leq i\leq n).
\end{equation}
For $p\geq 3$ one can see that a such sequence of weights defines a simple path of weights, then the methods defined in the section $6$ can be applied.\\
For $p=3$ one obtains a sequence of rigid Lie algebras $\mathfrak{a}_{3,n}$. On will develop in the next section the study of the direct filiation for the path of weights (\ref{e6.1}) for $p=4$. One will also study the path of weights $i\alpha_1$, $i=1,...,n.$
\subsection{Study of $\L_n^{\mathrm{T},\mathcal{A}_n}:\mathcal{A}_n=\left\{(2,4),(1,k):k\neq 3, 1<k<n\right\}$}
One considers the path of weights below: $$\alpha_1,2\alpha_1,3\alpha_1,\alpha_4,\alpha_4+\alpha_1,...,\alpha_4+k\alpha_1,...$$
-It is initialized for $n=4+k=6$ which corresponds to the algebra $\mathfrak{a}_{4,6}$ denoted by $\n_{6,7}$ in \cite{V}. One has only one admissible set
$$\mathcal{A}_6=\left\{(1,2),(1,4),(1,5),(2,4)\right\},$$
hence
$$\L_6^{\mathrm{T},\mathcal{A}_6}=\left\{\mathfrak{a}_{4,6}\right\},$$
and $\sum_6$ is equal to the orbit of $\mathfrak{a}_{4,6}$.\\
The local ring $\O_6$ at the point $\mathfrak{a}_{46}$ to $\L_n^{\mathrm{T},\mathcal{A}_6}$ is isomorphic to $\K$.\\
-For $n=7$ one has 3 new coordinates $X_{16},X_{25},X_{34}$ and one Jacobi's polynomial $\mathrm{J}_{124}$, i.e.
\begin{equation}\label{e7.3}
	X_{25}=X_{16}-X_{34}.
\end{equation}
The central extension is not trivial if $X_{16}\neq 0$, $X_{25}\neq 0$, or  $X_{34}\neq 0$. Hence there are $3$ admissible sets possible
$$\mathcal{A}_7=\mathcal{A}_6\cup\left\{(16)\right\},\mathcal{A'}_7=\mathcal{A}_6\cup\left\{(25)\right\}, \mathrm{or}     	    \quad\mathcal{A''}_7=\mathcal{A}_6\cup\left\{(34)\right\}.$$
One considers only the first admissible set. It permits to fix $X_{1,6}=1$ and Eq.(\ref{e7.3}) implies that
 \begin{equation}
	X_{1,6}=1,\quad X_{2,5}=1-t,\quad t:=X_{34}.
\end{equation}
Hence
 $$\L_7^{\mathrm{T},\mathcal{A}_7}=\left\{ \mathfrak{a}_{47}(t)\right\},$$ 
where $\mathfrak{a}_{47}(t)$ is the family defined by adding the following new brackets:
$$\left[e_1,e_6\right]=e_7,\quad\left[e_2,e_5\right]=(1-t)e_7,\quad\left[e_3,e_4\right]=te_7,$$
where $t\in\K$.\\
-For $n=8$ one has 3 new coordinates $X_{17},X_{26},X_{35}$ and two Jacobi's polynomials $\mathrm{J}_{125}$ and $\mathrm{J}_{134}$ , i.e.
\begin{equation}\label{e6.2}
	X_{26}=(1-t)X_{17}-X_{35},\quad X_{35}=tX_{17}.
\end{equation}
Ones takes $\mathcal{A}_8:=\mathcal{A}_7\cup\left\{(17)\right\}$ for $X_{17}\neq 0$ which permits to fix $X_{17}=1$. Eq. (\ref{e6.2}) is equivalent to
\begin{equation}\label{e6.3}
X_{1,7}=1,\quad X_{26}=1-2t,\quad X_{35}=t.
\end{equation}
Hence
$$\L_8^{\mathrm{T},\mathcal{A}_8}=\left\{\mathfrak{a}_{48}(t)\right\},$$
where $\mathfrak{a}_{48}(t)$ is the family defined by adding the following new brackets:
$$\left[e_1,e_7\right]=e_8,\quad\left[e_2,e_6\right]=(1-2t)e_8,\quad\left[e_3,e_5\right]=te_8,
$$
where $t\in\K$. \\
-For $n=9$ one has $3$ new coordinates $X_{18},X_{27},X_{36}$ and three Jacobi's polynomials $\mathrm{J}_{126}$ ,$\mathrm{J}_{135}$ and $\mathrm{J}_{234}$, i.e.
\begin{equation}\label{e6.4}
	X_{36}=(1-2t)X_{18}-X_{27},\quad X_{36}=tX_{18},\quad X_{36}=tX_{27}.
\end{equation}
On sets $\mathcal{A}_9=\mathcal{A}_8\cup\left\{(18)\right\}$ for $X_{18}\neq 0$ which permits to fix $X_{18}=1$. Eq. (\ref{e6.4}) gives
\begin{equation}\label{e6.5}
	X_{36}=(1-2t)-X_{27},\quad X_{36}=t,\quad X_{36}=tX_{27}.
\end{equation}
Hence
\begin{equation}\label{e6.6}
	X_{18}=1,\quad X_{27}=1-3t,\quad X_{36}=t,\quad -3t^2=0.
\end{equation}
One has only a geometric point $\mathfrak{a}_{49}(0)$ given by the value $t=0$ in $\K$ since $t^2=0$. From Theorem \ref{t6.1}, the local ring at this point of $\L_9^{\mathrm{T},\mathcal{A}_9}$ is given by $\O_9\simeq\O_8/(t^2)\cong\K[u]/(u^2)=\K+\K t$ with $t^2=0$. The continuous family is reduced by extension to one point $\mathfrak{a}_{49}$; the parameter $t$ which was free for $n=7,8$ became nilpotent for $n=9$.  
\begin{thm}\label{t7.1}The scheme $\L^{\mathrm{T},\mathcal{A}_n}_n$ where $\mathcal{A}_n=\left\{(2,4),(1,k):k\neq 3, 1<k<n\right\}$ is given by
\begin{enumerate}
	\item for $n=6$, the unique law $\mathfrak{a}_{4,6}$ and $\O_6\simeq\K$;
	\item for $n=7,8$, the continuous family $\mathfrak{a}_{4,n}(t)$ defined by the brackets\\
	$\left[e_1,e_i\right]=e_{i+1}$, $1<i<n$ and $i\neq 3$,\\ 
	$\left[e_2,e_i\right]=(1-(i-4)t)e_{i+2}$, $3<i<n-1$,\\
	$\left[e_3,e_i\right]=te_{i+3}$, $3<i<n-2$, \\
where $t\in\K$ and $\O_n\simeq\left\{\frac{p}{q}:p,q\in\K\left[u\right],q(0)\neq 0\right\}\subset\K(u)$ at each point;
\item for $n\geq 9$, the unique law $\mathfrak{a}_{4,n}$ defined for $t=0$. The canonical deformation in the schemes writes as $\mathfrak{a}_{4,n}(t)$ (above formula) with one nilpotent element $t$ with order two. $\O_n=\K+\K\cdot t\simeq\K[u]/(u^2)$.
\end{enumerate}
\end{thm}
\begin{proof} 3. One reasons by induction on $n\geq 9$. It is verified for $n=9$ (see above). One assumes that the property is satisfied for $n-1\geq 9$. For $n$ one has  three new coordinates $X_{1,n-1}$,$X_{2,n-2}$, and $X_{3,n-3}$ and Jacobi's polynomials $\mathrm{J}_{12,n-3}$, $\mathrm{J}_{13,n-4}$ and $\mathrm{J}_{23,n-5}$, i.e.
$$
X_{3,n-3}-(1-(n-7)t)X_{1,n-1}+X_{2,n-2}=0,\quad X_{3,n-3}=tX_{1,n-1}$$
and
$$(1-(n-9)t)X_{3,n-3}=tX_{2,n-2}.$$
It is necessary that $X_{1,n-1}\neq 0$, and so ones sets $\mathcal{A}_n=\mathcal{A}_{n-1}\cup \left\{(1,n-1)\right\}$ which permits to fix $X_{1,n-1}=1$. The above equations imply that
\begin{equation}
	X_{3,n-3}=t,\quad X_{2,n-2}=1-(n-6)t,\quad t^2=0.
\end{equation}
Using the induction hypothesis one gets the result.
\end{proof}
\begin{cor} The Lie algebra $\mathfrak{a}_{4,n}$ (resp. $\mathrm{T}\ltimes\mathfrak{a}_{4,n}$) is rigid in $\L^{\mathrm{T}}_n$ (resp. $\L_{n+2}$) for $n\geq 9$.
\end{cor}
\begin{proof} It follows that from Theorem \ref{t7.1} and the reduction theorem.
\end{proof}
\begin{rem}There is here a filiation which is not direct since the open set $\sum_7$ of $\L^{\mathrm{T}_2}_7$ is formed by the family $\mathfrak{a}_ {47}(t)$ and the laws $2.1 (ii), 2.1 (iii), 2.1 (iv)$ and $2.1 (v)$ in \cite{C5}. The Lie algebras quotient of the laws  $2.1 (ii), 2.1 (iii), 2.1 (iv)$ and $2.1 (v)$  by their centers are isomorphic to $\n_{66}$, $\K\times\n_{53}$ and $\n_{61}$ (we use the notation in \cite{V}); these Lie algebras admit a maximal torus of dimension $3$ and so do not belong to $\sum_6$.
\end{rem}
\subsection{Construction of the nilpotent element in the ring which defines the scheme $\L^{\mathrm{T}}_n, n\geq 9$}
The existence of a nilpotent element in the local ring $\O_n$ of the scheme $\L^{\mathrm{T}}_n$ at the point $\mathfrak{a}_{4,n}$ for $n\geq 9$, comes from a nilpotent element in the ring $\I_n=\K\left[X_{ij}\right]/\J_n$ which defines the scheme $\L^{\mathrm{T}}_n, n\geq 9$. \\
The method consists to eliminate the variables $X_{ij}$, $(i<j,i+j\leq 9)$ by using Jacobi's polynomials $\mathrm{J}_{ijk}$, $(i<j<k, i+j+k\leq 9)$ modulo the ideal $\J_n$ and preserving $X_{18},X_{17}$ and the coordinate $X_{34}$ which corresponds to the parameter. One has the Jacobi's polynomials:
\begin{equation}
	\mathrm{J}_{126}=X_{12}X_{36}-X_{26}X_{18}+X_{16}X_{27}.
\end{equation}
 \begin{equation}
	\mathrm{J}_{135}=-X_{35}X_{18}+X_{15}X_{36}.
\end{equation}
\begin{equation}
	\mathrm{J}_{234}=-X_{34}X_{27}+X_{24}X_{36}.
\end{equation}
One eliminates $X_{36}$ by the way below :
\begin{equation}\label{e6.14}
	X_{15}\J_{126}-X_{12}\J_{135}=X_{12}X_{35}X_{18}-X_{15}X_{26}X_{18}+X_{15}X_{16}X_{27}
\end{equation}
\begin{equation}\label{e6.15}
	X_{15}\J_{234}-X_{24}\J_{135}=-X_{15}X_{34}X_{27}+X_{24}X_{35}X_{18}
\end{equation}
One eliminates $X_{27}$ by the following way:
$$X_{34}\times(\ref{e6.14})+X_{16}\times(\ref{e6.15})=$$
\begin{equation}
	X_{18}(X_{34}X_{12}X_{35}-X_{34}X_{15}X_{26}+X_{16}X_{24}X_{35})=:f_1\in\J_n
\end{equation}
In dimension $8$, one also has Jacobi's polynomials:
\begin{equation}
	\J_{134}=-X_{34}X_{17}+X_{14}X_{35}
\end{equation}
\begin{equation}
	\J_{125}=X_{12}X_{35}-X_{25}X_{17}+X_{15}X_{26}
\end{equation}
\begin{equation}\label{e6.19}
	\J_{124}=X_{12}X_{34}-X_{24}X_{16}+X_{14}X_{25}
\end{equation}
One eliminates $X_{35}$ by multiplying $f_1$ by $X_{14}$ and by using Jacobi's polynomial $\J_{134}$, one obtains:
\begin{equation}\label{e6.20}
	f_2:=X_{18}X_{34}(X_{12}X_{34}X_{17}-X_{15}X_{26}X_{14}+X_{16}X_{24}X_{17})\in\J_n
\end{equation}
One eliminates $X_{35}$ by the way below:
\begin{equation}\label{e6.21}
	f_3:=X_{14}\J_{125}-X_{12}\J_{134}=-X_{14}X_{25}X_{17}+X_{14}X_{15}X_{26}+X_{12}X_{34}X_{17}
\end{equation}
One eliminates $X_{26}$ by using (\ref{e6.20}) and (\ref{e6.21}), we have
\begin{equation}\label{e6.22}
	f_4:=X_{18}X_{34}X_{17}(2X_{12}X_{34}-X_{14}X_{25}+X_{16}X_{24})\in\J_n
\end{equation}
By (\ref{e6.19}) and (\ref{e6.22}) we have
\begin{equation}\label{e6.23}
	X_{12}X_{17}X_{18}(X_{34})^2\in\J_n.
\end{equation}
If one denotes by $\overline{X}_{ij}$ the quotient modulo $\J_n$, one then obtains:
\begin{thm}\label{n6.2} We have $\overline{X}_{12}\overline{X}_{17}\overline{X}_{18}(\overline{X}_{34})^2 =0$ in the ring $\K[X_{ij}]/\J_n$ for all $n\geq 9$ and there is not null product with a less nombre of factors. In particular, the product $\overline{X}_{12}\overline{X}_{17}\overline{X}_{18}\overline{X}_{34}\neq 0 $ is nilpotent with order $2.$
\end{thm}
\begin{proof} It remains to verify that each product $Z$, if we remove exactly one of the factor $\overline{X}_{12},\overline{X}_{17},\overline{X}_{18}$, and $\overline{X}_{34} $ in $\overline{X}_{12}\overline{X}_{17}\overline{X}_{18}(\overline{X}_{34})^2 $, defines with the equation $Z\neq 0$ a nonempty Zariski open set.\\
$i)$ If we remove $\overline{X}_{12}$, then this open contains the following Lie algebras $\got{b}_n$ defined by the following brackets:\\
$\left[e_1,e_k\right]=e_{k+1}$, $(3<k<n)$, $\left[e_2,e_k\right]=e_{k+2}$, $(3<k<n-1)$,\\
$\left[e_3,e_k\right]=e_{k+3}$, $(3<k<n-2)$.\\
$ii)$ If we retire $\overline{X}_{17}$, then the open contains the following Lie algebra for $t\neq 0$:\\
$\left[e_1,e_k\right]=e_{k+1}$, $(1<i<n)$ and $i\neq 3,7$, $\left[e_2,e_4\right]=e_6,$ $\left[e_2,e_5\right]=(1-t)e_7,$\\ $\left[e_3,e_4\right]=te_7$.\\ 
$iii)$ If we retire $\overline{X}_{18}$, then the open contains the  Lie algebras $\mathfrak{a}_{4,8}(t)\times\K^{n-8}$ where $\K^{n-8}$ is the abelian ideal $\K e_9\oplus\cdots\K e_n$.\\
$iv)$ $Z=\overline{X}_{1,2}\overline{X}_{1,7}\overline{X}_{1,8}\overline{X}_{3,4}$, we have $Z^2=0$ and the open $Z\neq 0$ contains $\mathfrak{a}_{4,n}$. For $n>12$ by using Theorem \ref{t7.1} and the fact that $\overline{X}_{3,4}$ gives a nonzero nilpotent element in the local ring at $\mathfrak{a}_{4,n}$.
\end{proof}
\begin{cor} The variety $\L^{\mathrm{T}}_n(\K)$ contains at least $4$ irreducible components for $n\geq 9$.
\end{cor}
\begin{proof} For $n\geq 9$ the Lie algebras $\mathfrak{a}_{4,n}$ and $\got{b}_n$ are rigid, contained in $\sum_n$ and do not belong to the open set defined by $\overline{X}_{12}\overline{X}_{34}^2\neq 0$. It follows that from (\ref{e6.23}) the geometric points of $\L^{\mathrm{T}}_n(\K)$ which belong to this open set verify $\overline{X}_{17}=0$ or $\overline{X}_{18}=0$. This give two disjoint open set $\overline{X}_{17}\neq 0$ $(\overline{X}_{18}=0)$ and $\overline{X}_{18}\neq 0$ $(\overline{X}_{17}=0)$ which contain the laws $ii)$ and $iii)$ in the proof of Theorem \ref{n6.2} for $t\neq 0$. 
\end{proof}
\subsection{Study of $\L^{\mathrm{T},\mathcal{A}_n}_n:\mathcal{A}_n=\left\{(23),(1k),1<k<n\right\}$}
The study of simple paths of weights $\alpha_k=k\alpha_1$ with $k>0$ is appeared in (\cite{Br},\cite{C2}) and under the scheme form announced in \cite{C3}. There is another method by generators and relations \cite{F}.\\
The initialization is made for $n=5$ since the Lie algebra\\
$\n_{56}:\left[e_1,e_k\right]=e_{k+1}, (1<k<5), \left[e_2,e_3\right]=e_5$\\
is the first and unique Lie algebra with a maximal torus given by the weights $k\alpha_1$ for $n=5$. One has $\L^{\mathrm{T},\mathcal{A}_5}_5(\K)=\left\{\n_{56}\right\}$ with $\mathcal{A}_5=\left\{(23),(1k), 1<k<5\right\}$.\\
-For $n=6$ on obtains two new coordinates $X_{15},X_{24}$ and one relation $\J_{123}$, i.e.
\begin{equation}
	X_{24}=X_{15}.
\end{equation}
One sets $\mathcal{A}_6=\mathcal{A}_5\cup\left\{(15)\right\}$ for $X_{24}=X_{15}\neq 0$ which permits to fix $X_{24}=X_{15}=1.$
The scheme $\L^{\mathrm{T},\mathcal{A}_6}_6$ is consisted of the only Lie algebra denoted by $\n_{619}$ in \cite{V}.\\
-For $n=7$ we have three new coordinates $X_{16},X_{25},X_{34}$ and one Jacobi polynomial $\J_{124}$ which implies
\begin{equation}\label{e6.25}
	X_{34}-X_{16}+X_{25}=0
\end{equation}
On sets $\mathcal{A}_7=\mathcal{A}_6\cup\left\{(16)\right\}$ for $X_{16}\neq 0$ which permits to fix $X_{16}=1.$ We set $X_{34}=t$, it follows from Eq. (\ref{e6.25}) that $X_{25}=1-t$. We obtain a family of filiform Lie algebras $f_7(t)$ which describes the points of $\L^{\mathrm{T},\mathcal{A}_7}_7(\K)$ with $t\in\K$. The scheme $\L^{\mathrm{T},\mathcal{A}_7}_7$ may be identified with $\mathrm{Spec}(\K[u])$ where $u$ is an indeterminate.\\
-For $n=8$ we have three new coordinates $X_{17},X_{26},X_{35}$ and two Jacobi's polynomials $\J_{125}$ and $\J_{134}$,i.e.
\begin{equation}
	X_{35}-(1-t)X_{17}+X_{26}=0,
\end{equation}
 \begin{equation}
	-tX_{17}+X_{35}=0.
\end{equation}
One sets $\mathcal{A}_8=\mathcal{A}_7\cup\left\{(17)\right\}$ for $X_{17}\neq 0$ which permits to fix $X_{17}=1.$ The variety $\L^{\mathrm{T},\mathcal{A}_8}_8(\K)$ is formed of a family $f_8(t)$:
\begin{equation}
	X_{17}=1, X_{26}=1-2t, X_{35}=t.
\end{equation}
The scheme $\L^{\mathrm{T},\mathcal{A}_8}_8$ is identified with $\mathrm{Spec}(\K[u])$.\\
-For $n=9$ we have four new coordinates $X_{18},X_{27},X_{36},X_{45}$ and three Jacobi's polynomials $\J_{126},\J_{135},\J_{234}$, i.e.
\begin{equation}
	X_{36}-(1-2t)X_{18}+X_{27}=0
\end{equation}
\begin{equation}
	X_{45}-tX_{18}+X_{36}=0
\end{equation}
\begin{equation}
	-X_{45}-tX_{27}+X_{36}=0.
\end{equation}
One sets $\mathcal{A}_9=\mathcal{A}_8\cup\left\{(18)\right\}$ for $X_{18}\neq 0$ which permits to fix $X_{18}=1.$ Hence for $t\neq -2$ we have
\begin{equation}\label{e6.32}
	X_{18}=1, \ X_{27}=\frac{2-5t}{2+t}, \ X_{36}=\frac{2t-2t^2}{2+t}, \ X_{45}=\frac{3t^2}{2+t}.
\end{equation}
(\ref{e6.32}) is valid in the local ring at each point of $\L^{\mathrm{T},\mathcal{A}_9}_9(\K)$ since the projection of $2+t$ on $\K$ is different from zero and $2+t$ is invertible if $t\neq -2$.
\begin{rem} Concerning the condition $X_{18}\neq 0$ which corresponds to the choice of $\mathcal{A}_9$, we see that (\ref{e6.32}) does not admit a solution if \ $t=-2$ and the Lie algebra $f_8(-2)$ does not have a central extension in $\L^{\mathrm{T},\mathcal{A}_9}_9$. However if we choose another admissible set $\mathcal{A'}_9=\mathcal{A}_8\cup\left\{(27)\right\}$ then the central extension exists and the Lie algebra obtained $f_9(-2)$ satisfies
$X_{18}=0,X_{27}=1, X_{36}=-1, X_{45}=1.$
\end{rem}For $n=10$ we have four new coordinates $X_{19},X_{28},X_{37},X_{46}$ and four Jacobi's polynomials $\J_{127},\J_{136},\J_{145},\J_{235}$, i.e.
\begin{equation}
	X_{37}-\frac{2-5t}{2+t}X_{19}+X_{28}=0
\end{equation}
\begin{equation}
	X_{46}-\frac{2t-2t^2}{2+t}X_{19}+X_{37}=0
\end{equation}
\begin{equation}
	\frac{-3t^2}{2+t}X_{19}+X_{46}=0
\end{equation}
\begin{equation}
	-tX_{28}+(1-t)X_{37}=0.
\end{equation}
Ones sets $\mathcal{A}_{10}=\mathcal{A}_9\cup\left\{(19)\right\}$ for $X_{19}\neq 0$ which permits to fix $X_{19}=1$ and gives
\begin{equation}\label{e6.37}
	X_{19}=1, \ X_{28}=\frac{2-7t+5t^2}{2+t}, \ X_{37}=\frac{2t-5t^2}{2+t}, \ X_{46}=\frac{3t^2}{2+t}.
\end{equation}
in the associated local ring at each value $\mathrm{pr}(t)\in\K-\left\{-2\right\}$.\\
-For $n=11$ we have five new coordinates $X_{110},X_{29},X_{38},X_{47},X_{56}$ and five Jacobi's polynomials $J_{128},J_{137},J_{146},J_{236},J_{245}$ which give
\begin{equation}
	X_{38}-(-17+5t+\frac{36}{2+t})X_{110}+X_{29}=0
\end{equation}
\begin{equation}
	X_{47}-(12-5t-\frac{24}{2+t})X_{110}+X_{38}=0
\end{equation}
\begin{equation}
	X_{56}-(-6+3t+\frac{12}{2+t})X_{110}+X_{47}=0
\end{equation}
\begin{equation}
	X_{56}-(6-2t-\frac{12}{2+t})X_{29}+(1-2t)X_{38}=0
\end{equation}
\begin{equation}
-X_{56}-(-6+3t+\frac{12}{2+t})X_{29}+(1-t)X_{47}=0
\end{equation}
One sets $\mathcal{A}_{11}=\mathcal{A}_{10}\cup\left\{(110)\right\}$ for $X_{110}\neq 0$ permits to fix $X_{110}=1$ and gives
\begin{align}\label{e6.43}
	X_{110}=1,X_{29}=\frac{2-10t+16t^2-5t^3}{2(1-t^2)},X_{38}=\frac{4t-16t^2+8t^3-5t^4}{2(2+t)(1-t^2)},\\X_{47}=\frac{6t^2-12t^3+15t^4}{2(2+t)(1-t^2)},X_{56}=\frac{12t^3-21t^4}{2(2+t)(1-t^2)}.
\end{align}
in the associated local ring at each value $\mathrm{pr}(t)\in\K-\left\{-2,\pm 1\right\}$.\\
-For $n=12$ we have five new coordinates $X_{111},X_{210},X_{39},X_{48},X_{57}$ and seven Jacobi's polynomials $\J_{129},\J_{138},\J_{147},\J_{156},\J_{237},\J_{246},\J_{345}$, i.e.
\begin{equation}
	X_{39}-\frac{2-10t+16t^2-5t^3}{2(1-t^2)}X_{111}+X_{210}=0
\end{equation}
\begin{equation}
	X_{48}-\frac{4t-16t^2+8t^3-5t^4}{2(2+t)(1-t^2)}X_{111}+X_{39}=0
\end{equation}
\begin{equation}
X_{57}-\frac{6t^2-12t^3+15t^4}{2(2+t)(1-t^2)}X_{111}+X_{48}=0
\end{equation}
\begin{equation}
	X_{57}-\frac{12t^3-21t^4}{2(2+t)(1-t^2)}X_{111}=0
\end{equation}
\begin{equation}
X_{57}-(\frac{2t-5t^2}{2+t})X_{210}+\frac{2-5t}{2+t}X_{39}=0
\end{equation}
\begin{equation}
-\frac{3t^2}{2+t}X_{210}+(1-2t)X_{48}=0
\end{equation}
\begin{equation}
-tX_{57}-\frac{3t^2}{2+t}X_{39}+tX_{48}=0
\end{equation}
Set $\mathcal{A}_{12}=\mathcal{A}_{11}\cup\left\{(111)\right\}$ for $X_{110}\neq 0$ permits to fix $X_{111}=1$ and gives
\begin{align}\label{e6.51}
	X_{210}=\frac{4-22t+44t^2-26t^3+36t^4}{2(2+t)(1-t^2)},X_{39}=\frac{4t-22t^2+32t^3-41t^4}{2(2+t)(1-t^2)},\\X_{48}=\frac{6t^2-24t^3+36t^4}{2(2+t)(1-t^2)},X_{57}=\frac{12t^3-21t^4}{2(2+t)(1-t^2)}
\end{align}
\begin{equation}
	\J_{237}=-\J_{246}=\J_{345}=\frac{9t^5(10t-1)}{(2+t)^2(1-t)^2}=0.
\end{equation}
This shows that the scheme $\L^{\mathrm{T},\mathcal{A}_{12}}_{12}$ is consisted of the only two rigid Lie algebras denoted by $\mathfrak{f}_{12}$ for $t=0$ and the Witt Lie algebra $\mathfrak{w}_{12}$ for $t=1/10$. The local ring at point $\mathfrak{w}_{12}$ is trivial $(\simeq\K)$ and those of $\mathfrak{f}_{12}$ is isomorph to $\K[u]/u^5$ with a nilpotent element $t$ of order $5$ ($t^5=0$ and $t^4\neq 0$). We deduce that $\mathrm{H}^2(\mathfrak{w}_{12},\mathfrak{w}_{12})^{\mathrm{T}}=0$ and $\mathrm{H}^2(\mathfrak{f}_{12},\mathfrak{f}_{12})^{\mathrm{T}}\simeq\K$ (we can also verify by direct compute in \cite{C2}). We obtain the
\begin{thm}The scheme $\L^{\mathrm{T},\mathcal{A}_{n}}_{n}$  is given by
\begin{enumerate}
	\item for $n=5,6$, a rigid Lie algebra. The scheme is reduced and the local ring is trivial $(\simeq\K)$.
	\item For $7\leq n\leq 11$, we have a continuous family of Lie algebras $\mathfrak{f}_n(t)$ where $t$ runs through $\K$ except for a finite numbers of points defined by the brackets above.
	\item For $n\geq 12$, two rigid Lie algebras $\mathfrak{f}_n$ if $t=0$ and $w_n$ if $t=1/10$. The scheme is reduced at the point $\mathfrak{w}_n$ and its local ring is trivial $(\simeq\K)$. The scheme is not reduced at the point $\mathfrak{f}_n$ and possesses a nilpotent element $t$ with order $5$. The local ring may be written as
	$$\oplus_{k=0}^4\K t^k\simeq\K[u]/u^5.$$
The canonical deformation of $\mathfrak{f}_n$ in this scheme is defined by the following brackets:\\
$X_{1i}=1, (1<i<n)$\\
$X_{23}=X_{24}=1$,$ X_{25}=1-t$, $X_{26}=1-2t$,\\
$X_{2,m-2}=1+(6-m)t+\frac{3m^2-45m+168}{4} t^2+\frac{-4m^3+105m^2-923m+2712}{8} t^3$\\
$\frac{5m^4-192m^3+2812m^2-18579m+46608}{16}t^4$ for $n\geq m\geq 9$,\\
$X_{34}=X_{35}=t$, $X_{36}=t-\frac{3}{2}t^2+\frac{3}{4}t^3-\frac{3}{8}t^4$, $X_{37}=t-3t^2+\frac{3}{2}t^3-\frac{3}{4}t^4$,\\
$X_{3,m-3}=t+\frac{3}{2}(8-m)t^2+\frac{6m^2-111m+516}{4} t^3+\frac{-10m^3+303m^2-3110m+10794}{8} t^4$\\
for $n\geq m \geq 11,$\\
$X_{45}=X_{46}=\frac{3}{2}t^2-\frac{3}{4}t^3+\frac{3}{8}t^4,$\\
$X_{4m-4}=\frac{3}{2}t^2+\frac{-12m+117}{4}t^3+\frac{30m^2-636m+3423}{8}t^4$, for $n\geq m\geq 11$\\
$X_{56}=X_{57}=3t^3-\frac{27}{4}t^4$\\
$X_{5m-5}=3t^3+\frac{-30m+333}{4}t^4$, for $n\geq m\geq 12$,\\
$X_{6m-6}=\frac{15}{2}t^4$ for $m\geq 13$ and $X_{ij}=0$, for $6<i<j$.
\end{enumerate}
\end{thm}
\begin{proof} To describe the scheme at neighbourhood of $\mathfrak{f}_{12}$, it suffices to write the canonical deformation expanding the structure constants yet obtained (up to the order $4$ and by using the equation $t^5=0$). For $n+1\geq 13$ we separately study the central extensions of $\mathfrak{f}_n$ and of $\mathfrak{w}_n$ proceeding by induction. Let $Y_{ij}$ denote the ancient structure constants for $i+j\leq n$ and $X_{ij}$ the new variables for $i+j= n+1$. The following relations
\begin{equation}
	\J_{1jn-j}=X_{j+1,n-j}-Y_{j,n-j}+X_{j,n+1-j}=0
\end{equation}
give
\begin{equation}\label{e6.56}
	X_{j,n+1-j}=\sum_{k=j}^{j'}(-1)^{k-j}Y_{k,n-k}+(-1)^{j'-j+1}X_{j'+1,n-j'}
\end{equation}
and it suffices to write another relation for fixing $X_{ij}$ in function of the $Y_{ij}$. More particularly the relation $\J_{24n-5}=0$ fixes $X_{6n-5}$.\\
For $\mathfrak{w}_n$, the obtained Lie algebra will be only the Witt Lie algebra $\mathfrak{w}_{n+1}$ and there will be not new parameters; the others Jacobi relations are then automatically satisfied.\\
For $\mathfrak{f}_n$, the relation (\ref{e6.56}) and $\J_{24n-5}=0$ fix the structure constants at neighbourhood of $\mathfrak{f}_{n+1}$ which verify the given relations in the theorem. Hence there are not new parameter. It remains to verify that the remained Jacobi relations do not reduce the order of nilpotence of the parameter $t$. It is obviously to see that the Jacobi polynomials $\J_{ijk}$ $(i<j<k)$ are vanished for $j\geq 6$ and it remains to see that $\J_{ijk}$ for $2\leq i\leq j\leq 5$ are null too.
\end{proof}
\subsubsection* {Acknowledgments}
\noindent The authors would like thank Carmen M$\mathrm{\acute{a}}$rquez Garcia at the university of Sevilla, Spain, which checked the calculations of examples in a compute.


\begin{thebibliography}{9}
\bibitem{A} M.~Artin, On the solutions of analytic equations, Inv. Math. 5, 277-291 (1968).
\bibitem {B} N.~Bourbaki, Alg\`ebre, chap. 4 \`a 7, Masson (1981).
\bibitem{Br} F.~Bratzlavsky, Sur les alg\`ebres admettant un tore d'automorphismes donn\'e. J. of Algebra, vol.30 (1974), pp. 305-316. 
\bibitem{BTM}M.~Bordemann, A.~Makhlouf, T.~Petit, D\'{e}formation par quantification et 
rigidit\'{e}  des alg\`{e}bres enveloppantes, J. Algebra,  Vol 285/2 , pp 623-648, 2005.
%\bibitem{Bn} S.~Benayadi, Certaines propri\'et\'es d'une classe d'alg\`ebres de Lie qui g\'en\'eralisent les alg\`ebres de Lie semi-simples. Ann. Fac. Sci. Toulouse, vol. XII, 1, 1991, pp. 29-53.
\bibitem{C1} R.~Carles, Sur la structure des alg\`ebres de Lie rigides, Ann. Inst. Fourier, Grenoble, vol. 34,3 (1984), pp. 65-82.
%\bibitem{C6} R.~Carles, Sur certaines classes d'alg\`ebres de Lie rigides, Annales Institut Fourier XXIV-3-1984
\bibitem{C} R.~Carles, Sur les alg\`ebres de Lie caract\'eristiquement nilpotentes, Pr\'epublication Poitiers no 5, (1984).
\bibitem{C2} R.~Carles, Sur certaines classes d'alg\`ebres de Lie rigides, Math ann, 272, 477-488, (1985).
\bibitem{C7} R.~Carles, Introduction aux d\'eformations d'alg\`ebres de Lie de dimension finie, Preprint Universit\'{e} de Poitiers, n°19, 1986.
\bibitem{C3} R.~Carles, D\'eformations et \'el\'ements nilpotents dans les sch\'emas d\'efinis par les identit\'es de Jacobi. CRAS,312,671, 1991.
\bibitem{C4} R.~Carles, Construction des alg\`ebres de Lie compl\`etes, CRAS, 318, 711-714 (1994).
\bibitem{C5} R.~Carles, Weight systems, Poitiers, Preprint Universit\'{e} de Poitiers, 96 (1996).
\bibitem{C8} R.~Carles, Sur la cohomologie d'une nouvelle classe d'alg\`ebres de Lie qui g\'en\'eralisent les sous-alg\`ebres de Borel, J of Algebra Vol 154, n 2, 1993.
\bibitem{C9} R.~Carles, D\'eformations dans les sch\'emas d\'efinis par les identit\'es de Jacobi. Ann. Math. Blaise Pascal 3 (1996), no. 2, 33--62.
\bibitem{CC}R.~Carles, C.M$\mathrm{\acute{a}}$rquez Garcia, Four different methods for the study of obstructions in the schemes of Jacobi, in preparation.
\bibitem{CP} R.~Carles, T.~Petit, On deformation in Poisson polynomial scheme, in preparation.
\bibitem{Di} J.~Dixmier and W.G Lister, Derivations of nipotent Lie algebras. Proc. Amer. Soc. 8 (1957), pp. 155-157. 
%\bibitem{D}J.~Dixmier, Cohomologie des alg\'ebres de Lie nilpotent, 
%\bibitem{Do} J.~Dozias, Sur les d\'erivations des alg\`ebres de Lie. C. R. Acad. Sci. Paris Sér. I. Math. 259, p. 2748-2750 ( 1964).
%\bibitem{E} D.~Eisenbud, Commutative Algebra with a View Toward Algebraic, 1996, Springer.
\bibitem{Fa} G.~Favre, Syst\`emes de poids sur une alg\`ebres de Lie nilpotente. Manuscripta. Math. 9 (1973), pp. 53-90.
\bibitem{F} A.~Fialowski, On the classification of graded Lie algebras with two generators. Vestnyik MGU, (1983), pp. 62-64 (in Russian). English translation: Moscow Univ. Math. Bull., 38 (1983), no. 2, 76 - 79 
%\bibitem{F2} A.~Fialowski, Post, G., "Versal deformation of the Lie algebra L2", Journal of Algebra, Jan. 2001, 236(1), pp. 93-109
%\bibitem{FF}Fialowski A., Fuchs D.B., Construction of Miniversal Deformations of Lie Algebras. 
%Journal of Functional Analysis, 161/1 (1999), pp. 76-110
\bibitem{G} M.~Gerstenhaber, On the deformations of rings and algebras II. Ann. of Math., vol. 79 (1964), 59-103.
\bibitem{HS} G.~Hochschild, J-P.~Serre, Cohomology of Lie algebras, Ann.~Math. 57 (1953), 72-144.
\bibitem{K} M.Kuranishi, On the locally complete families of complex analytic structures, Ann. of Math., vol. 75 (1962), 536-577.
\bibitem{L} E.M.~Luks, What is the Typical nilpotent Lie algebras, Computers in Nonassociative Rings and Algebras, Acad. Press, 1977.
\bibitem{Mo} G.D.~Mostow, Fully reductive subalgebras of algebraic groups, Amer. J. Math. 68, p. 220-306 (1956).
\bibitem{NR} A.~Nijenhuis and R.W.Richardson,Cohomology and deformations in graded Lie algebras, Bull. Amer. Math. Sco., vol 72 (1966), pp. 1-29.
\bibitem{P}T.Petit, On The Strong Rigidity of Solvable Lie algebras. Noncommutative Algebra and Geometry, Lecture Notes in Pure and Applied Mathematics, vol.243, p.162-174 Dekker, New York, 2006.
\bibitem{R} G.~Rauch, Effacement et d\'eformation. Ann.Ins.Fourier, vol. 22 (1972), 239-269.
%\bibitem{Sc} E.~Schenkman, A theory of subinvariant Lie algebras. Amer. J. Math. 73, p. 453-474 (1951).
\bibitem{S} M.~Schlessinger, Functors of Artin rings, Trans. Amer. Math. Sco. vol. 130 (1968), pp. 208-222.
\bibitem{V} M.~Vergne, Cohomologie des alg\`ebres de Lie nipotentes. Application \`a l'\'etude de la vari\'et\'e des alg\`ebres de Lie nipotentes. Bull. Soc. Math. France, 98 (1970), p. 81-116.
\end{thebibliography}
\end{document}